\documentclass{memo-l}

\usepackage[latin1]{inputenc}
\usepackage[T1]{fontenc}
\usepackage{amsfonts}
\usepackage{amsmath}
\usepackage{amssymb}
\usepackage{eurosym}
\usepackage{mathrsfs}
\usepackage{palatino}
\usepackage{color}
\usepackage{esint}
\usepackage{url}
\usepackage{verbatim}

\usepackage{enumerate}

\newcommand{\R}{\mathbb{R}}
\newcommand{\C}{\mathbb{C}}

\newcommand{\Z}{\mathbb{Z}}
\newcommand{\E}{\mathbb{E}}
\newcommand{\calP}{\mathcal{P}}
\newcommand{\calB}{\mathcal{B}}
\newcommand{\calG}{\mathcal{G}}
\newcommand{\calM}{\mathcal{M}}
\newcommand{\calS}{\mathcal{S}}
\newcommand{\scrS}{\mathscr{S}}

\newcommand{\calQ}{\mathcal{Q}}

\newcommand{\frakD}{\mathfrak{D}}

\newcommand{\scrD}{\mathscr{D}}
\newcommand{\bbP}{\mathbb{P}}

\newcommand{\diam}{\operatorname{diam}}
\newcommand{\scrA}{\mathscr{A}}
\newcommand{\calJ}{\mathcal{J}}
\newcommand{\Qi}{Q_{1,i}}
\newcommand{\Qj}{Q_{2,j}}
\newcommand{\Ri}{R_{1,i}}
\newcommand{\Rj}{R_{2,j}}
\newcommand{\bei}{\beta_{1,i}}
\newcommand{\bej}{\beta_{2,j}}
\newcommand{\phii}{\varphi_{1,i}}
\newcommand{\phij}{\varphi_{2,j}}
\newcommand{\wi}{w_{1,i}}
\newcommand{\wj}{w_{2,j}}
\newcommand{\bla}{\big \langle}
\newcommand{\bra}{\big \rangle}
\newcommand{\inte}[0]{\operatorname{int}}

\newcommand{\DI}{\calD_0(\omega_1)}
\newcommand{\DII}{\calD_0(\omega_2)}
\newcommand{\DIII}{\calD_0(\omega_3)}

\newcommand{\ud}[0]{\,\mathrm{d}}

\newcommand{\dist}[0]{\operatorname{dist}}

\newcommand{\BMO}[0]{\operatorname{BMO}}
\newcommand{\supp}[0]{\operatorname{spt}}
\newcommand{\loc}[0]{\operatorname{loc}}

\newcommand{\good}[0]{\operatorname{good}}
\newcommand{\bad}[0]{\operatorname{bad}}
\newcommand{\ch}[0]{\operatorname{ch}}

\newcommand{\calD}[0]{\mathcal{D}}
\newcommand{\rest}{{\lfloor}}
\newcommand{\wt}[1]{{\widetilde{#1}}}

\newtheorem{theorem}{Theorem}[chapter]
\newtheorem{lemma}[theorem]{Lemma}
\newtheorem{corollary}[theorem]{Corollary}
\newtheorem{proposition}[theorem]{Proposition}

\theoremstyle{definition}

\theoremstyle{remark}
\newtheorem{remark}[theorem]{Remark}

\numberwithin{section}{chapter}
\numberwithin{equation}{chapter}

\makeindex

\begin{document}

\frontmatter

\title{Dyadic-probabilistic methods in bilinear analysis}

\author{Henri Martikainen}
\address[H.M.]{Department of Mathematics and Statistics, University of Helsinki, P.O.B. 68, FI-00014 University of Helsinki, Finland}
\email{henri.martikainen@helsinki.fi}
\thanks{H.M. is supported by the Academy of Finland through the grants 294840, 306901 and 266262, and is a member of the Finnish Centre of Excellence in Analysis and Dynamics Research.}

\author{Emil Vuorinen}
\address[E.V.]{Department of Mathematics and Statistics, University of Helsinki, P.O.B. 68, FI-00014 University of Helsinki, Finland}
\email{emil.vuorinen@helsinki.fi}
\thanks{E.V. is partially supported by T. Hyt\"onen's ERC Starting Grant  Analytic-probabilistic methods for borderline singular integrals, and is a member of the Finnish Centre of Excellence in Analysis and Dynamics Research.}



\maketitle

\tableofcontents

\begin{abstract}
We demonstrate and develop dyadic--probabilistic methods in connection with non-homogeneous bilinear operators, namely singular integrals and square functions.
We develop the full non-homogeneous theory of bilinear singular integrals using a modern point of view. The main result is a new global $Tb$ theorem for
Calder\'on--Zygmund operators in this setting.
Our main tools include maximal truncations, adapted Cotlar type inequalities and suppression and big piece methods. 

While proving our bilinear results we also advance and refine the linear theory of Calder\'on--Zygmund operators by improving techniques and results. For example, we simplify and make more efficient some non-homogeneous summing arguments
appearing in $T1$ type proofs. As a byproduct, we can manage with ease quite general modulus of continuity in the kernel estimates. Our testing conditions are also quite general by virtue of the big piece method of proof.
\end{abstract}

\mainmatter
\chapter{Introduction}

The best known boundedness results for usual Calder\'on--Zygmund operators and square functions
are proved using dyadic analysis and probabilistic methods. The general philosophy was originally
introduced by Nazarov--Treil--Volberg to deal with non-homogeneous measures, see e.g. \cite{NTV}, \cite{NTVa}. Such techniques have then been widely used
and refined to multiple directions. See e.g. Azzam--Hofmann--Martell--Mayboroda--Mourgoglou--Tolsa--Volberg \cite{AHMMMTV} (rectifiability of the harmonic measure),
Hyt\"onen \cite{Hy} ($A_2$ theorem and the representation of singular integrals), Lacey--Sawyer--Uriarte-Tuero--Shen \cite{LSUS} and Lacey \cite{La1} (two weight inequality for the Hilbert transform), Lacey--Martikainen \cite{LM:CZO} (non-homogeneous local $Tb$ theorem with $L^2$ testing conditions), Martikainen \cite{Ma1} (representation of bi-parameter singular integrals), and Tolsa \cite{ToBook}, \cite{To1} (general account of non-homogeneous theory and the Painlev\'e's problem).

We study the dyadic martingale structure behind \textbf{bilinear} operators, mainly
bilinear Calder\'on--Zygmund and square function operators, in the setting of non-homogeneous analysis.
The main result is a new global $Tb$ theorem for Calder\'on--Zygmund operators in this non-homogeneous bilinear setup.
However, we also continue to push some of the most recent and advanced techniques further. Thus,
our proof techniques also yield some new insight about the basic linear theory, for example, by replacing certain previous methods
by the wider use of maximal truncations, various Cotlar type inequalities, and the big piece and suppression methods.
These four aspects really are the cornerstone to our approach, and we shall carefully explain their place in the proof when we encounter them. 

Firstly, we take a look at what dyadic--probabilistic methods can do, and how to use them, in the multilinear world. It appears that these methods have not really been used before
in this setting. This entails recording the general probabilistic martingale proof structure for bilinear operators. It also turns out that the treatment of non-homogeneous
singular integrals is quite involved in this setting with relevant technical challenges for example at the diagonal part.
Secondly, we exploit and develop certain very recent techniques. In particular, the big pieces method -- see e.g. Martikainen--Mourgoglou-Tolsa \cite{MMT1}
and Martikainen--Mourgoglou--Vuorinen \cite{MMV1}
 -- is extended to bilinear operators and is used to improve some integrability assumptions. In fact, we prove a global $Tb$ theorem with a new proof and
weaker testing conditions than usual (we can go below $L^1$). In particular, we don't use the RBMO$(\mu)$ space in our proof at all (as e.g. in \cite{NTV}).

We consider the following philosophy, which we also follow in this paper, to be the most important new realisations in the $Tb$ world: testing conditions involving
the maximal truncations $T_{\sharp}$, instead of the original Calder\'on--Zygmund operator $T$, are much easier to exploit (via suppression methods). This is the case even if the
testing conditions involving $T_{\sharp}$ are extremely weak. This idea seems to originate from the paper by Nazarov, Treil and Volberg \cite{NTV:Vit},
where they prove a special big piece type $Tb$ theorem for Cauchy integral type operators in connection with Vitushkin's conjecture.
The full potential of this approach was not immediately used in the $Tb$ circles, rather it has really started to become clear only recently -- 
see e.g. Hyt\"onen--Nazarov \cite{HN} and Martikainen--Mourgoglou-Tolsa \cite{MMT1}. One of the fundamental problems is that we only want to assume conditions involving $T$ itself,
and the passage to conditions involving $T_{\sharp}$ can be tricky. This requires some kind of adapted Cotlar's inequality i.e. Cotlar's inequality
which only uses the assumed testing conditions instead of some form of a priori boundedness. There are no such problems in the context of square functions, which
explains, in part, why square functions are so much simpler to handle.

Let us get back to the multilinear theory.
The literature on multilinear analysis is certainly vast. We mention just some closely related papers here.
Recent papers concerning multinear $T1$ or $Tb$ type theorems are e.g. Grafakos--Oliveira \cite{GO},
J. Hart \cite{Ha1, Ha2, Ha3} and Kova\v c--Thiele \cite{KT}.  Some formulations
related to multilinear \emph{local} $Tb$ theorems appear in Grau de la Herr\'an--Hart--Oliveira \cite{GHO} and Mirek--Thiele \cite{MT}, but
the multilinear local theory seems to have various restrictions which require further understanding.
This is one of our motivations also, but we refrain from touching that part of the theory too much in this paper. Indeed, things
already get quite technical.
\subsection*{Definitions and the main theorem}
We now formulate the setting and our main theorem. Of course, many of the needed results, the proofs, and the big picture of the proof are as interesting
as the main theorem. We will lay down the structure of the proof later.

A function
$$
K\colon (\R^n \times \R^n \times \R^n) \setminus \Delta \to \C, \qquad \Delta := \{(x,y,z) \in \R^n \times \R^n \times \R^n\colon\, x = y = z\},
$$
is called a standard bilinear $m$-dimensional Calder\'on--Zygmund kernel if for some $\alpha \in (0,1]$ and $C_K < \infty$ it holds that
$$
|K(x,y,z)| \le \frac{C_K}{(|x-y| + |x-z|)^{2m}},
$$ 
$$
|K(x,y,z)-K(x',y,z)| \le C_K\frac{|x-x'|^{\alpha}}{(|x-y| + |x-z|)^{2m+\alpha}}
$$
whenever $|x-x'| \le \max(|x-y|, |x-z|)/2$,
$$
|K(x,y,z)-K(x,y',z)| \le C_K\frac{|y-y'|^{\alpha}}{(|x-y| + |x-z|)^{2m+\alpha}}
$$
whenever $|y-y'| \le \max(|x-y|, |x-z|)/2$, and
$$
|K(x,y,z)-K(x,y,z')| \le C_K\frac{|z-z'|^{\alpha}}{(|x-y| + |x-z|)^{2m+\alpha}}
$$
whenever $|z-z'| \le \max(|x-y|, |x-z|)/2$.

Given two Radon measures $\nu_1, \nu_2$ on $\R^n$, possibly complex, we define, whenever the right hand side makes sense, that
$$
T_{\varepsilon}(\nu_1, \nu_2)(x) = \iint_{\max(|x-y|, |x-z|) > \varepsilon} K(x,y,z)\,d\nu_1(y)\,d\nu_2(z), \qquad x \in \R^n,\, \varepsilon > 0.
$$
The defining integral is absolutely convergent e.g. if $|\nu_i|(\R^n) < \infty$ for $i=1,2$. The truncations could
also be defined as
$$
\tilde T_{\varepsilon}(\nu_1, \nu_2)(x) = \iint_{|x-y|^2 + |x-z|^2 > \varepsilon^2} K(x,y,z)\,d\nu_1(y)\,d\nu_2(z), \qquad x \in \R^n,\, \varepsilon > 0.
$$
We prefer $T_{\varepsilon}$ over $\tilde T_{\varepsilon}$ as it seems somewhat easier to work with in connection with some pointwise estimates.
However, our main theorem, Theorem \ref{thm:main}, can also be stated using $\tilde T_{\varepsilon}$. In fact, such a version follows from the one with $T_{\varepsilon}$
using that
$$
|T_{\varepsilon}(\nu_1, \nu_2)(x) - \tilde T_{\varepsilon}(\nu_1, \nu_2)(x)| \lesssim M_m \nu_1(x) M_m\nu_2(x),
$$
where
$$
M_m\nu(x) = \sup_{r > 0} \frac{|\nu|(B(x,r))}{r^m}, \qquad x \in \R^n.
$$

For us a bilinear $m$-dimensional SIO (singular integral operator) $T$ is simply the collection $(T_{\varepsilon})_{\varepsilon > 0}$ in the sense
that we are only interested in some uniform in $\varepsilon > 0$ boundedness properties of the operators $T_{\varepsilon}$. This is all simply determined by the given kernel $K$.

We continue to define the maximal truncations as follows:
\begin{align*}
T_{\sharp, \delta}(\nu_1, \nu_2)(x) &= \sup_{\varepsilon > \delta} |T_{\varepsilon}(\nu_1,\nu_2)(x)|, \qquad \delta \ge 0; \\
T_{\sharp}(\nu_1,\nu_2)(x) &=  T_{\sharp, 0}(\nu_1,\nu_2)(x).
\end{align*}

A positive Radon measure $\mu$ on $\R^n$ is said to be of order $m$ if for some constant $C_{\mu} < \infty$ we have $$\mu(B(x,r)) \le C_{\mu}r^m$$
for all $x \in \R^n$ and $r > 0$. We set
\begin{align*}
T_{\mu, \varepsilon}(f,g)(x) &= T_{\varepsilon}(f\,d\mu, g\,d\mu)(x) \\
&= \iint_{\max(|x-y|, |x-z|) > \varepsilon} K(x,y,z) f(y)g(z)\,d\mu(y)\,d\mu(z).
\end{align*}
The above is well-defined as an absolutely convergent integral if e.g. $f \in L^{p_1}(\mu)$ and $g \in L^{p_2}(\mu)$ for some $p_1, p_2 \in [1, \infty)$, since
then
\begin{equation}\label{eq:basicbound}
\begin{split}
\iint_{\max(|x-y|, |x-z|) > \varepsilon} |K(x,y,z)& f(y)g(z)|\,d\mu(y)\,d\mu(z) \\ &\lesssim \frac{1}{\varepsilon^{m(1/p_1+1/p_2)}} \|f\|_{L^{p_1}(\mu)} \|g\|_{L^{p_2}(\mu)}. 
\end{split}
\end{equation}
We also set
\begin{align*}
T_{\mu, \sharp, \delta}(f,g)(x) &= \sup_{\varepsilon > \delta} |T_{\mu, \varepsilon}(f,g)(x)|, \qquad \delta \ge 0; \\
T_{\mu, \sharp}(f,g)(x) &= T_{\mu, \sharp, 0}(f,g)(x).
\end{align*}
The notation $T^{1*}$ and $T^{2*}$ stand for the adjoints of a bilinear operator $T$, i.e.
$$
\langle T(f, g), h\rangle = \langle T^{1*}(h, g), f\rangle = \langle T^{2*}(f, h),g\rangle.
$$

To state the main theorem we still need the concept of cubes with small boundary. See also the end of the introduction for additional notation, which
is rather standard.
Given $t > 0$ we say that a cube $Q \subset \R^n$ has $t$-small boundary with respect to the measure $\mu$ if
\begin{displaymath}
\mu(\{x \in 2Q\colon\, \dist(x,\partial Q) \le \lambda\ell(Q)\}) \le t\lambda\mu(2Q)
\end{displaymath}
for every $\lambda > 0$.
The definition of the suppressed operators $T_{\Phi}$ can be found in Chapter \ref{sec:suppopp}. They appear in the formulation of the theorem
in connection with a weak boundedness property -- at this point one should simply understand that it is a purely diagonal condition,
a necessary condition (as we will show), and is automatically satisfied should $K$ possess some antisymmetry.
See also Corollary \ref{cor:MainThm} below.

\begin{theorem}\label{thm:main}
Let $\mu$ be a measure of order $m$ on $\R^n$ and $T$ be a bilinear $m$-dimensional SIO. Let $t$ be a large enough dimensional constant, $s, c_b > 0$
and $C_b, C_W, C_{\textup{test}} < \infty$.
Let the functions $b_i$, $i = 1,2,3$, be such that $\|b_i\|_{L^{\infty}(\mu)} \le C_b$ and
$$
|\langle b_i \rangle_Q^{\mu}| \ge c_b \qquad \textup{for all cubes } Q.
$$
We assume the weak boundedness property in the form that
\begin{equation}\label{eq:WeakBdd}
\sup_{\delta > 0} |\langle T_{\mu, \Phi, \delta}(1_Qb_1, 1_Qb_2), 1_Qb_3\rangle_{\mu}| \le C_W\mu(5Q) \qquad \textup{for all cubes } Q
\end{equation}
whenever $\Phi\colon\, \R^n \to [0, \infty)$ is a $1$-Lipschitz function. Suppose that for every cube $Q$ with $t$-small boundary we have
the following three testing conditions
\begin{equation}\label{eq:ThmMainTesting}
\begin{split}
\sup_{\delta > 0} \sup_{\lambda > 0} \lambda^s \mu(\{x \in Q\colon\, |T_{\mu, \delta}(b_11_Q, b_21_Q)(x)| > \lambda\}) \le C_{\textup{test}} \mu(2Q), \\
\sup_{\delta > 0} \sup_{\lambda > 0} \lambda^s \mu(\{x \in Q\colon\, |T^{1*}_{\mu, \delta}(b_31_Q, b_21_Q)(x)| > \lambda\}) \le C_{\textup{test}} \mu(2Q), \\
\sup_{\delta > 0} \sup_{\lambda > 0} \lambda^s \mu(\{x \in Q\colon\, |T^{2*}_{\mu, \delta}(b_11_Q, b_31_Q)(x)| > \lambda\}) \le C_{\textup{test}} \mu(2Q).
\end{split}
\end{equation}
Then for all $1 <p, q < \infty$ and $1/2 < r < \infty$ satisfying
$1/p + 1/q = 1/r$ we have that
$$
\|T_{\mu, \sharp}\|_{L^p(\mu) \times L^q(\mu) \to L^r(\mu)} \lesssim 1
$$
with a constant depending on $p,q,r$, the above fixed constants, and the constants appearing in the definitions of $T$ and $\mu$.
\end{theorem}

The weak boundedness property is only needed in the regime $s <1$. We record the following corollary regarding this.

\begin{corollary}\label{cor:MainThm}
If in Theorem \ref{thm:main}  in place of the weak type testing conditions \eqref{eq:ThmMainTesting}
one assumes the strong type testing conditions with exponent $s=1$, that is, the conditions
\begin{equation*}
\begin{split}
\sup_{\delta > 0} \int_Q |T_{\mu, \delta}(b_11_Q, b_21_Q)|d \mu \le C_{\textup{test}} \mu(2Q), \\
\sup_{\delta > 0} \int_Q |T^{1*}_{\mu, \delta}(b_31_Q, b_21_Q)|d \mu \le C_{\textup{test}} \mu(2Q), \\
\sup_{\delta > 0} \int_Q |T^{2*}_{\mu, \delta}(b_11_Q, b_31_Q)|d \mu \le C_{\textup{test}} \mu(2Q),
\end{split}
\end{equation*}
then the weak boundedness assumption \eqref{eq:WeakBdd} can be dropped.
\end{corollary}

\begin{remark}
We are unable to see any obvious obstructions to proving everything also in the $k$-linear, $k > 2$, context. Usually this is more or less an exercise, although the
notation can get really complicated. Here, due to the quite complicated proof involving suppression and surgery arguments, we refrain from explicitly claiming this
having not checked everything in detail in the $k$-linear setting.
\end{remark}

\subsection*{The structure of the proof}
Corollary \ref{cor:MainThm} will be proved in Chapter \ref{sec:syn}.
The proof of Theorem \ref{thm:main} does not involve any type of bilinear interpolation. Rather, it has the following steps:
\begin{enumerate}
\item \emph{Prove testing conditions for $T_{\sharp}$.} This is Corollary \ref{cor:cot} and it requires the adapted Cotlar's inequality Proposition \ref{prop:cot}.
This step is motivated by the techniques used in the very recent proof of a certain (linear) local $Tb$ theorem by Martikainen--Mourgoglou-Tolsa \cite{MMT1}.
See also Hyt\"onen--Nazarov \cite{HN} for the first instance of such an adapted Cotlar's inequality in the Lebesgue case.

\item \emph{Introduce suppressed operators $T_{\Phi}$.} This is done in Chapter \ref{sec:suppopp}.
These originate (in some special case) from \cite{NTV:Vit}, and have since then been used quite a lot.
The idea is that if $\Phi$ is chosen suitably, then $T_{\Phi}$ behaves better than $T$, but also
agrees with $T$ on the set $\{\Phi = 0\}$ -- which is arranged to be relatively large. This is called suppression.
These operators are needed in the next step -- the big piece $Tb$.
The bilinear suppression details do not differ too much from the linear ones, but we believe our presentation should be logical and nicely readable.

The following technical thing needs to be noted. The suppressed operators already make an appearance in the statement of the main theorem, because we use a
somewhat non-classical formulation of the weak boundedness property. However,
it is still a purely diagonal condition, a necessary condition (as we will show), and is automatically satisfied should $K$ possess some antisymmetry.
We seem to need this since our proof strategy goes through a new formulation of the big piece $Tb$.
The upshot is that we can allow $s < 1$ in the main testing conditions.

\item \emph{Prove a version of the big piece $Tb$ that can be applied to prove the main theorem.} This is Theorem \ref{thm:bigpieceTb}.
The formulation is necessarily relatively technical, and cannot be fully described here. Briefly, the testing conditions for $T_{\sharp}$ (proved in step 1) allow us to do the suppression
from Step 2. The big piece $Tb$ allows us to conclude that the operators $T_{\varepsilon}$ are, uniformly in $\varepsilon > 0$, bounded on a big piece of a given nice cube $Q$ (of small boundary).

The proof of this contains the bilinear dyadic--probabilistic $Tb$ argument, and so forms the technical core of the paper. There are many details
here to be noted. The bilinear framework seems to complicate at least the treatment of the non-homogeneous paraproduct and the diagonal.
We also simplify quite a few details from the linear theory with some new summation arguments -- we, for example, can make do without some standard
matrix summation lemmas previously extensively used in these arguments.

Moreover, the above noted improvements in the summation arguments allow us to replace, with almost the same proof, the modulus of continuity $t \mapsto t^{\alpha}$
by a modulus of continuity $\psi$ satisfying the modified Dini--type condition
$$
\int_0^1 \psi(t) \Big( 1 + \log \frac{1}{t} \Big) \frac{dt}{t} < \infty.
$$
That is, we almost recover the best known modulus of continuity
$$
\int_0^1 \psi(t) \Big( 1 + \log \frac{1}{t} \Big)^{1/2} \frac{dt}{t} < \infty
$$
for free (most of the other theory works with weaker assumptions but the $T1$ does not). We give these modifications in Chapter \ref{sec:dini}.
For a deeper understanding of these issues see the paper by Grau de la Herr\'an and Hyt\"onen \cite{GH} of which our paper is completely independent of.
\item \emph{Prove weak type end point estimates for $T_{\sharp}$.} This is done in Chapter \ref{sec:endpoint}, and is laborious and technical in our generality.
We chose to give the full details. This part is separate from the rest of the steps as it is essentially basic theory of non-homogeneous bilinear Calder\'on--Zygmund operators -- but we need
to write it here since we are unaware of references operating in our generality. A more standard formulation (than in step 1) of Cotlar's inequality is also needed here.
The conclusions need to involve $T_{\sharp}$, since the good lambda method (step 6) requires it.

\item \emph{Prove a bilinear adaptation of the good lambda method of Tolsa}. This is Theorem \ref{thm:goodlambda}. It is to be noted that the bilinear version
is quite straightforward to prove mimicking the linear proof. The good lambda method is a glue which yields global boundedness from local big piece type boundedness.
It is an extremely flexible tool due to the fact that one has a lot of freedom in the statement -- for example, the good lambda only requires very nice cubes (doubling and of small boundary).

\item \emph{Synthesis.} In Chapter \ref{sec:syn} we give the proof of the main theorem, Theorem \ref{thm:main}. It is a very short argument using the steps 1-6. It is probably instructive
to take a look at this proof to get an idea of the big picture before looking at all the details. In Chapter \ref{sec:syn} we also prove Corollary \ref{cor:MainThm}.
\end{enumerate}
In addition to the above, we demonstrate and make very heavy use of various dyadic $L^p$, $p \ne 2$, techniques. This is because we choose to prove
our big piece $Tb$, Theorem \ref{thm:bigpieceTb}, directly with general exponents $1 < p, q, r < \infty$. For the proof of our main theorem, we only need
this big piece $Tb$ with $r=2$ and $p=q=4$. It is possible that the proof could be somewhat simpler in this case, but the space $L^4$ would appear anyway.
We prefer the general $L^p$ techniques.

Let us mention that we briefly discuss the much simpler case of square functions in Chapter \ref{sec:SF}.

\subsection*{More about the contributions to the linear theory}
We discuss in more detail a few aspects of the proof that also have relavance to the linear theory.

One aspect is that in the non-homogeneous setting the natural spaces $\BMO_p(\mu)$, $1 \le p < \infty$, where
$$
\BMO_p(\mu) := \Big\{b \in L^1_{\loc} \colon \sup_Q \Big( \frac{1}{\mu(2Q)} \int_Q |b - \langle b \rangle^{\mu}_Q|^p \ud \mu\Big)^{1/p}  < \infty\Big\},
$$
need not be equivalent. 
If one wants to prove a global non-homogeneous $Tb$ theorem in the linear setting with the assumption
$Tb \in \BMO_1(\mu)$, one runs into the problem that in the proof one would need $Tb \in \BMO_2(\mu)$ -- see Nazarov--Treil--Volberg \cite{NTV}.
(To understand the discussion notice that $\sup_Q \frac{1}{\mu(Q)} \int_Q |T(b1_Q)| \ud \mu < \infty$ implies that $Tb \in \BMO_1(\mu)$ so that local testing conditions, such as those that appear in this paper, are related to global BMO testing discussed here.)

Although it can be that $$\BMO_2(\mu) \subsetneq \BMO_1(\mu),$$ it turns out that under the assumptions of the $Tb$ theorem
one can prove that $Tb \in \BMO_1(\mu)$ implies $Tb \in \BMO_2(\mu)$. This exploits the fact that the function under question
is of the special form $Tb$, and requires a relatively complicated argument which needs to go through the regularised BMO space
$\operatorname{RBMO}(\mu)$ of Tolsa. Indeed, it is first proved that $Tb \in \operatorname{RBMO}(\mu)$, and then it is used that $\operatorname{RBMO}(\mu)$ is a space that enjoys the John--Nirenberg property.

Our proof can handle testing conditions like $$\sup_Q \frac{1}{\mu(2Q)} \int_Q |T(b1_Q)| \ud \mu < \infty,$$ even with exponents below $1$ and the supremum
running only over cubes with small boundary, but still no such RBMO arguments like above appear. This is because of the big piece type proof method discussed already above.
The big piece $Tb$ arguments, that originate from the paper by Nazarov--Treil--Volberg \cite{NTV:Vit}, are in the spirit that 
under substantially weaker assumptions than usual, one achieves the boundedness of the operator on a big piece of the reference space -- but not globally.
In deep results on the intersection of harmonic analysis and geometric measure theory, it has been important to use such non-standard $Tb$ arguments --
see Tolsa's book \cite{ToBook}. However, here (and previously in \cite{MMV1, MMT1}) we are stressing the point that
these big piece arguments can actually be pushed to yield the full boundedness of the operator, which is the usual conclusion that one expects of a $Tb$ theorem. That is, using the proof method described above we can reduce the proofs of $Tb$ theorems to certain big piece type $Tb$ theorems. How difficult a big piece $Tb$ theorem is needed depends on the context. We stress that to prove a global $Tb$ theorem, like we do in this paper, a fully general big piece $Tb$ theorem (such as the extremely complicated Theorem 5.1 in \cite{ToBook}) is not needed -- this is key for keeping this approach reasonable.

The magic is that in the proof of a big piece $Tb$ theorem no BMO arguments appear as the related
suppression arguments essentially give that $Tb \in L^{\infty}$. Where one needs to work harder than usual, however, is in proving the related
adapted Cotlar's inequality that is needed in the reductions discussed above.

Another aspect is the summing arguments appearing in the core part of the $Tb$ proof. Experts will recognise
that usual non-homogeneous proofs, such as the original one \cite{NTV} and essentially all the others after that, have employed e.g. the fact that the numbers
$$
\delta(Q,R) := \frac{\ell(Q)^{\alpha/2} \ell(R)^{\alpha/2}}{D(Q, R)^{m+\alpha}}, \qquad D(Q,R) := \ell(Q) + \ell(R) + d(Q,R),
$$
where $Q,R$ are dyadic cubes, satisfy $\ell^2$ estimates like
$$
\sum_{Q, R} \delta(Q,R) x_Q y_Q \lesssim \Big( \sum_Q x_Q^2 \Big)^{1/2} \Big( \sum_Q y_Q^2 \Big)^{1/2}.
$$
We are no longer relying on such estimates. Instead, we are relying on easier summing arguments -- for these see the proof of Theorem \ref{thm:bigpieceTb} (see
e.g. the easy argument for the separated sum).
This is convenient and interesting in and on itself, but is also key to getting an easy access to the more general modulus of continuity
$$
\int_0^1 \psi(t) \Big( 1 + \log \frac{1}{t} \Big) \frac{dt}{t} < \infty.
$$

\subsection*{Additional notation}
We write $$A \lesssim B,$$ if there is a constant $C>0$ (depending only on some fixed constants like $m,n, \alpha$ etc.) so that $A \leq C B$.
Moreover, $$A \lesssim_{\tau} B$$ means that the constant $C$ can also depend on some relevant given parameter $\tau > 0$.
We may also write $$A \sim B$$ if $B \lesssim A \lesssim B$.

We then define some notation related to cubes. If $Q$ and $R$ are two cubes we set:
\begin{itemize} 
\item $\ell(Q)$ is the side-length of $Q$;
\item If $a>0$, we denote by $aQ$ the cube that is concentric with $Q$ and has sidelength $a\ell(Q)$; 
\item The cube $Q$ is called $(\alpha,\beta)$-doubling for a given Radon measure $\mu$ if $$\mu(\alpha Q) \le \beta\mu(Q);$$
\item $d(Q,R) = \dist(Q,R)$ denotes the distance between the cubes $Q$ and $R$;
\item $D(Q,R):=d(Q,R)+\ell(Q)+\ell(R)$ is the long distance;
\item $\text{ch}(Q)$ denotes the dyadic children of $Q$;
\item $\mu \rest Q$ denotes the measure $\mu$ restricted to $Q$;
\item If $Q$ is in a dyadic grid, then $Q^{(k)}$ denotes the unique dyadic cube $S$ in the same grid so that $Q \subset S$ and $\ell(S) = 2^k\ell(Q)$;
\item If $\calD$ is a dyadic grid, then $$\calD_k = \{Q \in \calD\colon\, \ell(Q) = 2^{-k}\};$$
\item $\langle f \rangle_Q^{\mu} = \mu(Q)^{-1}\int_Q f\,d\mu = \langle f \rangle_Q$ (when the measure is clear from the context).
One can interpret this to equal zero if $\mu(Q) = 0$.
\end{itemize}
The notation $\langle f, g\rangle_{\mu}$ stands for the pairing $\int fg\,d\mu$.

The following maximal functions are also used:
\begin{align*}
M_{\mu, \calD} f(x) &= \sup_{Q \in \calD} 1_Q(x) \langle |f| \rangle_Q^{\mu} \qquad (\calD \textup{ is a dyadic grid}); \\
M_{\mu, m} f(x) &= \sup_{r>0} \frac{1}{r^m} \int_{B(x,r)} |f|\,d\mu; \\
M_{\mu} f(x) &= \sup_{r>0} \, \langle |f| \rangle_{B(x,r)}^{\mu}; \\ 
M_{\mu}^{\mathcal Q} f(x) &= \sup_{r>0} \, \langle |f| \rangle_{Q(x,r)}^{\mu}; \\
N_{\mu} f(x) &= \sup\Big\{ \frac{1}{\mu(5B)} \int_{B} |f|\,d\mu\colon\, B \textup{ is a ball containing the point } x\Big\}.
\end{align*}
Here $Q(x,r)$ stands for the open cube with center $x$ and side length $2r$, while $$B(x,r) = \{y\colon\, |x-y| < r\}.$$
Given $s > 0$ we define the $s$-adapted maximal functions as in $$M_{\mu,s} f(x) = M_{\mu} ( |f|^s )(x)^{1/s}.$$ The bilinear
variants are defined in the natural way, e.g.
$$
M_{\mu}(f,g)(x) = \sup_{r>0} \, \langle |f| \rangle_{B(x,r)}^{\mu} \langle |g| \rangle_{B(x,r)}^{\mu}.
$$

We can also hit complex measures $\nu$ with these maximal functions -- simply replace the appearing integrals $\int_A |f|\,d\mu$
with $|\nu|(A)$, e.g. $$M_{\mu} \nu(x) =  \sup_{r>0}  \frac{|\nu|(B(x,r))}{\mu(B(x,r))}.$$ 

The following additional notation for singular integrals is occasionally useful.
We want to sometimes be able to e.g. hit $f \otimes g(y,z):= f(y)g(z)$ instead of the pair $(f,g)$ -- to enable this we use the notation $\tilde T$ as below.
For every $\varepsilon>0$ and measure $\sigma$ in $\R^{2n}$ we formally define $\tilde{T}_{\varepsilon} \sigma$ by setting  
\begin{equation*}
\tilde{T}_{\varepsilon} \sigma(x)
:=  \iint \displaylimits_{\max(|x-y|,|x-z|)>\varepsilon} K(x,y,z)\, d\sigma(y,z), \qquad x \in \R^n.
\end{equation*}
All the other notions involving $\tilde T$ are defined analogously. 

Lastly, we record here the following standard estimate that we shall have frequent use for.
\begin{lemma}\label{lem:basic}
Let $x \in \R^n$ and $t > 0$. The following estimate holds
$$
\iint \frac{d|\nu_1|(y)\,d|\nu_2|(z)}{(t+|x-y|+|x-z|)^{2m+\alpha}} \lesssim t^{-\alpha} M_{m}(\nu_1,\nu_2)(x)
$$
for all appropriate complex measures $\nu_1$ and $\nu_2$.
\end{lemma}
\begin{proof}
Simply split the domain of integration to $$\max(|x-y|, |x-z|) < t$$ and $$2^{k-1}t \le \max(|x-y|, |x-z|) < 2^kt, \qquad k \ge 1,$$ and estimate
in a straightforward way.
\end{proof}

\subsection*{Acknowledgements} 
We thank Professor Hyt\"onen for suggesting to us that our new summation method in the $Tb$ argument should require only the modified Dini-type regularity of the kernel as explained in Chapter 8. We also thank the anonymous referee for his or her comments.

\chapter{Adapted Cotlar type inequality and testing condition for $T_{\mu, \sharp}$}

The purpose of this chapter is to show that uniform testing conditions concerning $T_{\mu, \delta}(1_Qb_1, 1_Qb_2)$, $\delta > 0$,
imply testing conditions for $T_{\mu, \sharp}(1_Qb_1, 1_Qb_2)$. This is important
when we want to apply the big pieces type $Tb$ theorem. We achieve this improved testing via the following version
of Cotlar's inequality. It is extremely important to note that this version of Cotlar only uses the assumed testing conditions -- not some form of a priori boundedness.
\begin{proposition}\label{prop:cot}
Let $\mu$ be a measure of degree $m$ on $\R^n$, $T$ be a bilinear $m$-dimensional SIO and $s, \delta > 0$.
Let $b$ and $t$ be large enough constants depending only on the dimension $n$. Suppose
$b_i \in L^{\infty}(\mu)$, $i = 1, 2$, satisfy for every $(2,b)$-doubling cube $R$ with $t$-small boundary that
$$
\sup_{\lambda >0} \lambda^s \mu(\{x \in R\colon\, |T_{\mu, \delta}(1_Rb_1, 1_Rb_2)(x)| > \lambda\}) \lesssim \mu(R).
$$
Suppose $Q \subset \R^n$ is a fixed cube and $\tau > 0$. Then uniformly for every $\varepsilon > \delta$ and $x \in (1-\tau)Q$ there holds that
$$
|T_{\mu, \varepsilon}(1_Qb_1, 1_Qb_2)(x)| \lesssim C(\tau) +  M_{\mu, s/4}^{\mathcal{Q}}(1_QT_{\mu, \delta}(1_Qb_1, 1_Qb_2))(x).
$$
\end{proposition}
\begin{proof}
Fix $x \in (1-\tau)Q$  and $\varepsilon_0 > \delta$.
Let $C(n)$ be a large dimensional constant.
In what follows we will implicitly need that $b$ is sufficiently much larger than $C(n)$, say $b > C(n)^{n+1} \ge C(n)^{m+1}$.
Choose the smallest $k$ so that $B(x, C(n)^k\varepsilon_0)$ is $((C(n), b)$-doubling. Set $\varepsilon = C(n)^k\varepsilon_0$. 
Notice that
\begin{align*}
|T_{\mu, \varepsilon_0}(1_Qb_1, 1_Qb_2)(x) - T_{\mu, \varepsilon}(1_Qb_1, 1_Qb_2)(x)|  
\lesssim 1.
\end{align*}
The last estimate is seen using a standard calculation based on the choice of $\varepsilon$ (see Lemma 2.15 in \cite{ToBook}). This calculation is performed carefully
in somewhat more generality at the beginning of the proof of Proposition \ref{prop:cotbas}.
Therefore, it suffices to estimate $|T_{\mu, \varepsilon}(1_Qb_1, 1_Qb_2)(x)|$.

If epsilon happens to be large enough compared to $\ell(Q)$, this is easy. Indeed, for $\varepsilon > c_{\tau}\ell(Q)$, say, we have
$$
|T_{\mu, \varepsilon}(1_Qb_1, 1_Qb_2)(x)| \lesssim \mu(Q) \int \frac{d\mu(z)}{(\varepsilon + |x-z|)^{2m}} 
\lesssim_{\tau} \frac{\mu(Q)}{\ell(Q)^m} \lesssim 1. 
$$
We can therefore assume that $\varepsilon \le c_{\tau}\ell(Q)$ for a sufficiently small constant $c_{\tau} > 0$.
Choose (using Lemma 9.43 of \cite{ToBook}) a cube $R$ centred at $x$ so that it has $t$-small boundary with respect to $\mu$, and
\begin{displaymath}
B(x, \varepsilon) \subset R \subset B(x, C_n\varepsilon) \subset Q.
\end{displaymath}
The last inclusion holds if $c_{\tau}$ is fixed small enough. Notice that $R$ is $(2,b)$-doubling as
$$
\mu(2R) \le \mu(B(x, C(n)\varepsilon)) \le b\mu(B(x, \varepsilon)) \le b\mu(R),
$$
where we used that $C(n)$ was chosen large in the beginning. In particular, we have for some $C_0 < \infty$ that
\begin{equation}\label{eq:cf1} 
\mu(\{w \in R\colon\, |T_{\mu, \delta}(1_Rb_1, 1_Rb_2)(w)| > \lambda\}) \le C_0\frac{\mu(R)}{\lambda^s}
\end{equation}
for every $\lambda > 0$.

Write for fixed $w \in R$ the equality
\begin{align*}
T_{\mu, \varepsilon}(1_Qb_1, 1_Qb_2)(x) =  T_{\mu, \varepsilon}&(1_Qb_1, 1_Qb_2)(x) - T_{\mu, \delta}(1_{(2R)^c}1_Qb_1, 1_Qb_2)(w) \\
&+ T_{\mu, \delta}(1_Qb_1, 1_Qb_2)(w) - T_{\mu, \delta}(1_{2R}1_Qb_1, 1_Qb_2)(w).
\end{align*}
Notice that $(2R)^c \subset B(w, \delta)^c \cap B(x,\varepsilon)^c$ so that
$$
|T_{\mu, \varepsilon}(1_Qb_1, 1_Qb_2)(x) - T_{\mu, \delta}(1_{(2R)^c}1_Qb_1, 1_Qb_2)(w)|
$$
can be dominated by the sum of
$$
\int_{(2R)^c} \int_{\R^n} |K(x,y,z) - K(w,y,z)|\,d\mu(z)\,d\mu(y)  \lesssim \ell(R)^{\alpha} \int_{R^c} \frac{d\mu(y)}{|x-y|^{m+\alpha}} \lesssim 1
$$
and
$$
\mathop{\iint_{\max(|x-y|, |x-z|) > \varepsilon}}_{y \in 2R}  |K(x,y,z)|\,d\mu(z)\,d\mu(y)
\lesssim \int_{2R} \int \frac{d\mu(z)\,d\mu(y)}{(\varepsilon + |x-z|)^{2m}}
\lesssim \frac{\mu(2R)}{\varepsilon^m} \lesssim 1.
$$
Therefore, we have
\begin{align*}
|T_{\mu, \varepsilon}(1_Qb_1, 1_Qb_2)(x)| \lesssim 1 + |T_{\mu, \delta}(1_Qb_1, 1_Qb_2)(w)| + |T_{\mu, \delta}(1_{2R}1_Qb_1, 1_Qb_2)(w)|.
\end{align*}
It follows from this by raising to the power $s/4$, averaging over $w \in R$ and raising to power $4/s$ that
\begin{align*}
|T_{\mu, \varepsilon}(1_Qb_1, 1_Qb_2)(x)| \lesssim 1 &+ M_{\mu, s/4}^{\mathcal{Q}}(1_QT_{\mu, \delta}(1_Qb_1, 1_Qb_2))(x) \\ 
&+ \bigg(\frac{1}{\mu(R)} \int_R |T_{\mu, \delta}(1_{2R}1_Qb_1, 1_Qb_2)(w)|^{s/4}\,d\mu(w)\bigg)^{4/s}.
\end{align*}
To get the maximal function bound, we also used that $R$ is a cube centred at $x$ and $R \subset Q$. 

We have that
$$
 \int_R |T_{\mu, \delta}(1_{2R \setminus R}1_Qb_1, 1_Qb_2)(w)|^{s/4}\,d\mu(w) \lesssim \int_R \bigg[ \int_{2R \setminus R} \frac{d\mu(y)}{|w-y|^m}\bigg]^{s/4}\,d\mu(w)
 \lesssim \mu(R)
$$
using that $R$ has small boundary and is doubling (see Lemma 9.44 in \cite{ToBook}).

 Next, we bound
\begin{align*}
\int_R |T_{\mu, \delta}(1_{R}b_1, 1_Qb_2)(w)|^{s/4}\,d\mu(w) \lesssim \int_R &|T_{\mu, \delta}(1_{R}b_1, 1_Rb_2)(w)|^{s/4}\,d\mu(w) \\
&+ \int_R |T_{\mu, \delta}(1_{R}b_1, 1_{2R \setminus R}1_Qb_2)(w)|^{s/4}\,d\mu(w) \\
&+ \int_R |T_{\mu, \delta}(1_{R}b_1, 1_{(2R)^c}1_Qb_2)(w)|^{s/4}\,d\mu(w).
\end{align*}
Notice that
$$
 \int_R |T_{\mu, \delta}(1_{R}b_1, 1_{2R \setminus R}1_Qb_2)(w)|^{s/4}\,d\mu(w) \lesssim \int_R \bigg[ \int_{2R \setminus R} \frac{d\mu(z)}{|w-z|^m}\bigg]^{s/4}\,d\mu(w) \lesssim \mu(R)
$$
and
\begin{align*}
\int_R |T_{\mu, \delta}(1_{R}b_1, 1_{(2R)^c}1_Qb_2)(w)|^{s/4}\,d\mu(w) &\lesssim \int_R \bigg[ \int_{R} \int_{(2R)^c} \frac{d\mu(z)}{|z-w|^{2m}} \,d\mu(y)\bigg]^{s/4}\,d\mu(w) \\
&\lesssim  \mu(R) \bigg[ \mu(R) \int_{R^c} \frac{d\mu(z)}{|z-x|^{2m}} \bigg]^{s/4} \\
&\lesssim \mu(R) \bigg[ \frac{\mu(R)}{\ell(R)^m} \bigg]^{s/4} \lesssim \mu(R).
\end{align*}

It only remains to show that $I \lesssim \mu(R)$ for the term
$$
I := \int_R |T_{\mu, \delta}(1_{R}b_1, 1_Rb_2)(w)|^{s/4}\,d\mu(w).
$$
The point simply is that weak type testing implies strong type testing for strictly smaller exponents.
Indeed, using \eqref{eq:cf1} we see that
\begin{align*}
I &= \frac{s}{4} \int_0^{\infty} \lambda^{s/4-1} \mu(\{w \in R\colon\, |T_{\mu, \delta}(1_{R}b_1, 1_Rb_2)(w)| > \lambda\})\,d\lambda \\
&\le \frac{s}{4} \bigg[  \int_0^1 \lambda^{s/4-1}\,d\lambda + C_0 \int_1^{\infty} \lambda^{-3s/4-1} \,d\lambda\bigg]\mu(R) \lesssim \mu(R).
\end{align*}
The desired bound
$$
|T_{\mu, \varepsilon}(1_Qb_1, 1_Qb_2)(x)| \lesssim C(\tau) +  M_{\mu, s/4}(1_QT_{\mu, \delta}(1_Qb_1, 1_Qb_2))(x)
$$
has now been proved.
\end{proof}

The following corollary contains the improved testing i.e. testing for $T_{\mu, \sharp}$.
\begin{corollary}\label{cor:cot}
Let $\mu$ be a measure of degree $m$ on $\R^n$, $T$ be a bilinear $m$-dimensional SIO and $s > 0$.
Let $t$ be a large enough constant depending only on the dimension $n$. Suppose
$b_i \in L^{\infty}(\mu)$, $i = 1, 2$, satisfy for every cube $R$ with $t$-small boundary that
\begin{equation}\label{eq:WeakCorCot}
\sup_{\delta > 0} \sup_{\lambda >0} \lambda^s \mu(\{x \in R\colon\, |T_{\mu, \delta}(1_Rb_1, 1_Rb_2)(x)| > \lambda\}) \lesssim \mu(2R).
\end{equation}
Let $t_0$ be another small boundary parameter. Then for every cube $Q$ with $t_0$-small boundary we have that
\begin{equation}\label{eq:CorCotConcl}
 \int_{Q} [T_{\mu, \sharp}(1_Qb_1, 1_Qb_2)]^{s/2}\,d\mu \lesssim \mu(5Q).
\end{equation}

\end{corollary}

\begin{proof}
Fix a cube $Q$ with $t_0$-small boundary. Let $a \in [2,2.2]$ be such that the cube $aQ$ has $t$-small boundary (such $a$ exists by Lemma 9.43 of \cite{ToBook}). First do the splitting
$$
T_{\mu, \sharp}(1_{Q}b_1,1_{Q}b_2)
\le T_{\mu, \sharp}(1_{aQ}b_1,1_{aQ}b_2)
+T_{\mu, \sharp}(1_{aQ\setminus Q} b_1,1_Qb_2)
+T_{\mu, \sharp}(1_{aQ}b_1,1_{aQ\setminus Q} b_2).
$$
We show that each of these three terms satisfies the desired estimate.

Because $Q$ has $t_0$-small boundary there holds  (see Lemma 9.44 in \cite{ToBook})
\begin{equation}\label{eq:SmallBoundComp}
\begin{split}
\int_Q T_{\mu, \sharp}(1_{aQ\setminus Q} b_1,1_Qb_2)^{s/2} d \mu
&\lesssim \int_Q \Big[\int_{aQ \setminus Q} \int_Q \frac{d\mu(z) d\mu(y) }{(|x-y|+|x-z|)^{2m}} \Big]^{s/2} d\mu(x) \\
&\lesssim \int_Q \Big[\int_{aQ \setminus Q} \frac{d\mu(y) }{|x-y|^{m}} \Big]^{s/2} d\mu(x)  \\
& \lesssim t_0 \mu(2Q).
\end{split}
\end{equation}
Similarly we have
$$
\int_Q T_{\mu, \sharp}(1_{aQ}b_1,1_{aQ\setminus Q} b_2)^{s/2} d \mu
\lesssim t_0 \mu (2Q).
$$

Fix some $\delta > 0$. We apply Proposition \ref{prop:cot} with the cube $aQ$ and the parameter $\tau=1/2$ to have that
$$
T_{\mu, \sharp, \delta}(1_{aQ}b_1, 1_{aQ}b_2)(x) \lesssim 1 + M_{\mu, s/4}^{\mathcal{Q}}(1_{aQ}T_{\mu, \delta}(1_{aQ}b_1, 1_{aQ}b_2))(x)
$$
for every $x \in (a/2)Q$, especially for all $x \in Q$. Therefore, we get
\begin{equation}\label{eq:ApplyCot}
\begin{split}
\int_{Q } [T_{\mu, \sharp, \delta}&(1_{aQ}b_1, 1_{aQ}b_2)]^{s/2}\,d\mu \\
&\lesssim \mu(Q) + \int_{\R^n} M_{\mu}^{\mathcal{Q}}(1_{aQ}|T_{\mu, \delta}(1_{aQ}b_1, 1_{aQ}b_2)|^{s/4})^2\,d\mu \\
&\lesssim \mu(Q) +  \int_{aQ} |T_{\mu, \delta}(1_{aQ}b_1, 1_{aQ}b_2)|^{s/2}\,d\mu 
\lesssim  \mu(5Q).
\end{split}
\end{equation}
The last estimate used the calculation in Proposition \ref{prop:cot} showing that weak type testing implies strong type testing for strictly smaller exponents.
Letting $\delta \to 0$ yields by monotone convergence that
$$
 \int_{Q} [T_{\mu, \sharp}(1_{aQ}b_1, 1_{aQ}b_2)]^{s/2}\,d\mu \lesssim \mu(5Q).
$$
This concludes the proof.
\end{proof}

\begin{remark}\label{rem:StrongCorCot}
If in Corollary \ref{cor:cot} one replaces the weak type assumption \eqref{eq:WeakCorCot} with the strong type condition
\begin{equation}\label{eq:StrongCorCot}
\sup_{\delta > 0} \int_R |T_{\mu, \delta}(1_Rb_1, 1_Rb_2)|^s d \mu \lesssim \mu(2R),
\end{equation}
then one gets the conclusion with the same exponent $s$, that is, instead of \eqref{eq:CorCotConcl} one has the conclusion
\begin{equation*}
 \int_{Q} [T_{\mu, \sharp}(1_Qb_1, 1_Qb_2)]^{s}\,d\mu \lesssim \mu(5Q).
\end{equation*}
This is proved similarly as above, except that in \eqref{eq:ApplyCot} one directly gets to apply the testing condition \eqref{eq:StrongCorCot} without the need to dominate the integral by the weak type testing in the last step. 
\end{remark}
\chapter{Suppressed bilinear singular integrals}\label{sec:suppopp}

Given a $1$-Lipschitz function $\Phi \colon \R^n \to [0,\infty)$ we define
$$
A_{\Phi}(x,y,z) =  \frac{(|x-y|+|x-z|)^{3\beta}}{ (|x-y|+|x-z|)^{3\beta}+ \Phi(x)^\beta \Phi(y)^\beta \Phi(z)^\beta},
$$
where $\beta = \beta(m) = \max(1, 2m/3)$. Given a standard $m$-dimensional bilinear Calder\'on--Zygmund kernel $K$ we define the suppressed kernel
$$
K_{\Phi}(x,y,z) = A_{\Phi}(x,y,z)K(x,y,z).
$$
It is important to understand that $A_{\Phi}(x,y,z) = 1$ if $\Phi(x) = 0$, say.
\begin{lemma}
The function $K_{\Phi}$ is a standard $m$-dimensional bilinear Calder\'on--Zygmund kernel with constants independent of the choice of the
$1$-Lipschitz function $\Phi$. Moreover, $K_{\Phi}$ satisfies the following improved
size condition
\begin{equation}\label{eq:impsize}
|K_{\Phi}(x,y,z)| \lesssim \frac{1}{ (|x-y|+|x-z|+ \Phi(x)+\Phi(y)+\Phi(z))^{2m}}.
\end{equation}
\end{lemma}

\begin{proof}

We begin with the size condition \eqref{eq:impsize}. Let $x,y,z \in \R^n$. We show that
\begin{equation}\label{eq:key_to_size}
\begin{split}
(|x-y|&+|x-z|)^{3\beta}+ \Phi(x)^\beta \Phi(y)^\beta \Phi(z)^\beta \\
& \gtrsim (|x-y|+|x-z|)^{3\beta}+ \Phi(x)^{3\beta} +\Phi(y)^{3\beta} +\Phi(z)^{3\beta},
\end{split}
\end{equation}
which clearly holds if
$$
\max(|x-y|,|x-z|,|y-z|) \geq \frac{1}{2} \max ( \Phi(x),\Phi(y),\Phi(z)).
$$
Suppose for example that 
$$
\max(|x-y|,|x-z|,|y-z|) \le \frac{1}{2} \max ( \Phi(x),\Phi(y),\Phi(z))= \frac{1}{2} \Phi(y).
$$
In this case
$$
\Phi(x) \geq \Phi(y)-|x-y| \geq \Phi(y)/2,
$$
and similarly $\Phi(z) \geq \Phi(y)/2$. Hence $$ \Phi(x)^\beta \Phi(y)^\beta \Phi(z)^\beta \gtrsim \Phi(y)^{3\beta},$$ whence
\begin{align*}
(|x-y|&+|x-z|)^{3\beta}+ \Phi(x)^\beta \Phi(y)^\beta \Phi(z)^\beta  \\
& \gtrsim (|x-y|+|x-z|)^{3\beta}+ \Phi(x)^{3\beta}+\Phi(y)^{3\beta}+\Phi(z)^{3\beta}.
\end{align*}

Using \eqref{eq:key_to_size} we have 
\begin{equation*}
\begin{split}
|K_ \Phi(x,y,z) |
& \lesssim \frac{ (|x-y|+|x-z|)^{3\beta-2m}}{ (|x-y|+|x-z|)^{3\beta}+ \Phi(x)^{3\beta} +\Phi(y)^{3\beta} +\Phi(z)^{3\beta}} \\
&\lesssim  \frac{1}{ (|x-y|+|x-z|+ \Phi(x)+\Phi(y)+\Phi(z))^{2m}},
\end{split}
\end{equation*}
where we applied the fact $\beta \geq 2m/3 $.

We turn to the H\"older conditions. Let $x,x',y,z \in \R^n$ be such that $|x-x'| \le \max(|x-y|, |x-z|)/2$. We have
\begin{equation*}
\begin{split}
|K_ \Phi(x',y,z)-K_\Phi(x,y,z)|
& \leq \big|K(x',y,z) \big( A_\Phi (x',y,z)-A_\Phi(x,y,z)\big) \big| \\
&+ \big|\big(K(x',y,z)-K(x,y,z)\big)A_\Phi (x,y,z) \big| \\
&=: I+ II.
\end{split}
\end{equation*}

We may use the $x$-H\"older condition of $K$ to get
\begin{equation*}
\begin{split}
 II
 \lesssim \frac{|x-x'|^\alpha}{(|x-y|+|x-z|)^{2m+\alpha}}, 
\end{split}
\end{equation*}
since $|A_\Phi| \leq 1$.

Consider then $I$. The size estimate of $K$  gives
\begin{equation}\label{eq:usegrad}
\begin{split}
I
\lesssim
\frac{1}{(|x-y|+|x-z|)^{2m}} \big| A_\Phi (x',y,z)-A_\Phi(x,y,z)\big|.
\end{split}
\end{equation}
Define the mapping $\gamma \colon [0,1] \to \R^n$ by setting 
$$
\gamma (t) =tx'+(1-t)x.
$$
We can write the difference to be estimated as
\begin{equation}\label{eq:FundCalc}
 A_\Phi (x',y,z)-A_\Phi(x,y,z)
=  \int_0^1 \frac{d}{dt} A_\Phi(\gamma(t),y,z) dt .
\end{equation}
Using the Lipschitz property of $\Phi$ one can check that $t \mapsto A_\Phi(\gamma(t),y,z)$ is absolutely continuous, whence this formula is valid.

Define
$$
a(t):= (|\gamma(t)-y|+|\gamma(t)-z|)^{3\beta}
$$
and 
$$
b(t):=  \Phi(\gamma(t))^\beta \Phi(y)^\beta \Phi(z)^\beta,
$$
which gives
$$
A_\Phi (\gamma(t),y,z)= \frac{a(t)}{a(t)+b(t)}.
$$
The $t$-derivative  can be written as 
$$
\frac{d}{dt} A_\Phi(\gamma(t),y,z)= \frac{a'(t)b(t)-a(t)b'(t)}{(a(t)+b(t))^2}.
$$

Computation of the derivatives gives
$$
a'(t)= 3\beta (|\gamma(t)-y|+|\gamma(t)-z|)^{3\beta-1}\Big(\frac{\gamma(t)-y}{|\gamma(t)-y|}\cdot (x'-x)+ \frac{\gamma(t)-z}{|\gamma(t)-z|} \cdot(x'-x)\Big)
$$
and 
$$
b'(t)
= \beta \Phi(\gamma(t))^{\beta-1} \frac{d}{d t}\Phi(\gamma(t)) \Phi(y)^\beta \Phi(z)^\beta.
$$
Notice also that $$
|\frac{d}{d t}\Phi(\gamma(t))| \leq |x'-x|
$$
by the  $1$-Lipschitz property of $\Phi$. These give us the estimates
\begin{equation*}
\begin{split}
| a'(t) b(t)|
&\lesssim (|\gamma(t)-y|+|\gamma(t)-z|)^{3\beta-1}|x'-x| \Phi(\gamma(t))^\beta \Phi(y)^\beta \Phi(z)^\beta \\
& \leq \big( |\gamma(t)-y|+|\gamma(t)-z|+ \Phi(\gamma(t))+\Phi(y) +\Phi(z)\big)^{6\beta-1}|x'-x|
\end{split}
\end{equation*}
and
\begin{equation*}
\begin{split}
|a(t) b'(t)| 
&\lesssim (|\gamma(t)-y|+|\gamma(t)-z|)^{3\beta}\Phi(\gamma(t))^{\beta-1}|x'-x| \Phi(y)^\beta \Phi(z)^\beta \\
& \leq \big( |\gamma(t)-y|+|\gamma(t)-z|+ \Phi(\gamma(t))+\Phi(y) +\Phi(z)\big)^{6\beta-1}|x'-x|,
\end{split}
\end{equation*}
where we used the fact that $\beta -1 \geq 0$. 

Combining the above estimates and using \eqref{eq:key_to_size}, we have shown that
\begin{equation}\label{eq:estpartial}
\begin{split}
\Big|\frac{d}{dt} &A_\Phi(\gamma(t),y,z)\Big| \\
&\lesssim 
\frac{\big( |\gamma(t)-y|+|\gamma(t)-z|+ \Phi(\gamma(t))+\Phi(y) +\Phi(z)\big)^{6\beta-1}|x'-x|}{\big((|\gamma(t)-y|+|\gamma(t)-z|)^{3\beta}+ \Phi(\gamma(t))^{3\beta} +\Phi(y)^{3\beta} +\Phi(z)^{3\beta}\big)^2}  \\
&\sim \frac{|x'-x|}{|x-y|+|x-z|+ \Phi(x)+\Phi(y) +\Phi(z)}.
\end{split}
\end{equation}
Applying this in \eqref{eq:usegrad} and \eqref{eq:FundCalc} leads to
\begin{equation*}
\begin{split}
I 
& \lesssim 
\frac{1}{(|x-y|+|x-z|)^{2m}} \frac{|x'-x|}{|x-y|+|x-z|+ \Phi(x)+\Phi(y) +\Phi(z)}  \\
& \le \frac{|x'-x|}{(|x-y|+|x-z|)^{2m+1}} \\
&\leq\frac{|x'-x|^\alpha}{(|x-y|+|x-z|)^{2m+\alpha}}.
 \end{split}
\end{equation*}
Hence $K_\Phi$ satisfies the $x$-H\"older estimate. In the same way one shows that $K_\Phi$ satisfies the other H\"older estimates.

\end{proof}

We define in the natural way
\begin{align*}
T_{\mu, \Phi, \varepsilon}(f,g)(x) &= \iint_{\max(|x-y|, |x-z|) > \varepsilon}  K_{\Phi}(x,y,z) f(y) g(z)\,d\mu(y)\,d\mu(z); \\
T_{\mu, \Phi, \sharp, \delta}(f,g)(x) &= \sup_{\varepsilon > \delta} |T_{\mu, \Phi, \varepsilon}(f,g)(x)|; \\
T_{\mu, \Phi, \sharp}(f,g)(x) &=  T_{\mu, \Phi, \sharp, 0}(f,g)(x).
\end{align*}
The following proposition is one of the key reasons why the suppressed operators are useful.
\begin{proposition}\label{prop:suppmain}
Let $\mu$ be a measure of degree $m$ on $\R^n$ and $T$ be a bilinear $m$-dimensional SIO. For a given $1$-Lipschitz function $\Phi$
there holds that
\begin{equation}\label{eq:suppmain}
T_{\mu, \Phi, \sharp}(f,g)(x) \le T_{\mu, \sharp, \Phi(x)}(f,g)(x) + CM_{\mu} f(x) M_{\mu} g(x).
\end{equation}
\end{proposition}
\begin{proof}
Fix $\delta > 0$ for which we will control $|T_{\mu, \Phi, \delta}(f,g)(x)|$ with a bound independent of $\delta$.
Assume first that $\Phi(x) \ge \delta$. Then we have that
\begin{align*}
T_{\mu, \Phi, \delta}(f,g)(x) &= \iint_{\max(|x-y|, |x-z|) > 2\Phi(x)}   K_{\Phi}(x,y,z) f(y) g(z)\,d\mu(y)\,d\mu(z) \\
&+ \iint_{\delta < \max(|x-y|, |x-z|) \le 2\Phi(x)}  K_{\Phi}(x,y,z) f(y) g(z)\,d\mu(y)\,d\mu(z).
\end{align*}
Notice that 
\begin{align*}
 & \iint_{\max(|x-y|, |x-z|) \le 2\Phi(x)}  |K_{\Phi}(x,y,z) f(y) g(z)|\,d\mu(y)\,d\mu(z)  \\
 &\lesssim \int_{\bar B(x, 2\Phi(x))}  \int_{\R^n} \frac{|f(y) g(z)|}{(\Phi(x) + |x-z|)^{2m}}\,d\mu(z) \,d\mu(y) \\
& \lesssim M_{\mu}  g(x) \cdot \frac{1}{\Phi(x)^m} \int_{B(x, 3\Phi(x))} |f(y)|\,d\mu(y)
\lesssim M_{\mu}  f(x) M_{\mu}  g(x).
\end{align*}
We also bound
\begin{align*}
\bigg|& \iint_{\max(|x-y|, |x-z|) > 2\Phi(x)}   K_{\Phi}(x,y,z) f(y) g(z)\,d\mu(y)\,d\mu(z) \bigg| \\
&\le \sup_{\varepsilon > \Phi(x)} \bigg| \iint_{\max(|x-y|, |x-z|) > \varepsilon} K_{\Phi}(x,y,z) f(y) g(z)\,d\mu(y)\,d\mu(z) \bigg| = T_{\mu, \Phi, \sharp, \Phi(x)}(f,g)(x).
\end{align*}
If it happens that $\Phi(x) < \delta$ we obviously have the bound
$$
|T_{\mu, \Phi, \delta}(f,g)(x)| \le  T_{\mu, \Phi, \sharp, \Phi(x)}(f,g)(x).
$$
So we have shown that for every $x \in \R^n$ there holds that
$$
|T_{\mu, \Phi, \delta}(f,g)(x)| \le  T_{\mu, \Phi, \sharp, \Phi(x)}(f,g)(x) + CM_{\mu}  f(x) M_{\mu}  g(x).
$$
Therefore, we are done if we show that
\begin{equation}\label{eq:showthis}
T_{\mu, \Phi, \sharp, \Phi(x)}(f,g)(x) \le T_{\mu, \sharp, \Phi(x)}(f,g)(x) + CM_{\mu}  f(x) M_{\mu}  g(x).
\end{equation}

To this end, we fix $\varepsilon > \Phi(x)$ and shall control $|T_{\mu, \Phi, \varepsilon}(f,g)(x)|$ with a bound independent of $\varepsilon$. 
Since now
$$
\bigg| \iint_{\max(|x-y|, |x-z|) > \varepsilon}  K(x,y,z) f(y) g(z)\,d\mu(y)\,d\mu(z) \bigg| \le T_{\mu, \sharp, \Phi(x)}(f,g)(x),
$$
the equation \eqref{eq:showthis} follows from showing that
$$
\iint_{\max(|x-y|, |x-z|) > \varepsilon}  |K(x,y,z)-K_{\Phi}(x,y,z)| |f(y)| |g(z)|\,d\mu(y)\,d\mu(z) \lesssim M_{\mu} f(x) M_{\mu} g(x).
$$
Notice that
$$
1-A_{\Phi}(x,y,z) \le \frac{\Phi(x)^\beta \Phi(y)^\beta \Phi(z)^\beta}{(|x-y| + |x-z|)^{3\beta}}
\lesssim \sum_{j=1}^3 \frac{\varepsilon^{j\beta}}{(|x-y| + |x-z|)^{j\beta}} 
$$
by using the definition of $A_{\Phi}(x,y,z)$, the $1$-Lipschitz property of $\Phi$ and the fact that $\Phi(x) < \varepsilon$. This implies that
\begin{align*}
&\iint_{\max(|x-y|, |x-z|) > \varepsilon}  |K(x,y,z)-K_{\Phi}(x,y,z)| |f(y)| |g(z)|\,d\mu(y)\,d\mu(z)\\
&\lesssim \sum_{j=1}^3 \varepsilon^{j\beta} \iint_{\max(|x-y|, |x-z|) > \varepsilon}  \frac{|f(y)| |g(z)|}{(|x-y|+|x-z|)^{2m+j\beta}} \,d\mu(y) \,d\mu(z)
\lesssim M_{\mu} f(x) M_{\mu} g(x),
\end{align*}
completing the proof.
\end{proof}

\begin{remark}
It follows that for every $1$-Lipschitz function the operator $T_{\mu, \Phi, \sharp}$ is bounded $$L^p(\mu) \times L^q(\mu) \to L^r(\mu),$$ if $T_{\mu, \sharp}$ is.
In particular, the weak boundedness property involving $T_{\mu, \Phi, \delta}$ is a reasonable condition.
\end{remark}
\subsection{$L^{\infty}$ suppression}
We now indicate how the key estimate \eqref{eq:suppmain} allows us -- in a proper sense -- to extend $L^{\infty}$ properties from a given set to the whole space.
In this chapter let $\mu$ be a \emph{finite} measure of order $m$.

Suppose $f_0,g_0$ are some fixed functions satisfying $|f_0|, |g_0| \le 1$.
Notice that $T_{\mu, \varepsilon}(f_0, g_0)(x)$ is, for every $x \in \R^n$ and $\varepsilon >0$, well-defined as an absolutely convergent integral (since $\mu$ is finite).
Let $S_0$ consist of those $x \in \R^n$ for which it holds that
$$
T_{\mu, \sharp}(f_0,g_0)(x) > \lambda_0.
$$
Here $\lambda_0 > 0$ is some fixed constant.
This means that $1_{\R^n \setminus S_0} T_{\mu, \sharp}(f_0,g_0) \le \lambda_0$ -- a property which we would like to have also in $S_0$. Of course,
just the opposite holds in $S_0$! However, if we choose $\Phi$ appropriately, then for some absolute constant $C$ we have
$$T_{\mu, \Phi, \sharp}(f_0,g_0) \le \lambda_0 + C$$ everywhere. Let us see how come. Notice that in $\R^n \setminus S_0$ everything is fine with any choice of $\Phi$.
Indeed, if $x \in \R^n \setminus S_0$, then simply
$$
T_{\mu, \Phi, \sharp}(f_0,g_0)(x) \le T_{\mu, \sharp}(f_0,g_0)(x) + C \le \lambda_0 + C.
$$
When controlling what happens in $S_0$ the choice of $\Phi$ becomes very relevant.

Define
$$
\varepsilon(x) = \sup\{\varepsilon > 0\colon\, |T_{\mu, \varepsilon}(f_0, g_0) (x)| > \lambda_0\}.
$$
If $x \in S_0$ then obviously $\varepsilon(x) > 0$. It also holds that $\varepsilon(x) < \infty$, since we have
$\lim_{\varepsilon \to \infty} T_{\mu, \varepsilon}(f_0, g_0) (x) = 0$ by monotone convergence. Define
$$
S = \bigcup_{x \in S_0} B(x, \varepsilon(x)).
$$
If $\Phi$ is any $1$-Lipschitz function satisfying that $\Phi(x) \ge d(x, S^c)$, then we are in business. To see this simply note that if $x \in S_0$, then
$\Phi(x) \ge \varepsilon(x)$, and so
$$
T_{\mu, \Phi, \sharp}(f_0,g_0)(x) \le T_{\mu, \sharp, \Phi(x)}(f_0,g_0)(x) + C \le \lambda_0 + C.
$$
We have shown that $T_{\mu, \Phi, \sharp}(f_0,g_0) \le \lambda_0 + C$ everywhere.

This is certainly extremely convenient. Of course, only if $S$ is not some horribly large set! To prevent this from happening
we need some additional, but rather weak, assumptions. It is enough that for some $s >0$ and $C_0 < \infty$ we have
$$
\sup_{\lambda > 0} \lambda^s \mu(\{x \in \R^n \setminus H \colon\, T_{\mu, \sharp}(f_0, g_0)(x) > \lambda\}) \le C_0 \mu(\R^n)
$$
for some set $H$ satisfying that $$\mu(H) \le \eta_0\mu(\R^n), \qquad \eta_0 < 1.$$ We also need to choose $\lambda_0 \lesssim 1$ large enough.
Indeed, we will show that for all large enough $\lambda_0$ we have
\begin{equation}\label{eq:Ssmall}
S \subset \{ x \in \R^n \colon T_{\mu, \sharp} (f_0, g_0)(x) > \lambda_0/2  \},
\end{equation}
and so
$$
\mu(S \setminus H) \le  \frac{2^sC_0}{\lambda_0^s}\mu(\R^n).
$$
This allows us to make sure that for large enough $\lambda_0$ we have $\mu(S^c) \sim \mu(\R^n)$, since $\mu(H) \le \eta_0 \mu(\R^n)$.

Let us show \eqref{eq:Ssmall}.
Let $x \in S$. Then there exists a point $x_0 \in S_0$ and a radius $\varepsilon_0$ such that $$
x \in B(x_0,\varepsilon_0) \textup{ and } |T_{\mu,\varepsilon_0}(f_0,g_0)(x_0)| > \lambda_0.
$$
The claim follows once we show that
$$ \big | T_{\mu,\varepsilon_0} (f_0,g_0)(x)-T_{\mu,\varepsilon_0} (f_0,g_0)(x_0) \big | \lesssim 1.$$ We have
\begin{equation*}
\begin{split}
\big | T_{\mu,\varepsilon_0} &(f_0,g_0)(x)-T_{\mu,\varepsilon_0} (f_0,g_0)(x_0) \big | \\
&\leq \big | T_{\mu,\varepsilon_0} (1_{B(x_0,2\varepsilon_0)}f_0,g_0)(x) \big | 
+\big | T_{\mu,\varepsilon_0} (1_{B(x_0,2\varepsilon_0)}f_0,g_0)(x_0) \big | \\
&+ \big | T_{\mu,\varepsilon_0} (1_{B(x_0,2\varepsilon_0)^c}f_0,g_0)(x)-T_{\mu,\varepsilon_0} (1_{B(x_0,2\varepsilon_0)^c}f_0,g_0)(x_0) \big |. 
\end{split}
\end{equation*}
Applying the size condition of the kernel there holds
\begin{equation*}
\begin{split}
\big | T_{\mu,\varepsilon_0} (1_{B(x_0,2\varepsilon_0)}f_0,g_0)(x) \big | 
&\lesssim \iint \frac{1_{B(x_0,2\varepsilon_0)}(y)}{\big(\varepsilon_0+|x-z|\big)^{2m}} \,d \mu(y) \,d \mu(z) \\
& \lesssim \frac{\mu(B(x_0,2\varepsilon_0))}{\varepsilon_0^m} 
\lesssim 1.
\end{split}
\end{equation*}
The corresponding term evaluated at $x_0$ is estimated in the same way.
The difference can be estimated as follows
\begin{equation*}
\begin{split}
\big | T_{\mu,\varepsilon_0} &(1_{B(x_0,2\varepsilon_0)^c}f_0,g_0)(x)-T_{\mu,\varepsilon_0} (1_{B(x_0,2\varepsilon_0)^c}f_0,g_0)(x_0) \big | \\
&= \Big | \iint \big(K(x,y,z)-K(x_0,y,z)\big) 1_{B(x_0,2\varepsilon_0)^c}(y)f_0(y) g_0(z) \,d\mu(y) \,d\mu(z) \Big | \\
& \lesssim \iint \frac{ \varepsilon_0 ^\alpha 1_{B(x_0,2\varepsilon_0)^c }(y)}{\big( |x_0-y|+|x_0-z| \big)^{2m+\alpha}} \,d\mu(y) \,d\mu(z)
\lesssim 1.
\end{split}
\end{equation*}
We have shown \eqref{eq:Ssmall}. 

In practice, we will need to do the above with slightly more generality. We formulate this as a separate proposition.
\begin{proposition}\label{prop:Phi_0}
Suppose $\mu$ is a finite measure of order $m$ and $T^1, \ldots, T^k$, $k \lesssim 1$, are bilinear $m$-dimensional SIO.
Let the functions $(f_0^1, g_0^1), \ldots, (f_0^k, g_0^k)$ satisfy
$$
|f_0^i|, |g_0^i| \lesssim 1, \qquad i = 1, \ldots, k.
$$
Suppose that for some $s >0$, $C_0 < \infty$ and for some set $H$ satisfying that $\mu(H) \le \eta_0\mu(\R^n)$, $\eta_0 < 1$, we have for every $i = 1, \ldots, k$ that
$$
\sup_{\lambda > 0} \lambda^s \mu(\{x \in \R^n \setminus H \colon\, T^i_{\mu, \sharp}(f_0^i, g_0^i)(x) > \lambda\}) \le C_0 \mu(\R^n).
$$
Then there exists a $1$-Lipschitz function $\Phi_0$ with
$$
\mu(\{\Phi_0 = 0\}) \sim \mu(\R^n),
$$ 
so that for all $1$-Lipschitz functions  $\Phi \geq \Phi_0$ there holds
$$
T^i_{\mu, \Phi, \sharp}(f_0^i,g_0^i)(x) \lesssim 1, \qquad x \in \R^n,
$$
for every $i = 1, \ldots, k$.

\end{proposition}
\begin{proof}
Let $S_i$ be the suppression sets with parameter $\lambda_0$ when we apply the above suppression procedure with $T$ replaced by $T^i$ and $f_0, g_0$ replaced with $f_0^i, g_0^i$.
Define $$
\Phi_0(x) = d(x, (S_1 \cup \cdots \cup S_k)^c),
$$
and let $\Phi$ be a $1$-Lipschitz function such that $\Phi \geq \Phi_0$. In particular, $\Phi(x) \ge d(x, S_i^c)$ and so
$T^i_{\mu, \Phi, \sharp}(f_0^i,g_0^i)(x) \le \lambda_0 + C$ for every $x$ and $i = 1, \ldots, k$.
Suppose that $\lambda_0 \lesssim 1$ is fixed large enough.
Then
$$
\mu(S_i \setminus H) \le  \frac{2^sC_0}{\lambda_0^s}\mu(\R^n) \le \frac{1-\eta_0}{2k}\mu(\R^n)
$$
and we get
$$
\mu( (S_1 \cup \cdots S_k)^c ) \ge \frac{1-\eta_0}{2}\mu(\R^n).
$$
\end{proof}
\begin{remark}
Given a bilinear $m$-dimensional SIO $T$
this proposition will be later applied with $k=3$, $T^1 = T$, $T^2 = T^{1*}$, $T^3 = T^{2*}$, and to some accretive $L^{\infty}$ functions.
\end{remark}

This completes our explanation of the $L^{\infty}$ suppression techniques. These ideas and the above calculations will be used concretely
to prove a certain bilinear big pieces $Tb$ theorem, which is the key to proving our main $Tb$ theorem.
\chapter{The big piece $Tb$}
In this chapter we prove Theorem \ref{thm:bigpieceTb} -- a very particular big piece type $Tb$ theorem adapted to our needs.
\begin{remark}\label{rem:matrixlemma}
In previous literature many summation arguments in $Tb$ theorems were based on summing the numbers
$$
\delta(J,R) := \frac{\ell(J)^{\alpha/2} \ell(R)^{\alpha/2}}{D(J, R)^{m+\alpha}}, \qquad D(J,R) := \ell(J) + \ell(R) + d(J,R),
$$
over all dyadic cubes. This was done in the $\ell^2$ sense by Nazarov--Treil--Volberg. We can also prove the $L^p$ analog, which goes as follows.
Let $\mu$ be a measure of order $m$ on $\R^n$, and $\calD$, $\calD'$ be two dyadic grids on $\R^n$.
For every $s \in (1,\infty)$ and $x_J \ge 0$, $J \in \calD$, we have
$$
\bigg\|\bigg( \sum_{R \in \calD'}1_R \bigg [ \sum_{J \in \calD} \delta(J,R) \mu(J) x_J  \bigg]^2 \bigg)^{1/2} \bigg\|_{L^s(\mu)}
\lesssim \bigg\|\bigg( \sum_{J \in \calD} x_J^2 1_J\bigg)^{1/2} \bigg\|_{L^s(\mu)}.
$$
However, we noticed a new simpler way to sum all the relevant parts in the $Tb$ argument, and no longer rely on this result.
\end{remark}
\begin{theorem}\label{thm:bigpieceTb}
Let $\mu$ be a measure of order $m$ such that $\mu(\R^n \setminus Q_0) = 0$ for some cube $Q_0 \subset \R^n$. Assume also that for some $t_0 < \infty$ we have for every
$\lambda > 0$ that
$$
\mu(\{x \in Q_0\colon\, d(x, \partial Q_0) \le \lambda\ell(Q_0)\}) \le t_0 \lambda \mu(Q_0).
$$
Let $T$ be a bilinear $m$-dimensional SIO, and let $b_i \in L^{\infty}(\mu)$, $i = 1,2,3$, be such that
$$
|\langle b_i \rangle_Q^{\mu}| \gtrsim 1 \qquad \textup{for all cubes } Q \subset Q_0.
$$
We assume the weak boundedness property in the form that
$$
|\langle T_{\mu, \Phi, \delta}(1_Qb_1, 1_Qb_2), 1_Qb_3\rangle_{\mu}| \lesssim \mu(5Q) \qquad \textup{for all cubes } Q \textup{ satisfying } 5Q \subset Q_0
$$
uniformly over the choice of the $1$-Lipschitz function $\Phi\colon\, \R^n \to [0, \infty)$ and the truncation parameter $\delta > 0$.
Let $s > 0$. Assume that there is a set $H \subset \R^n$ so that $\mu(H) \le \eta_0\mu(Q_0)$ for some $\eta_0 < 1$ and so that the following three testing conditions hold:
$$
\sup_{\lambda > 0} \lambda^s \mu(\{x \in Q_0 \setminus H\colon\, S_{\mu, \sharp}(b, b')(x) > \lambda\}) \lesssim \mu(Q_0)
$$
for all the choices $(S, b, b') \in \{(T, b_1, b_2), (T^{1*}, b_3, b_2), (T^{2*}, b_1, b_3)\}$.

Then there is a set $G \subset Q_0$ so that $\mu(G) \sim \mu(Q_0)$ and the following holds. For every $1 < p, q, r < \infty$ satisfying $1/p + 1/q = 1/r$
we have uniformly for functions $f \in L^p(\mu), g \in L^q(\mu)$ and   $h \in L^{r'}(\mu)$  supported in $G$ that
\begin{equation}\label{eq:bigpiece_conclusion}
\sup_{\varepsilon>0}|\langle T_{\mu, \varepsilon}(f,g), h\rangle_{\mu}| \lesssim  \| f \|_{L^p(\mu)}\| g \|_{L^q(\mu)}\| h \|_{L^{r'}(\mu)}.
\end{equation}
\end{theorem}

\begin{proof}

We begin by reducing the desired estimate to the boundedness of a certain suppressed operator. Let us apply Proposition \ref{prop:Phi_0} with $T_{\mu,\sharp}(b_1,b_2), T_{\mu,\sharp}^{1*}(b_3,b_2)$ and $T_{\mu,\sharp}^{2*}(b_1,b_3)$, and let $\Phi_0$ be the resulting $1$-Lipschitz function. For reasons that will become clear later, we have to modify the function $\Phi_0$ a little. Fix a small number $\lambda_0>0$ so that
$$
t_0 \lambda_0 \mu(Q_0) \leq \frac{\mu(\{ \Phi_0 =0 \})}{2}.
$$  
Clearly, we can choose $\lambda_0$ so that it only depends on the constants in our assumptions, since $\mu(\{ \Phi_0 =0 \}) \sim \mu(Q_0)$.
The choice of $\lambda_0$ implies by the small boundary assumption of $Q_0$ that
\begin{equation}\label{eq:smallboundary}
\mu(\{x \in Q_0\colon\, d(x, \partial Q_0) \le \lambda_0\ell(Q_0)\}) \leq \frac{\mu(\{ \Phi_0 =0 \})}{2}.
\end{equation}

Let $\phi$ be the $1$-Lipschitz function 
$$
\phi(x):= \max\big(\lambda_0 \ell(Q_0) - d(x, \partial Q_0), 0 \big).
$$
We have by \eqref{eq:smallboundary} that
$$
\mu(\{\phi \not=0)\} \leq \frac{\mu(\{ \Phi_0 =0 \})}{2}.
$$
Moreover, 
$$
\phi(x) \geq \frac{\lambda_0 \ell(Q_0)}{2} \qquad \textup{if } d(x, \partial Q_0)\leq \lambda_0  \ell(Q_0)/2.
$$

Define $\Phi_1:= \max(\Phi_0, \phi).$ Then $\Phi_1$ is a $1$-Lipschitz function, and the properties that we just verified for $\phi$ give that
$$
\mu(\{\Phi_1 =0\}) \geq \frac{\mu(\{ \Phi_0 =0 \})}{2}
$$
and 
\begin{equation}\label{eq:supprboundary}
\Phi_1(x) \geq \frac{\lambda_0 \ell(Q_0)}{2} \qquad \text{if } d(x, \partial Q_0) \leq \frac{\lambda_0 \ell(Q_0)}{2}.
\end{equation}
We define the set $G$ by setting
$$
G = \{x \in Q_0 \colon \Phi_1(x)=0\}.
$$
Since $\mu(G) \sim \mu(Q_0)$ it only remains to check \eqref{eq:bigpiece_conclusion}.

To this end, fix an arbitrary truncation parameter $\varepsilon>0$ and exponents $p,q,r \in (1, \infty)$ so that $1/p+1/q=1/r$.  The $1$-Lipschitz function that we will use in the suppression is
defined by
$$
\Phi = \max( \varepsilon, \Phi_1).
$$
Since $\Phi \geq \Phi_0$, Proposition \ref{prop:Phi_0} shows that 
$$
T_{\mu,\Phi,\sharp}(b_1,b_2)+ T_{\mu,\Phi,\sharp}^{1*}(b_3,b_2) +T_{\mu,\Phi,\sharp}^{2*}(b_1,b_3)\lesssim 1.
$$
We write $T_{\mu,\Phi}= T_{\mu,\Phi,0}$ which makes sense because $\Phi(x) \geq \varepsilon$ for every $x$.

Suppose $f \in L^p(\mu),g \in L^q(\mu)$ and let $x \in G$. Since $\Phi(x) = \varepsilon$, we have
\begin{equation*}
\begin{split}
\big|  T_{\mu, \varepsilon}(f,g)(x) - T_{\mu,\Phi} (f,g)(x) \big |
&\leq \big|  T_{\mu, \varepsilon}(f,g)(x) - T_{\mu,\Phi, \Phi(x)} (f,g)(x) \big | \\
&+ \big|  T_{\mu,\Phi, \Phi(x)} (f,g)(x) - T_{\mu,\Phi} (f,g)(x) \big | \\
& \lesssim M_{\mu} f(x) M_\mu g(x).
\end{split}
\end{equation*}
The required estimates for the final step can easily be read from the proof of Proposition \ref{prop:suppmain}.
This shows that if   $f \in L^p(\mu), g \in L^q(\mu)$ and $h \in L^{r'}(\mu)$ are functions supported in $G$, then
\begin{equation*}
\begin{split}
\big| \langle T_{\mu, \varepsilon}(f,g),h \rangle_\mu \big |
&\leq C \| M_{\mu} f M_\mu g\|_{L^r(\mu)} \| h \|_{L^{r'}(\mu)}
+\big| \langle T_{\mu,\Phi} (f,g),h \rangle_\mu \big | \\
& \leq C \| f \|_{L^p(\mu)} \| g \|_{L^q(\mu)}\| h \|_{L^{r'}(\mu)}
+\big| \langle T_{\mu,\Phi} (f,g),h \rangle_\mu \big |.
\end{split}
\end{equation*}
Thus, the required estimate \eqref{eq:bigpiece_conclusion} follows, if we show that there exists an absolute constant $C$, depending only on the constants in the assumptions, so that
\begin{equation}\label{eq:boundforsuppressed}
| \langle T_{\mu, \Phi}(f,g),h  \rangle_\mu |
\leq C \|f\|_{L^p(\mu)}\|g\|_{L^q(\mu)}\|h\|_{L^{r'}(\mu)}
\end{equation}
whenever $f \in L^p(\mu), g \in L^q(\mu)$ and $h \in L^{r'}(\mu)$ (not necessarily supported in $G$ anymore). We denote the best such constant $C$ by $\| T_{\mu, \Phi} \|$.  
Notice that we a priori know that $\| T_{\mu, \Phi} \| < \infty$, because $\Phi(x) \geq \varepsilon$ for every $x \in Q_0$. We are
after the quantitative bound.

\begin{remark}
Since $\Phi_1(x)=0$ for $x \in G$, we directly have the identity  
$$
\langle T_{\mu, \varepsilon}(f,g),h \rangle_\mu=\langle T_{\mu, \Phi_1, \varepsilon}(f,g),h \rangle_\mu
$$ 
if the functions are supported in $G$. However, we reduced to the operator $T_{\mu, \Phi}$ in order to get rid of the truncation
present in $T_{\mu, \Phi_1, \varepsilon}$. This is convenient -- see the proof of Lemma \ref{lem:zeroave_sep} to understand that zero average is
easier to utilise if there are no truncations in the integration.
\end{remark}

We start aiming towards \eqref{eq:boundforsuppressed}. For the moment, we fix three functions $f \in L^p(\mu), g \in L^q(\mu)$ and $h \in L^{r'}(\mu)$ with norm  at most $1$ so that
\begin{equation*}
| \langle T_{\mu, \Phi}(f,g),h  \rangle_\mu | 
\geq \frac{\| T_{\mu, \Phi} \|}{2}.
\end{equation*}
Define $$
Q_{0,\partial}:=\{x \in Q_0 \colon d(x,\partial Q_0) < \lambda_0 \ell(Q_0)/2\}.
$$
We may write $f = f_\partial + f_{\inte}$, where
$f_\partial := 1_{Q_{0,\partial}}f$ and $f_{\inte}:= f-f_\partial $. The functions $h_\partial,g_\partial, h_{\inte}$ and $g_{\inte}$ are defined similarly giving us the decomposition
\begin{equation*}
\begin{split}
\langle T_{\mu, \Phi}(f,g),h  \rangle_\mu 
&=  \langle T_{\mu, \Phi}(f_\partial,g),h  \rangle_\mu 
+  \langle T_{\mu, \Phi}(f_{\inte} ,g_\partial),h  \rangle_\mu \\
&+\langle T_{\mu, \Phi}(f_{\inte} ,g_{\inte}),h_\partial  \rangle_\mu 
+ \langle T_{\mu, \Phi}(f_{\inte} ,g_{\inte}),h_{\inte}  \rangle_\mu.
\end{split}
\end{equation*}
The last term is the main one. The first three terms are handled  using the suppression in the boundary region $Q_{0, \partial}$, where at least one of the appearing functions is supported in.
For example, we have by \eqref{eq:supprboundary}  and the improved size condition \eqref{eq:impsize} that
\begin{equation*}
\begin{split}
| \langle T_{\mu, \Phi}(f_{\inte} ,g_\partial),h  \rangle_\mu|
&\lesssim \frac{\| f_{\inte}\|_{L^1(\mu)}\| g_{\partial}\|_{L^1(\mu)}\| h\|_{L^1(\mu)}}{ \ell(Q_0)^{2m}} 
\lesssim 1.
\end{split}
\end{equation*}
The last step follows from H\"older's inequality using the fact that $\mu$ is of order $m$.
Thus, we have shown that there exists a constant $C$, depending only on the constants in our assumptions, so that
\begin{equation}\label{eq:nicesupports}
\frac{\| T_{\mu, \Phi} \|}{2} \leq C + \sup_{f,g,h} | \langle T_{\mu, \Phi}(f ,g),h  \rangle_\mu|,
\end{equation}
where the supremum is now over functions $f \in L^p(\mu), g \in L^q(\mu)$ and $h \in L^{r'}(\mu)$ with norm at most one and \emph{supported in $Q_0 \setminus Q_{0, \partial}$.}
The reason why we reduced to functions supported in $Q_0 \setminus Q_{0,\partial}$ is related to the fact that we can only use cubes which are inside $Q_0$ in our upcoming martingale decompositions. This is dictated by our assumptions.

We shall choose three such functions and split them using $b$-adapted martingales, and then perform the standard averaging argument of Nazarov--Treil--Volberg \cite{NTV} to reduce to ``good'' functions. 
\subsection*{Adapted martingales}
At this point we need to recall the random dyadic grids (these facts are essentially presented in this way by Hyt\"onen \cite{Hy}).
Let $\calD_{st}$ denote the standard dyadic grid, consisting of all the cubes of the form $2^{-k}(\ell + [0,1)^n)$, where $k \in \Z$ and $\ell \in \Z^n$. A generic dyadic grid, parametrized by
$$
\omega \in \Omega := (\{0,1\}^n)^\Z,
$$
is of the form
$$
\calD(\omega)= \cup_{k \in \Z} \calD_k(\omega), \textup{ where } \calD_k(\omega)=\{Q+x_k^{\omega}\colon\, Q \in \calD_{st, k}\} \textup{ and } x_k^{\omega}=\sum_{j > k} \omega_j 2^{-j}. 
$$
We get random dyadic grids by placing the natural product probability measure $\mathbb{P}_{\omega}$ on $\Omega = (\{0,1\}^n)^\Z$ (thus the coordinate functions $\omega_j$ are independent and $\mathbb{P}_{\omega}(\omega_j=\eta)=2^{-n}$ if $\eta \in \{0,1\}^n$). 

Let us consider some $\omega \in \Omega$ for the moment.
Choose the integer $u_0$ so that $2^{u_0} < \lambda_0 \ell(Q_0)/{4} \leq 2^{u_0+1}$.
We define the shorthand 
$$
\calD_0(\omega):= \{Q \in \calD(\omega) \colon Q \subset Q_0, \ell(Q) \leq 2^{u_0}\}.
$$ 

Suppose  $f$ is a locally $\mu$-integrable function. 
Define for every  $Q \in \calD_0(\omega)$  the numbers
$$
E^1_Qf:=\frac{\langle f \rangle^\mu_Q}{\langle b_1 \rangle^\mu_Q}.
$$
If $\ell(Q)=2^{u_0-1}$, we set 
$$
D^1_{Q}f:= E^1_Q f,
$$
and if $\ell(Q)<2^{u_0-1}$, then
$$
D^1_{Q}f:= E^1_Q f - E^1_{Q^{(1)}} f.
$$
Using these numbers the $b_1$-adapted martingale difference operators $\Delta^1_Q$ are defined for every $Q \in \calD_0(\omega)$ by setting
$$
\Delta^1_Qf:= \sum_{Q' \in \ch (Q)} (D^1_{Q'} f)1_{Q'} b_1.
$$
Notice the difference depending on  whether $\ell(Q)=2^{u_0}$ or $\ell(Q)<2^{u_0}$.

We will also need the adjoint operators of $\Delta^1_Q$. Namely, for $Q \in \calD_0(\omega)$ with $\ell(Q)=2^{u_0}$ define
$$
\Delta^{1*}_Q f
:=\sum_{Q' \in \ch (Q)}\frac{\langle fb_1 \rangle^\mu_{Q'}}{\langle b_1 \rangle^\mu_{Q'}}1_{Q'},
$$
and for $Q \in \calD_0(\omega)$ with $\ell(Q) < 2^{u_0}$ define
$$
\Delta^{1*}_Q f
:=\sum_{Q' \in \ch (Q)}\frac{\langle fb_1 \rangle^\mu_{Q'}}{\langle b_1 \rangle^\mu_{Q'}}1_{Q'}
-\frac{\langle f b_1\rangle^\mu_Q}{\langle b_1 \rangle^\mu_Q}1_{Q}.
$$

We list some properties of these martingale differences. Let $Q \in \calD_0(\omega)$. If $\ell(Q)<2^{u_0}$, then $\int \Delta^1_Q f \,d\mu =0$. For two functions $f$ and $g$  we have $\langle \Delta^1_Q f,g \rangle_\mu= \langle f, \Delta^{1*}_Q g \rangle_\mu$.
 If $R \in \calD_0 (\omega)$ is another cube, there holds
\begin{equation}\label{eq:orthog_mart}
\Delta^1_Q \Delta^1_R f =
\begin{cases}
0, \quad &\text{if } Q \not=R, \\
\Delta^1_Q f, \quad &\text{if } Q=R,
\end{cases}
\end{equation}
and
\begin{equation*}
\Delta^{1*}_Q \Delta^{1*}_R f =
\begin{cases}
0, \quad &\text{if } Q \not=R, \\
\Delta^{1*}_Q f, \quad &\text{if } Q=R.
\end{cases}
\end{equation*}

Suppose now $s \in (1,\infty)$. Let $f \in L^s(\mu)$ be a function whose support can be covered with cubes in $\calD_0(\omega)$ of side length $2^{u_0}$.  Then $f$ can be represented as  
$$
f= \sum_{Q \in \calD_0(\omega)} \Delta^1_Qf=\sum_{Q \in \calD_0(\omega)} \Delta^{1*}_Qf,
$$
where the convergence takes place unconditionally (that is, independently of the order) in $L^s(\mu)$.  
Define the $b_1$-adapted dyadic square functions 
\begin{equation*}
S^1_{\omega} f := \Big( \sum_{\substack{Q \in \calD_0(\omega) \\ \ell(Q) < 2^{u_0}}} | D^1_Qf|^21_Q \Big)^{1/2}
 \quad \text{and} \quad S^{1*}_{\omega} f := \Big( \sum_{Q \in \calD_0(\omega)} | \Delta^{1*}_Qf|^2 \Big)^{1/2}.
\end{equation*}
We have the standard estimates
\begin{equation}\label{eq:SFnorm}
\| f \|_{L^s(\mu)}
\lesssim  \| (S^1_{\omega} f) b_1\|_{L^s(\mu)}
\lesssim  \| S ^1_{\omega} f  \|_{L^s(\mu)}
\lesssim \| f \|_{L^s(\mu)}
\end{equation}
and
\begin{equation}\label{eq:SFDnorm}
\| f \|_{L^s(\mu)}
\sim  \| S^{1*}_{\omega} f  \|_{L^s(\mu)}.
\end{equation}
In \eqref{eq:SFnorm} the middle step is trivial using $b_1 \in L^\infty(\mu)$; the first and the last inequalities are the main facts. Notice the following consequence of this: If $\{a_Q\}_{Q \in \calD_0(\omega)}$ is a collection of real numbers so that 
$$
\Big(\sum_{Q \in \calD_0(\omega)} | a_Q \Delta^1_Q f|^2 \Big)^{1/2} 
$$ 
is in $L^s(\mu)$, then  $g:= \sum_{Q \in \calD_0(\omega)} a_Q \Delta^1_Q f$ is well defined in $L^s(\mu)$.  For every $Q \in \calD_0(\omega)$ there holds by \eqref{eq:orthog_mart} that $\Delta^1_Q g=a_Q \Delta^1_Q f$, and accordingly
\begin{equation}\label{eq:Burkholder}
\Big \| \sum_{Q \in \calD_0(\omega)} a_Q \Delta^1_Q f \Big \|_{L^s(\mu)}
\sim \Big \| \Big( \sum_{Q \in \calD_0(\omega)} | a_Q \Delta^1_Q f |^2 \Big)^{1/2} \Big \|_{L^s(\mu)}.
\end{equation}
The corresponding observation holds also with the operators $\Delta^{1*}_Q$.

For later use, we define the operators 
\begin{align*}
&E_{\omega, 2^k}^1 f := \sum_{\substack{Q \in \calD_0(\omega) \\ \ell(Q)=2^k}} 
(E^1_Qf) 1_{Q}b_1,    \quad \quad \ \ k \leq u_0 ,\\
&D^1_{\omega, 2^{u_0}} f:=E^1_{2^{u_0-1}}f, \\
&D^1_{\omega, 2^k}f:= E^1_{\omega, 2^{k-1}}f-E^1_{\omega, 2^k}f,   \quad \quad k <u_0.
\end{align*}

The corresponding definitions can be made using the functions $b_2$ and $b_3$ too.
Then we simply replace the super index $1$ above by $2$ or $3$ depending on the case.

Next, we recall the good and bad cubes of Nazarov-Treil-Volberg \cite{NTV}. If $\omega \in \Omega$, we say that a cube $Q$ is $\omega$-good with parameters $(\gamma,\sigma) \in (0,1) \times \Z_+$ if
$$
d(Q,R) > \ell(Q)^\gamma \ell(R)^{1-\gamma} \quad \text{for all } R \in \calD(\omega) \text{ with } \ell(R) \geq 2^{\sigma}\ell(Q);
$$
otherwise $Q$ is said to be $\omega$-bad. If $\omega_1, \omega_2 \in \Omega$, we say that $Q$ is $(\omega_1,\omega_2)$-good if it is \emph{both} $\omega_1$-good \emph{and} $\omega_2$-good; otherwise $Q$ is said to be $(\omega_1,\omega_2)$-bad. From now on, we fix $\gamma \in (0,1)$.
The parameter $\sigma$ will be fixed during the probabilistic argument below to be large enough but still $\lesssim 1$.

Let again $s \in (1, \infty)$ and $f \in L^s(\mu)$ with support in $Q_0 \setminus Q_{0, \partial}$. For $\omega =(\omega_1,\omega_2, \omega_3) \in \Omega \times \Omega \times \Omega$  we define 
$$
P_{\mathcal{B}}^1(\omega) f
:= \sum_{\substack{Q \in \calD_0(\omega_1) \\ Q \text{ is } (\omega_2,\omega_3)\text{-}\bad} }\Delta^1_Q f
$$
and
$$
P_{\mathcal{G}}^1(\omega) f
:= \sum_{\substack{Q \in \calD_0(\omega_1) \\ Q \text{ is } (\omega_2,\omega_3)\text{-}\good}}\Delta^1_Q f,
$$
where we keep in mind the dependence of these definitions on the (not yet fixed) goodness parameter $\sigma$.
We also define $P_{\mathcal{B}}^2(\omega)$ using the operators $\Delta^2_Q$ with cubes $Q \in \calD_0(\omega_2)$ that are $(\omega_1, \omega_3)$-bad
and $P_{\mathcal{G}}^3(\omega) $ using the operators $\Delta^3_Q $ with cubes $Q \in \calD_0(\omega_3)$ that are $(\omega_1, \omega_2)$-good and so on. 

On average, the norm of the bad part is small. The $L^2$ case is by Nazarov-Treil-Volberg \cite{NTV}.  The $L^p$ case is by Hyt\"onen \cite{Hy:vec}. An easier proof of the $L^p$ case
using interpolation and the $L^2$ case is by Lacey-V\"ah\"akangas \cite{LV}. We state the $L^p$-version adapted to our situation here:

\begin{lemma}\label{lem:ave_bad_norm}
Let $s \in (1,\infty)$ and $\omega_0 \in \Omega$. There exists a constant $\calP_\calB(\sigma,s)$, where $\sigma$ is the goodness parameter, such that
$$
\calP_{\calB}(\sigma,s) \to 0, \quad \text{as } \sigma \to \infty,
$$
and
$$
\E_\omega \Big \| \sum_{\substack{Q \in \calD_0(\omega_0) \\ Q \text{ is } \omega\text{-}\bad} } \Delta^1_Q f \Big\|_{L^s(\mu)}
\leq  \calP_{\calB}(\sigma,s) \| f \|_{L^s(\mu)}, \quad \text{for every } f \in L^s(\mu),
$$
where $\E_\omega$ denotes the expectation over $\omega \in \Omega$.
\end{lemma}

When it is clear from the context, we also denote by $\E_\omega:= \E_{\omega_1}\E_{\omega_2} \E_{\omega_3}$ the expectation over $\omega=(\omega_1, \omega_2, \omega_3) \in \Omega \times \Omega \times \Omega$. The result of Lemma \ref{lem:ave_bad_norm} can directly be extended to the case of three lattices. 
Indeed, suppose $f \in L^s(\mu)$. 
Since a $(\omega_2,\omega_3)$-bad cube is either $\omega_2$-bad or $\omega_3$-bad, we have 
$$
\Big(\sum_{\substack{Q \in \calD_0(\omega_1) \\ Q \text{ is } (\omega_2,\omega_3)\text{-}\bad} } |\Delta^1_Q f|^2 \Big)^{1/2}
\leq \Big( \sum_{\substack{Q \in \calD_0(\omega_1) \\ Q \text{ is } \omega_2\text{-}\bad} } |\Delta^1_Q f|^2 \Big)^{1/2}
+\Big(\sum_{\substack{Q \in \calD_0(\omega_1) \\ Q \text{ is } \omega_3\text{-}\bad} } |\Delta^1_Q f|^2 \Big)^{1/2}.
$$
Thus, there holds by \eqref{eq:Burkholder} that
\begin{equation*}
\begin{split}
\E_\omega \| P^1_\calB (\omega) f \|_{L^s(\mu)}
& \lesssim  \E_\omega \Big \| \sum_{\substack{Q \in \calD_0(\omega_1) \\ Q \text{ is } \omega_2\text{-}\bad} } \Delta^1_Q f \Big \|_{L^s(\mu)} 
+ \E_\omega \Big \| \sum_{\substack{Q \in \calD_0(\omega_1) \\ Q \text{ is } \omega_3\text{-}\bad} } \Delta^1_Q f \Big \|_{L^s(\mu)} \\
&= \E_{\omega_1}\E_{\omega_2}  \Big \| \sum_{\substack{Q \in \calD_0(\omega_1) \\ Q \text{ is } \omega_2\text{-}\bad} } \Delta^1_Q f \Big \|_{L^s(\mu)}
+\E_{\omega_1} \E_{\omega_3} \Big \| \sum_{\substack{Q \in \calD_0(\omega_1) \\ Q \text{ is } \omega_3\text{-}\bad} } \Delta^1_Q f \Big \|_{L^s(\mu)} \\
&\leq 2 \calP_{\calB}(\sigma,s) \| f \|_{L^s(\mu)},
\end{split}
\end{equation*}
where we applied Lemma \ref{lem:ave_bad_norm} in the last step. The same conclusion holds of course with the operators $P^2_\calB (\omega)$ and $P^3_\calB (\omega)$. 

Observe also that 
\begin{equation}\label{eq:norm_good}
\| P^i_\calG(\omega) f \|_{L^s(\mu)} \lesssim \| f \|_{L^s(\mu)}
\end{equation}
uniformly for $\omega \in \Omega \times \Omega \times \Omega$. This follows directly from Equation \eqref{eq:Burkholder}.

Having introduced the $b$-adapted martingales and the good and bad cubes, we continue with the proof of Theorem \ref{thm:bigpieceTb} from Equation \eqref{eq:nicesupports}. 
Consider three functions $f \in L^p(\mu), g \in L^q(\mu)$ and $h \in L^{r'}(\mu)$ with supports in $Q_0 \setminus Q_{0, \partial}$ and with norm at most one. 
For every $\omega \in \Omega \times \Omega \times \Omega$ we can split $\langle T_{\mu, \Phi}(f ,g),h  \rangle_\mu$, without denoting the dependence on $\omega$,  as
\begin{equation*}
\begin{split}
&\big \langle T_{\mu, \Phi}( P_{\mathcal{B}}^1 f ,g),h  \big \rangle_\mu
+\big \langle T_{\mu, \Phi}( P_{\mathcal{G}}^1 f , P_{\mathcal{B}}^2 g),h  \big \rangle_\mu \\
& \ \ + \big \langle T_{\mu, \Phi}( P_{\mathcal{G}}^1 f , P_{\mathcal{G}}^2 g),P_{\mathcal{B}}^3 h  \big\rangle_\mu
+  \big \langle T_{\mu, \Phi}(P_{\mathcal{G}}^1 f , P_{\mathcal{G}}^2 g),P_{\mathcal{G}}^3 h  \big \rangle_\mu.
\end{split}
\end{equation*}

On average the absolute value of the terms where there is at least one bad function involved is small. Indeed, applying Lemma \ref{lem:ave_bad_norm}  and Equation \eqref{eq:norm_good} we have for example that
\begin{equation*}
\begin{split}
\Big | \E_{\omega}& \bla T_{\mu, \Phi}( P_{\calG}^1(\omega) f , P_{\calG}^2(\omega) g),P_{\calB}^3(\omega) h  \bra_\mu \Big | \\
& \leq \E_{\omega} \| T_{\mu, \Phi} \| \| P_{\calG}^1(\omega) f \|_{L^p(\mu)} \| P_{\calG}^2(\omega) g\|_{L^q(\mu)} \| P_{\calB}^3(\omega) h\|_{L^{r'}(\mu)} \\
& \lesssim \| T_{\mu, \Phi} \| \| f \|_{L^p(\mu)} \|  g\|_{L^q(\mu)}  \E_\omega \| P_{\calB}^3(\omega) h\|_{L^{r'}(\mu)} \\
& \lesssim  \calP_{\calB}(\sigma,r') \| T_{\mu, \Phi}\|.
\end{split}
\end{equation*}
Thus, there exists a constant $C$ such that
\begin{equation}\label{eq:fraction+good}
\begin{split}
|\langle T_{\mu, \Phi}(f ,g),h  \rangle_\mu |
\leq C \calP_{\calB}(\sigma) \| T_{\mu, \Phi} \| + \Big| \E_\omega  \bla T_{\mu, \Phi}(P_{\mathcal{G}}^1 f , P_{\mathcal{G}}^2 g),P_{\mathcal{G}}^3 h \bra_\mu \Big|,
\end{split}
\end{equation}
where  $\calP_\calB(\sigma):= \max(\calP_\calB(\sigma,p),\calP_\calB(\sigma,q),\calP_\calB(\sigma,r'))$. 

By fixing the goodness parameter $\sigma$ to be big enough, there holds  $C\calP_{\calB}(\sigma) \leq 1/4$. Combining this with \eqref{eq:nicesupports} we have shown that
\begin{equation*}
\frac{ \| T_{\mu, \Phi}\|}{4} 
\leq C +  \sup_{f,g,h} \Big| \E_\omega  \bla T_{\mu, \Phi}(P_{\mathcal{G}}^1(\omega) f , P_{\mathcal{G}}^2 (\omega)g),P_{\mathcal{G}}^3(\omega) h \bra_\mu \Big|,
\end{equation*}
where the supremum is as in \eqref{eq:nicesupports}.
Now we fix three functions $f$, $g$ and $h$ as in the supremum, and turn to proving that
\begin{equation}\label{eq:reducedtogood}
\Big| \E_\omega  \bla T_{\mu, \Phi}(P_{\mathcal{G}}^1(\omega) f , P_{\mathcal{G}}^2 (\omega)g),P_{\mathcal{G}}^3(\omega) h \bra_\mu \Big| 
 \le C + \|T_{\mu, \Phi}\|/8.
\end{equation}
Once this is done, the proof of Theorem \ref{thm:bigpieceTb} is complete.

\subsection*{Proof of \eqref{eq:reducedtogood}}
For $\omega \in \Omega \times \Omega \times \Omega$ define
$f_\omega:=P_{\mathcal{G}}^1(\omega) f$, $g_\omega:= P_{\mathcal{G}}^2(\omega) g$ and $h_\omega:= P_{\mathcal{G}}^3(\omega) h$. Then there holds, for example, the identity
$$
f_\omega
= \sum_{Q \in \calD_0(\omega_1)}\Delta^1_Q f_\omega
= \sum_{\substack{Q \in \calD_0(\omega_1) \\ Q \text{ is } (\omega_2,\omega_3)\text{-}\good} }\Delta^1_Q f.
$$

We have
\begin{equation}\label{eq:sym1}
\begin{split}
\bla T_{\mu, \Phi}(f_\omega , g_\omega), h_\omega \bra_\mu 
&=\sum_{K \in \calD_0(\omega_3)} \sum_{\substack{I \in \calD_0(\omega_1) \\ \ell(K) \leq \ell(I)  }}
\sum_{\substack{J \in \calD_0(\omega_2) \\ \ell(K) \leq \ell(J) }} \bla T_{\mu, \Phi}(\Delta^1_I f_\omega , \Delta^2_Jg_\omega), \Delta^3_K h_\omega \bra_\mu \\
& + \sum_{I \in \calD_0(\omega_1)} \sum_{\substack{J \in \calD_0(\omega_2) \\ \ell(I) \leq \ell(J)  }}
\sum_{\substack{K \in \calD_0(\omega_3) \\ \ell(I) < \ell(K) }} \bla T^{1*}_{\mu, \Phi}(\Delta^3_K h_\omega , \Delta^2_Jg_\omega), \Delta^1_I f_\omega \bra_\mu \\
& +\sum_{J \in \calD_0(\omega_2)} \sum_{\substack{I \in \calD_0(\omega_1) \\ \ell(J) < \ell(I)  }}
\sum_{\substack{K \in \calD_0(\omega_3) \\ \ell(J) < \ell(K) }} \bla T^{2*}_{\mu, \Phi}(\Delta^1_I f_\omega , \Delta^3_Kh_\omega), \Delta^2_J g_\omega \bra_\mu.
\end{split}
\end{equation}
These three triple sums are essentially symmetric. We will concentrate on the first one.

Suppose $K \in \calD_0(\omega_3)$. The double sum $\sum_{\substack{I \in \calD_0(\omega_1) \\ \ell(K) \leq \ell(I)  }} \sum_{\substack{J \in \calD_0(\omega_2) \\ \ell(K) \leq \ell(J) }}$ 
can be organized as
$$
\sum_{\substack{I \in \calD_0(\omega_1) \\ \ell(K) \leq \ell(I)  }} \sum_{\substack{J \in \calD_0(\omega_2) \\ \ell(I) \leq \ell(J)  }}
+\sum_{\substack{J \in \calD_0(\omega_2) \\ \ell(K) \leq \ell(J) < 2^{u_0} }} \sum_{\substack{I \in \calD_0(\omega_1) \\ \ell(J) < \ell(I)  }}.
$$
Also, for  $I \in \calD_0(\omega_1)$  there holds
$$
\sum_{\substack{J \in \calD_0(\omega_2) \\ \ell(I) \leq \ell(J)  }} \Delta^2_J g_\omega= E_{\omega_2, \ell(I)/2}^2 g_\omega =: E_{\ell(I)/2}^2 g_\omega, 
$$
and for $J \in \calD_0(\omega_2)$ with  $\ell(J) <2^{u_0}$ we have
$$ 
\sum_{\substack{I \in \calD_0(\omega_1) \\ \ell(J) < \ell(I)  }} \Delta^1_I f_\omega= E_{\omega_1, \ell(J)}^1 f_\omega =: E_{\ell(J)}^1 f_\omega.
$$
Regarding the notation we make the following explanation. In what follows we always have some $\omega = (\omega_1, \omega_2, \omega_3)$ like here, and we
understand that e.g. $$E_{\omega_2, \ell(I)/2}^2 g_\omega = E_{\ell(I)/2}^2 g_\omega$$ i.e. that the superscript $2$ does not only
mean that we use the function $b_2$, but also that we use the dyadic grid $\calD_0(\omega_2)$. With this understanding we may suppress the additional subscripts denoting the dyadic grid used.
Now, combining the above with the bilinearity of $T_{\mu,\Phi}$ leads to
\begin{equation}\label{eq:sym2}
\begin{split}
\sum_{K \in \calD_0(\omega_3)}&\sum_{\substack{I \in \calD_0(\omega_1) \\ \ell(K) \leq \ell(I)  }}
\sum_{\substack{J \in \calD_0(\omega_2) \\ \ell(K) \leq \ell(J) }}  \bla T_{\mu, \Phi}(\Delta^1_I f_\omega , \Delta^2_Jg_\omega), \Delta^3_K h_\omega \bra_\mu \\
&= \sum_{K \in \calD_0(\omega_3)}\sum_{\substack{I \in \calD_0(\omega_1) \\ \ell(K) \leq \ell(I)  }} \bla T_{\mu, \Phi}(\Delta^1_I f_\omega , E^2_{\ell(I)/2}g_\omega), \Delta^3_K h_\omega \bra_\mu \\
&+ \sum_{\substack{K \in \calD_0(\omega_3) \\ \ell(K)< 2^{u_0}}}\sum_{\substack{J \in \calD_0(\omega_2) \\ \ell(K) \leq \ell(J) < 2^{u_0} }} \bla T_{\mu, \Phi}(E^1_{\ell(J)} f_\omega,\Delta^2_Jg_\omega), \Delta^3_K h_\omega \bra_\mu.
\end{split}
\end{equation}

Consider first the term
\begin{equation}\label{eq:reduced}
\E_\omega  \sum_{K \in \calD_0(\omega_3)}   \sum_{\substack{I \in \calD_0(\omega_1) \\ \ell(K) \leq \ell(I)  }} 
\bla T_{\mu, \Phi}(\Delta^1_I f_\omega , E^2_{\ell(I)/2}g_\omega), \Delta^3_K h_\omega \bra_\mu,
\end{equation}
where we note that because of goodness of the functions involved there are only good cubes in the summations.  
For every two cubes $Q,R \subset \R^n$ define the number $d_{Q,R}:= \max(2 \sqrt n \ell(Q), \ell(Q)^\gamma \ell(R)^{1-\gamma})$.
Related to the reduction into good cubes above, we can assume that the parameter $\sigma$ is so large that $2 \sqrt n \leq 2^{\sigma(1-\gamma)}$. In this case if $2^\sigma \ell(Q) \le \ell(R)$, then $d_{Q,R}=\ell(Q)^\gamma \ell(R)^{1-\gamma}$. Thus, if $K \in \DIII$ is $\omega_1$-good, $I \in \DI$, $2^\sigma \ell(K)\le \ell(I)$ and $d(K,I) \leq d_{K,I}$, then by goodness $K \subset I$. The reason why we use the numbers $d_{Q,R}$ is that they sometimes ensure that there is enough separation to use Lemma \ref{lem:zeroave_sep}.

Applying goodness, \eqref{eq:reduced} can be written as a sum of the following three terms: 
\begin{equation}\label{eq:2ndsymmetries}
\begin{split}
&I:= \E_\omega  \sum_{K \in \calD_0(\omega_3)}   \sum_{\substack{I \in \calD_0(\omega_1) \\ \ell(K) \leq \ell(I) \\ d(K,I) > d_{K,I} }} 
\bla T_{\mu, \Phi}(\Delta^1_I f_\omega , E^2_{\ell(I)/2}g_\omega), \Delta^3_K h_\omega \bra_\mu, \\
& II:= \E_\omega  \sum_{K \in \calD_0(\omega_3)}   \sum_{\substack{I \in \calD_0(\omega_1) \\ \ell(K) \leq \ell(I) \leq 2^{\sigma} \ell(K) \\ d(K,I) \leq d_{K,I} }} 
\bla T_{\mu, \Phi}(\Delta^1_I f_\omega , E^2_{\ell(I)/2}g_\omega), \Delta^3_K h_\omega \bra_\mu, \\
&III:= \E_\omega  \sum_{K \in \calD_0(\omega_3)}   \sum_{\substack{I \in \calD_0(\omega_1) \\ 2^{\sigma} \ell(K) < \ell(I) \\ K \subset I }} 
\bla T_{\mu, \Phi}(\Delta^1_I f_\omega , E^2_{\ell(I)/2}g_\omega), \Delta^3_K h_\omega \bra_\mu. 
\end{split}
\end{equation}
Goodness was needed to conclude that if $K \in \calD_0(\omega_3)$ is $\omega_1$-good and $I \in \calD_0(\omega_1)$ is such that $2^{\sigma} \ell(K) < \ell(I)$ and $d(K,I) \leq d_{K,I} =\ell(K)^\gamma \ell(I)^{1-\gamma}$, then $K \subset I$. In the term $II$ the average over $\omega \in \Omega \times \Omega \times \Omega$ will still be important, otherwise
we just estimate uniformly for every given $\omega \in \Omega \times \Omega \times \Omega$.

The other term in the right hand side of \eqref{eq:sym2} gives an essentially symmetric term as \eqref{eq:reduced}, and
can analogously be split into $I_{\text{sym}}$, $II_{\text{sym}}$ and  $III_{\text{sym}}$.
Previously (in older ArXiv versions) we wrote an argument where we directly estimated $I$, $II$ and $III$, 
which by symmetry
takes care of $I_{\text{sym}}$, $II_{\text{sym}}$ and  $III_{\text{sym}}$. 
However, in the followup work \cite{LMOV}
joint with Kangwei Li and Yumeng Ou, we noticed that  the paraproduct term arising from $III$ and the corresponding part from $III_{\text{sym}}$ can be combined to give one simple paraproduct. We will use this simpler way here. 

So the plan is a follows. We will estimate 
the terms $I$, $II$ and certain error terms related to $III$, which by symmetry takes care of
 $I_{\text{sym}}$, $II_{\text{sym}}$ and the error terms related to $III_{\text{sym}}$.
After these steps  we are left with two paraproduct type terms coming from $III$ and $III_{\text{sym}}$, which will be combined to give
one simple paraproduct.

\subsection*{An auxiliary estimate}
Before going into the analysis of the above three parts, let us record an easy estimate that is useful in what follows.
\begin{lemma}\label{lem:zeroave_sep}
Suppose $A \subset \R^n$ is a bounded set and $h_0$ is a function supported on $A$ such that $\int h_0 \,d\mu =0$.
Suppose also that $t \ge 2$ and $B \subset \R^{2n}$ is a set satisfying
$$
B \subset \{(y,z) \in \R^{2n}\colon\, \inf_{x \in A} \max(|x-y|, |x-z|) \ge t d(A)\}.
$$
Then we have for $f_0,g_0 \in L^1_{\text{loc}}(\mu)$ that
$$
\big |\bla \tilde{T}_{\mu,\Phi}(1_{B} f_0 \otimes g_0), h_0 \bra_\mu \big|
\lesssim \frac{1}{t^\alpha} \int M_{\mu, m}(f_0,g_0) | h_0 | \,d\mu.
$$
\end{lemma}
\begin{proof}
Applying the H\"older estimate in the $x$-variable, we have have for an arbitrary $x_A \in A$ that
\begin{equation*}
\begin{split}
&\big |\bla \tilde{T}_{\mu,\Phi}(1_{B} f_0 \otimes g_0), h_0 \bra_\mu \big| \\
&= \Big | \iiint \big(K_\Phi(x,y,z) -K_\Phi(x_A,y,z)\big)1_{B}(y,z) f_0(y) g_0(z) h_0(x) \,d\mu(x) \ud  \mu(y) \,d\mu(z) \Big| \\
& \lesssim \iiint \frac{d(A)^\alpha |f_0(y)g_0(z)h_0(x)|}{(td(A)+ |x-y|+|x-z|)^{2m+\alpha}} \,d\mu(y) \ud  \mu(z) \,d\mu(x) \\
& \lesssim \frac{1}{t^\alpha} \int M_{\mu, m}(f_0,g_0) | h_0 | \,d\mu.
\end{split}
\end{equation*}
The last estimate used Lemma \ref{lem:basic}.
\end{proof}

\subsection*{The separated sum} 
Now we begin with the term $I$.  Fix some $\omega \in  \Omega \times \Omega \times \Omega$. The sum over $K$ is  further divided into  the sum over those $K \in \calD_0(\omega_3)$ such that $\ell(K_0) < 2^{u_0}$ and those $K$ with $\ell(K)=2^{u_0}$. 

Suppose first that $K \in \calD_0(\omega_3)$ and $I \in \calD_0(\omega_1)$ are such that $\ell(K)=\ell(I)=2^{u_0}$ and $d(I,K) > \ell(K)^\gamma \ell(I)^{1-\gamma}= \ell(I)$. Then, applying directly the size condition of the kernel gives
\begin{equation*}
\begin{split}
|\bla T_{\mu, \Phi}(\Delta^1_I f_\omega , E^2_{\ell(I)/2}g_\omega), \Delta^3_K h_\omega \bra_\mu|
\lesssim \frac{\| \Delta_I^1f_\omega \| _{L^1(\mu)}\| E^2_{\ell(I)/2}g_\omega \| _{L^1(\mu)}\| \Delta_K^3h_\omega \| _{L^1(\mu)}}{\ell(I)^{2m}}.
\end{split}
\end{equation*}
Since $\| E^2_{\ell(I)/2}g_\omega \| _{L^1(\mu)}\lesssim \| g_\omega \|_{L^1(\mu)}$, summing over $I$ and $K$ leads to
\begin{equation}\label{eq:separated_top}
\begin{split}
\Big |\sum_{\substack{K \in \calD_0(\omega_3), I \in \calD_0(\omega_1)\\ \ell(I)=\ell(K)=2^{u_0} \\ d(K,I) > d_{K,I}  }} &  \bla T_{\mu, \Phi}(\Delta^1_I f_\omega , E^2_{\ell(I)/2}g_\omega), \Delta^3_K h_\omega \bra_\mu \Big| \\
&\lesssim   \sum_{\substack{I \in \calD_0(\omega_1) \\  \ell(I) =2^{u_0}}} \| \Delta_I^1f_\omega \| _{L^1(\mu)}
\cdot\frac{\|g_\omega \| _{L^1(\mu)}}{2^{2u_0m}} 
\cdot\sum_{\substack{K \in \calD_0(\omega_3)\\ \ell(K)=2^{u_0}} } \| \Delta_K^3h_\omega \| _{L^1(\mu)} \\
&\lesssim \frac{\|f_\omega \| _{L^1(\mu)}\|g_\omega \| _{L^1(\mu)}\| h_\omega \| _{L^1(\mu)}}{2^{2u_0m}}.
\end{split}
\end{equation}
H\"older's inequality gives 
\begin{equation*}
\begin{split}
\|f_\omega \| _{L^1(\mu)}\|g_\omega \| _{L^1(\mu)}\| h_\omega \| _{L^1(\mu)} 
&\leq \mu(Q_0)^2 \|f_\omega \| _{L^p(\mu)}\|g_\omega \| _{L^q(\mu)}\| h_\omega \| _{L^{r'}(\mu)}.
\end{split}
\end{equation*}
Because $2^{u_0} \sim \ell(Q_0)$ there holds $\mu(Q_0)^2/(2^{2u_0m}) \lesssim 1$. Therefore, the left hand side of \eqref{eq:separated_top} can be estimated as
\begin{equation}\label{eq:est_sep_top}
LHS\eqref{eq:separated_top} \lesssim   \|f_\omega \| _{L^p(\mu)}\|g_\omega \| _{L^q(\mu)}\| h_\omega \| _{L^{r'}(\mu)}
\lesssim \|f \| _{L^p(\mu)}\|g \| _{L^q(\mu)}\| h \| _{L^{r'}(\mu)}.
\end{equation}

Now we turn to those $K$ with $\ell(K)<2^{u_0}$. In this case we know that $\int \Delta^3_K h_\omega\,d\mu = 0$.
Define for the moment the shorthand 
$$
\varphi_{K,l}:= \sum_{\substack{I \in \calD_0(\omega_1) \\  2^l\ell(K) = \ell(I)  \\ d(K,I) >  d_{K,I}}} \Delta^1_I f_\omega,
$$
where $K \in \calD_0(\omega_3)$, $\ell(K_0) < 2^{u_0}$ and $l = 0, 1, \ldots, u_0 - \log_2 \ell(K)$.
Notice that $|\varphi_{K,l}| \leq |D^1_{2^l \ell(K)} f_\omega|$, and that $d(\supp \varphi_{K,l}, K) > d_{K,I}$. Lemma \ref{lem:zeroave_sep} gives 
\begin{equation}\label{eq:separ_begin}
\begin{split}
\Big| \Big \langle T_{\mu, \Phi} & ( \varphi_{K,l} , E^2_{2^{l-1}\ell(K)}g_\omega ), \Delta^3_K h_\omega \Big \rangle_\mu \Big| \\
& \lesssim   2^{-l(1-\gamma)\alpha} \int M_{\mu,m}( D_{2^l \ell(K)}^1 f_\omega) M_{\mu,m} M_{\mu, \calD_0(\omega_2)} g_\omega  | \Delta_K^3h_\omega | \,d\mu.
\end{split}
\end{equation}

Next, we show that
\begin{equation}\label{eq:separ_error}
\begin{split}
\sum_{K \in \calD_0(\omega_3)} \sum_{l=0}^{u_0-\log_2\ell(K)}2^{-l(1-\gamma)\alpha}
&\int M_{\mu,m}( D_{2^l \ell(K)}^1 f_\omega) M_{\mu,m} M_{\mu, \calD_0(\omega_2)} g_\omega  | \Delta_K^3h_\omega | \,d\mu \\
&\lesssim \| f\|_{L^p(\mu)} \| g\|_{L^q(\mu)}\| h\|_{L^{r'}(\mu)}.
\end{split}
\end{equation}
The left hand side of \eqref{eq:separ_error} can be reorganized as
\begin{equation*}
\begin{split}
\sum_{l=0}^\infty 2^{-l(1-\gamma)\alpha} 
&\sum_{k\colon k \leq u_0 - l} 
\sum_{\substack{K \in \calD_0(\omega_3) \\ \ell(K) =2^{k}}}
\int M_{\mu,m}( D_{2^{k+l}}^1 f_\omega) M_{\mu,m} M_{\mu, \calD_0(\omega_2)} g_\omega  | \Delta_K^3h_\omega | \,d\mu \\
&=\sum_{l=0}^\infty 2^{-l(1-\gamma)\alpha} 
\sum_{k\colon k \leq u_0 - l}  \int M_{\mu,m}( D_{2^{k+l}}^1 f_\omega) M_{\mu,m} M_{\mu, \calD_0(\omega_2)} g_\omega  | D^3_{2^k}h_\omega | \,d\mu.
\end{split}
\end{equation*}
Fix one $l \in \{0,1,\dots \}.$ We have
\begin{equation*}
\begin{split}
&\sum_{k\colon k \leq u_0 - l}  \int M_{\mu,m}( D_{2^{k+l}}^1 f_\omega) M_{\mu,m} M_{\mu, \calD_0(\omega_2)} g_\omega  | D^3_{2^k}h_\omega | \,d\mu \\
& \leq 
\Big \| 
\Big( 
\sum_{k\colon k \leq u_0 - l} 
M_{\mu,m}( D_{2^{k+l}}^1 f_\omega)^2 \Big)^{1/2} \Big\|_{L^p(\mu)}
\| M_{\mu,m} M_{\mu, \calD_0(\omega_2)} g_\omega    \|_{L^q(\mu)} \\
&  \ \ \ \!  \cdot \Big \| 
\Big( 
\sum_{k\colon k \leq u_0 - l}
| D^3_{2^k}h_\omega |^2 \Big)^{1/2} \Big\|_{L^{r'}(\mu)} \\
&\lesssim \| f_\omega\|_{L^p(\mu)} \|g_\omega\|_{L^q(\mu)} \| h_\omega \| _{L^{r'}(\mu)},
\end{split}
\end{equation*}
where we applied the Fefferman--Stein inequality for the radial maximal function $M_{\mu,m}$.
Such a version of the Fefferman--Stein inequality follows e.g. by using that it is at least known to be true for non-homogeneous dyadic maximal functions,
and dominating the radial maximal function $M_{\mu,m}$ by finitely many such dyadic maximal functions.
To prove \eqref{eq:separ_error} it only remains to sum the geometric series  $\sum_{l \ge 0} 2^{-l(1-\gamma)\alpha}$.

Equations \eqref{eq:separ_begin} and \eqref{eq:separ_error} combined show that for a fixed $\omega \in \Omega \times \Omega \times \Omega$  the part of $I$ that consists of those $K$ with $\ell(K) < 2^{u_0}$ satisfies the right bound.
Since the estimates we have done have been independent of $\omega$, this finishes the proof of
\begin{equation}\label{eq:resSEP}
|I| \lesssim \| f\|_{L^p(\mu)} \| g\|_{L^q(\mu)}\| h\|_{L^{r'}(\mu)} \le 1.
\end{equation}

\subsection*{Deeply contained cubes; error terms}
Here a part of the term $III$ is considered. Again a uniform estimate will be made for every $\omega \in \Omega \times \Omega \times \Omega$. We fix one $\omega$  now until we have estimated the whole term $III$. The goal is to reduce the estimate to a so-called paraproduct that involves the function $T_{\mu, \Phi}(b_1, b_2)$, which will then allow us to apply the property $|T_{\mu, \Phi}(b_1,b_2) | \lesssim 1.$ To achieve this, we must first estimate two error terms. The paraproduct is handled in the next subsection.

Let us first introduce some notation. 
Define
\begin{equation*}
\begin{split}
\calD_h(\omega_3)  &= \{K \in \calD_0(\omega_3)\colon\, \ell(K) < 2^{u_0 - \sigma}, \Delta^3_K h_\omega \ne 0\} \\
& = \{K \in \calD_0(\omega_3)\colon\, \ell(K) < 2^{u_0 - \sigma}, K \text{ is } (\omega_1, \omega_2) \text{-good}, \Delta^3_K h \ne 0\}.
\end{split}
\end{equation*}
Suppose $K \in \calD_h(\omega_3)$ and $l \in \Z$ are such that $2^{\sigma}\ell(K) \leq 2^l\ell(K) \leq 2^{u_0}$. Then, by the goodness of the cube $K$, there exist cubes $I \in \calD(\omega_1)$ and $J \in \calD(\omega_2)$ of side length $2^l\ell(K)$ containing $K$. Since moreover $\Delta_K^3 h \ne 0$, which implies 
$$
K \subset \{d(\cdot, \R^n \setminus Q_0) > \lambda_0\ell(Q_0) / 4\},
$$
the above cubes actually satisfy $I \in \calD_0(\omega_1)$ and $J \in \calD_0(\omega_2)$. 
We denote these unique cubes $I$ and $J$ by $I_{K,l}$ and $J_{K,l}$, and sometimes also by $I(K,l)$ and $J(K,l)$. So, in this context ``$I$'' refers to the lattice $\calD(\omega_1)$ and ``$J$'' refers to $\calD(\omega_2)$.  These definitions depend on $\omega$, but it does not matter since it is fixed for the moment. 

We write
\begin{align*}
& \sum_{K \in \calD_0(\omega_3)} \sum_{\substack{I \in \calD_0(\omega_1) \\ 2^{\sigma} \ell(K) < \ell(I) \\ K \subset I }} 
\bla T_{\mu, \Phi}(\Delta^1_I f_\omega , E^2_{\ell(I)/2}g_\omega), \Delta^3_K h_\omega \bra_\mu \\
&= \sum_{K \in \calD_h(\omega_3)} \sum_{\substack{l \in \Z \colon \\ 2^{\sigma}\ell(K) < 2^l\ell(K) \le 2^{u_0}}}
\bla T_{\mu, \Phi}(\Delta^1_{I_{K,l}} f_\omega , E^2_{2^{l-1}\ell(K)}g_\omega), \Delta^3_K h \bra_\mu.
\end{align*}
Notice that 
$\Delta^1_{I_{K,l}} f_\omega \otimes E^2_{2^{l-1}\ell(K)}g_\omega$ can be written as the sum
\begin{equation*}
\begin{split}
&1_{I_{K,l-1} \times J_{K,l-1}} \Delta^1_{I_{K,l}} f_\omega \otimes E^2_{2^{l-1}\ell(K)}g_\omega    
+ 1_{(I_{K,l-1} \times J_{K,l-1})^c} \Delta^1_{I_{K,l}} f_\omega \otimes E^2_{2^{l-1}\ell(K)}g_\omega,
\end{split}
\end{equation*}
where further
\begin{equation*}
\begin{split}
1_{I_{K,l-1} \times J_{K,l-1}} \Delta^1_{I_{K,l}} f_\omega \otimes E^2_{2^{l-1}\ell(K)}g_\omega
&=(D^1_{I_{K,l-1}}f_\omega) b_1 \otimes (E^2_{J_{K,l-1}}g_\omega) b_2 \\
&-1_{(I_{K,l-1} \times J_{K,l-1})^c} (D^1_{I_{K,l-1}}f_\omega) b_1 \otimes (E^2_{J_{K,l-1}}g_\omega) b_2.
\end{split}
\end{equation*}
Here $(D^1_{I_{K,l-1}}f_\omega) b_1 \otimes (E^2_{J_{K,l-1}}g_\omega) b_2$ leads to the main term, and the other two give the error terms.

The error terms can easily be handled using Lemma \ref{lem:zeroave_sep}. Notice that
\begin{equation}\label{eq:error1}
\begin{split}
\Big|\Big \langle \tilde{T}_{\mu,\Phi} &\big(1_{(I_{K,l-1} \times J_{K,l-1})^c} \Delta^1_{I_{K,l}} f_\omega \otimes E^2_{2^{l-1}\ell(K)}g_\omega\big), \Delta_K^3h \Big \rangle_\mu \Big| \\
&\lesssim 2^{-l(1-\gamma)\alpha} \int M_{\mu,m}(D^1_{2^{l}\ell(K)} f_\omega) M_{\mu,m}M_{\mu, \calD_0(\omega_2)}g_\omega | \Delta^3_K h| \,d\mu
\end{split}
\end{equation}
and
\begin{equation}\label{eq:error2}
\begin{split}
\Big|\Big \langle \tilde{T}_{\mu,\Phi} &\big(1_{(I_{K,l-1} \times J_{K,l-1})^c} (D^1_{I_{K,l-1}}f_\omega) b_1 \otimes (E^2_{J_{K,l-1}}g_\omega) b_2\big), \Delta_K^3h \Big \rangle_\mu \Big| \\
& \lesssim 2^{-l(1-\gamma)\alpha}  |D^1_{I_{K,l-1}}f_\omega| |E^2_{J_{K,l-1}}g_\omega| \int | \Delta^3_Kh| \,d\mu \\
& \lesssim 2^{-l(1-\gamma)\alpha} \int |D^1_{2^{l}\ell(K)} f_\omega | M_{\mu, \calD_0(\omega_2)} g_\omega | \Delta^3_Kh| \,d\mu.
\end{split}
\end{equation}
It is clear that the sum over $K$ and $l$ of these is bounded by $\| f\|_{L^p(\mu)} \|g\|_{L^q(\mu)} \| h \| _{L^{r'}(\mu)}$ -- see \eqref{eq:separ_error}.

\subsection*{Deeply contained cubes; paraproduct}
We begin this subsection by proving the following lemma.
\begin{lemma}\label{lem:parest}
Let $B$ be a function satisfying $|B| \lesssim 1$. Then for every $\omega \in \Omega$ and $\varphi \in L^r(\mu)$ we have
\begin{equation}\label{eq:general_para}
\Big\|\sum_{Q \in \calD_0(\omega)}
\langle \varphi \rangle _Q^\mu
\Delta^{1*}_Q B  \Big \|_{L^{r}(\mu)} 
\lesssim \|  \varphi \|_{L^r(\mu)}.
\end{equation}
\end{lemma}
\begin{proof}
We begin by constructing the collection $\scrS \subset \calD_0(\omega)$ of principal cubes for the function $\varphi$. Of course, this is a completely standard construction. Set 
$$
\scrS_0 := \{ Q \in \calD_0(\omega) \colon \ell(Q)=2^{u_0}\},
$$
and suppose $\scrS_0, \dots , \scrS_k$ are defined for some $k$. If $S \in \scrS_k$, we define $\ch_\scrS (S)$ to be the maximal cubes $Q \in \calD_0(\omega)$ such that $Q \subset S$ and 
$$
\langle |\varphi| \rangle _Q^\mu > 2\langle |\varphi| \rangle _S^\mu.
$$
Then, set $\scrS_{k+1}:= \bigcup_{S \in \scrS_{k}} \ch_{\scrS}(S)$, and finally $\scrS:= \bigcup_{k=0}^\infty \scrS_k$.
The collection $\scrS$ is obviously a Carleson family of cubes with respect to the measure $\mu$.
For every $Q \in \calD_0(\omega)$ there exists a a cube $S \in \scrS$ such that $Q \subset S$; the minimal such $S$ is denoted by $\pi_\scrS Q$. By the construction of $\scrS$ there holds
$$
\langle |\varphi|\rangle _Q^\mu \leq 2\langle |\varphi| \rangle _{\pi_\scrS Q}^\mu, \quad Q \in \calD_0(\omega).
$$

Turning to \eqref{eq:general_para}, we can use the principal cubes together with \eqref{eq:Burkholder} to have
\begin{equation}\label{eq:paraprinc}
\begin{split}
\Big\|\sum_{Q \in \calD_0(\omega)}
\langle \varphi \rangle _Q^\mu
\Delta^{1*}_Q B  \Big \|_{L^{r}(\mu)} 
\lesssim 
\Big \| \sum_{S \in \scrS} \langle |\varphi| \rangle _S^\mu
\sum_{\substack{Q \in \calD_0(\omega) \\ \pi_\scrS Q =S}} \Delta^{1*}_Q B \Big \|_{L^r(\mu)}.
\end{split}
\end{equation}

Suppose $S \in \scrS$, $S' \in \ch_\scrS (S)$ and $x \in S'$. If $\ell(S) < 2^{u_0}$, there holds 
$$
\Big| \sum_{\substack{Q \in \calD_0(\omega) \\ \pi_\scrS Q =S}}
\Delta^{1*}_QB(x) \Big|
= \Big|\frac{\bla b_1B \bra_{S'}^{\mu}}{\langle b_1 \rangle_{S'}^{\mu}}
-\frac{\bla b_1B \bra_{S}^{\mu}}{\langle b_1 \rangle_{S}^{\mu}} \Big|
\lesssim 1,
$$
while if $\ell(S)=2^{u_0}$, then
$$
\Big| \sum_{\substack{Q \in \calD_0(\omega) \\ \pi_\scrS Q =S}}
\Delta^{1*}_QB(x) \Big|
= \Big|\frac{\bla b_1B \bra_{S'}^{\mu}}{\langle b_1 \rangle_{S'}^{\mu}}
\Big|
\lesssim 1.
$$
On the other hand, if $\ell(S) <2^{u_0}$, for  $\mu$-almost every $x \in S \setminus \bigcup_{S' \in \ch_\scrS(S)} S'$ we have
$$
\Big| \sum_{\substack{ Q \in \calD_0(\omega) \\ \pi_\scrS Q =S}}
\Delta^{1*}_QB(x) \Big|
= \Big|B(x) -\frac{\bla b_1B \bra_{S}^{\mu}}{\langle b_1\rangle_{S}^{\mu}} \Big|
\lesssim 1,
$$
and if $\ell(S) = 2^{u_0}$, then
$$
\Big| \sum_{\substack{ Q \in \calD_0(\omega) \\ \pi_\scrS Q =S}}
\Delta^{1*}_QB(x) \Big|
= \big|B(x)  \big|
\lesssim 1.
$$

Hence, the right hand side of \eqref{eq:paraprinc} is dominated by
$$
\Big \| \sum_{S \in \scrS} \langle |\varphi| \rangle _S^\mu
1_S \Big \|_{L^r(\mu)} \lesssim \| \varphi \|_{L^r(\mu)},
$$
where we applied the Carleson embedding theorem in a form that appears at least in Theorem 2.2 of \cite{Vu1}.
\end{proof}

We continue with the estimation of $III$. Recall that after the error terms handled in the last section what is left is
\begin{equation}\label{eq:BefCol1}
\sum_{K \in \calD_h(\omega_3)} \sum_{l = \sigma+1}^{u_0 - \log_2 \ell(K)}
D^1_{I_{K,l-1}}f_\omega E^2_{J_{K,l-1}}g_\omega \bla T_{\mu, \Phi}(b_1, b_2), \Delta_K^3 h\bra_{\mu}.
\end{equation}
Now it is time to combine this with the corresponding part coming from the other term on the right hand side of \eqref{eq:sym2}, which is
\begin{equation}\label{eq:BefCol2}
\sum_{K \in \calD_h(\omega_3)} \sum_{l = \sigma+1}^{u_0 - \log_2 \ell(K)-1}
E^1_{I_{K,l}}f_\omega D^2_{J_{K,l-1}}g_\omega \bla T_{\mu, \Phi}(b_1, b_2), \Delta_K^3 h\bra_{\mu}
\end{equation}
(see the discussion after \eqref{eq:2ndsymmetries}).

Let $K \in \calD_h(\omega_3)$. If $l \in \{ \sigma +1, \dots u_0- \log_2 \ell(K)-1\}$, then
\begin{equation*}
\begin{split}
D^1_{I_{K,l-1}}f_\omega &E^2_{J_{K,l-1}}g_\omega
+E^1_{I_{K,l}}f_\omega D^2_{J_{K,l-1}}g_\omega \\
&= E^1_{I_{K,l-1}}f_\omega E^2_{J_{K,l-1}}g_\omega
-E^1_{I_{K,l}}f_\omega E^2_{J_{K,l}}g_\omega.
\end{split}
\end{equation*}
If $l=u_0- \log_2 \ell(K)$, then the different definition of the martingale differences on the top level says that
$$
D^1_{I_{K,l-1}}f_\omega E^2_{J_{K,l-1}}g_\omega
=E^1_{I_{K,l-1}}f_\omega E^2_{J_{K,l-1}}g_\omega.
$$
Hence, we  see that \eqref{eq:BefCol1} and \eqref{eq:BefCol2} summed together give
\begin{equation}\label{eq:SimplePara}
\sum_{K \in \calD_h(\omega_3)}
E^1_{I_{K,\sigma}}f_\omega E^2_{J_{K,\sigma}}g_\omega
\bla T_{\mu, \Phi}(b_1, b_2), \Delta_K^3 h\bra_{\mu}.
\end{equation}

Notice that if $K \in \calD_h(\omega_3)$, then
$$
\big|
E^1_{I_{K,\sigma}}f_\omega E^2_{J_{K,\sigma}}g_\omega \big|
\lesssim \bla M_{\mu, \calD_0(\omega_1)} f_\omega \cdot M_{\mu, \calD_0(\omega_2)} g_\omega \bra_K^\mu.
$$
Hence,  the absolute value of \eqref{eq:SimplePara} 
is dominated by
\begin{equation*}
\begin{split}
\Big \| \sum_{K \in \calD_h(\omega_3)}
 &E^1_{I_{K,\sigma}}f_\omega E^2_{J_{K,\sigma}}g_\omega 
\Delta^3_K T_{\mu, \Phi}(b_1, b_2) \Big\|_{L^r(\mu)} \| h \|_{L^{r'}(\mu)} \\
&\lesssim \Big \| \sum_{K \in \calD_h(\omega_3)}
\bla M_{\mu, \calD_0(\omega_1)} f_\omega \cdot M_{\mu, \calD_0(\omega_2)} g_\omega \bra_K^\mu
\Delta^3_K T_{\mu, \Phi}(b_1, b_2) \Big\|_{L^r(\mu)} \| h \|_{L^{r'}(\mu)} \\
& \lesssim
\| f_\omega \|_{L^p(\mu)} \| g _\omega \|_{L^q(\mu)} 
\| h \|_{L^{r'}(\mu)},
\end{split}
\end{equation*}
where we applied Equation \eqref{eq:Burkholder} in the first inequality and Lemma \ref{lem:parest} in the second.

This finishes the estimate for the term $III+III_{\text{sym}}$.

\subsection*{The diagonal}
Here we need to deal with the final term
\begin{equation}\label{eq:diag_collapsed}
II = \E_\omega  \sum_{K \in \calD_0(\omega_3)}   
\sum_{\substack{I \in \calD_0(\omega_1) \\ \ell(K) \leq \ell(I) \leq 2^{\sigma} \ell(K) \\ d(K,I) \leq d_{K,I} }}
\bla T_{\mu, \Phi}(\Delta^1_I f_\omega , E^2_{\ell(I)/2}g_\omega), \Delta^3_K h_\omega \bra_\mu.
\end{equation}
For two cubes $Q$ and $R$ in $\R^n$ let us write $Q \sim R$ if $\ell(Q) \leq \ell(R) \leq 2^\sigma \ell(Q)$ 
and $d(Q,R) \leq d_{Q,R}$. Notice the non-symmetry of this definition related to the sidelengths.

Notice that 
\begin{equation*}
\begin{split}
II = \E_\omega  \sum_{K \in \calD_0(\omega_3)}
\sum_{\substack{I \in \calD_0(\omega_1) \\ K \sim I}}
\sum_{\substack{J \in \calD_0(\omega_2) \\ \ell(J) \geq \ell(I)}} 
& \bla T_{\mu, \Phi}(\Delta^1_{I} f_\omega , \Delta^2_J g_\omega), \Delta^3_K h_\omega \bra_\mu.
\end{split}
\end{equation*}
We divide the argument into several steps.

\subsection*{Step 1}
If $\omega \in \Omega$ we agree that  $\calD_0'(\omega)$ denotes the cubes  $Q \in \calD_0(\omega)$ with $\ell(Q)<2^{u_0}$.
In this step we control
$$
II_1 = 
 \E_\omega  \sum_{K \in \calD_0(\omega_3)}
\sum_{\substack{I \in \calD_0(\omega_1) \\ K \sim I}}
\sum_{\substack{J \in \calD_0(\omega_2) \\ \ell(J) \geq \ell(I) \\ d(J,K) > d_{K,J}}} 
\bla T_{\mu, \Phi}(\Delta^1_{I} f_\omega , \Delta^2_J g_\omega), \Delta^3_K h_\omega \bra_\mu,
$$
which we further split as
$$
II_1' = 
 \E_\omega  \sum_{K \in \calD_0'(\omega_3)}
\sum_{\substack{I \in \calD_0(\omega_1) \\ K \sim I}}
\sum_{\substack{J \in \calD_0(\omega_2) \\ \ell(J) \geq \ell(I) \\ d(J,K) > d_{K,J}}} 
\bla T_{\mu, \Phi}(\Delta^1_{I} f_\omega , \Delta^2_J g_\omega), \Delta^3_K h_\omega \bra_\mu
$$
and
$$
II_1'' =
 \E_\omega  \mathop{\sum_{K \in \calD_0(\omega_3)}}_{\ell(K) = 2^{u_0}}
\sum_{\substack{I \in \calD_0(\omega_1) \\ K \sim I}}
\sum_{\substack{J \in \calD_0(\omega_2) \\ \ell(J) \geq \ell(I) \\ d(J,K) > d_{K,J}}} 
\bla T_{\mu, \Phi}(\Delta^1_{I} f_\omega , \Delta^2_J g_\omega), \Delta^3_K h_\omega \bra_\mu.
$$
This will be done by estimating uniformly in $\omega \in \Omega \times \Omega \times \Omega$.

We handle $II_1''$ first. Notice that here we have
$\ell(I) = \ell(J) = \ell(K) = 2^{u_0}$ and $d_{K,J} > 2^{u_0}$.
Therefore, it holds
\begin{align*}
|II_1''| &\le  \E_\omega  \mathop{\sum_{K \in \calD_0(\omega_3)}}_{\ell(K) = 2^{u_0}}
\sum_{\substack{I \in \calD_0(\omega_1) \\ \ell(I) = 2^{u_0}}}
\sum_{\substack{J \in \calD_0(\omega_2) \\ \ell(J) = 2^{u_0}}} 
\frac{\| \Delta^1_{I} f_\omega \|_{L^1(\mu)} \| \Delta^2_{J} g_\omega \|_{L^1(\mu)} \| \Delta^3_{K} h_\omega \|_{L^1(\mu)}}{2^{2mu_0}} \\
&\lesssim \E_\omega \frac{\| f_\omega \|_{L^1(\mu)} \|  g_\omega \|_{L^1(\mu)} \|  h_\omega \|_{L^1(\mu)}}{2^{2mu_0}} \\
&\le  \E_\omega \frac{\mu(Q_0)^{2}}{2^{2mu_0}}  \| f_\omega \|_{L^p(\mu)} \| g_\omega \|_{L^q(\mu)} \|  h_\omega \|_{L^{r'}(\mu)} \lesssim
\| f \|_{L^p(\mu)} \| g \|_{L^q(\mu)} \|  h \|_{L^{r'}(\mu)},
\end{align*}
where we used that $2^{u_0} \sim \ell(Q_0)$. This finishes the estimate for the term $II_1''$.

We now consider $II_1'$. Fix $l \ge 0$ and write
$$
\phi_{K,I,l} := \mathop{\mathop{{\sum_{J \in \calD_0(\omega_2)}}}_{d(J,K) > d_{K,J}}}_{\ell(J) = 2^l\ell(I)} \Delta^2_J g_\omega.
$$
We can agree that this function vanishes for those $l$ for which $2^l\ell(I) > 2^{u_0}$. 
Notice that $|\phi_{K,I,l}| \le |D_{2^l\ell(I)}^2 g_{\omega}| \lesssim M_{\mu, \calD_0(\omega_2)} g_{\omega}$.
Using Lemma \ref{lem:zeroave_sep}
we see that here
\begin{align*}
| \bla T_{\mu, \Phi}(\Delta^1_{I} &f_\omega,\phi_{K,I,l}), \Delta^3_K h_\omega \bra_\mu| \\
&\lesssim  2^{-l(1-\gamma)\alpha}
\int M_{\mu,m}( \Delta^1_{I} f_\omega) M_{\mu,m}(M_{\mu, \calD_0(\omega_2)} g_{\omega})  | \Delta_K^3h_\omega | \,d\mu.
\end{align*}
Therefore, we have
\begin{align*}
\sum_{l \ge 0}& \sum_{K \in \calD_0'(\omega_3)} \sum_{\substack{I \in \calD_0(\omega_1) \\ K \sim I}} | \bla T_{\mu, \Phi}(\Delta^1_{I} f_\omega,\phi_{K,I,l}), \Delta^3_K h_\omega \bra_\mu| \\
&\lesssim \sum_{K \in \calD_0(\omega_3)} \sum_{\substack{I \in \calD_0(\omega_1) \\ K \sim I}} 
\int M_{\mu,m}( \Delta^1_{I} f_\omega) M_{\mu,m}(M_{\mu, \calD_0(\omega_2)} g_{\omega})  | \Delta_K^3h_\omega | \,d\mu \\
&\lesssim \Big\| \Big( \sum_{I \in \calD_0(\omega_1)} M_{\mu,m}( \Delta^1_{I} f_\omega)^2 \Big)^{1/2}\Big\|_{L^p(\mu)} 
\| M_{\mu,m}(M_{\mu, \calD_0(\omega_2)} g_{\omega}) \|_{L^q(\mu)} \\
&  \quad \quad \quad \times \Big\| \Big( \sum_{K \in \calD_0(\omega_3)} |\Delta^3_{K} g_\omega|^2 \Big)^{1/2}\Big\|_{L^{r'}(\mu)} \\
&\lesssim \| f \|_{L^p(\mu)} \| g \|_{L^q(\mu)} \|  h \|_{L^{r'}(\mu)},
\end{align*}
where we used that given $I$ there are only $\lesssim 1$ cubes $K$ so that $K \sim I$, and the usual bounds for maximal functions and square functions.
We have now controlled $II_1'$, and so are done with Step 1:
\begin{equation}\label{eq:resStep1}
|II_1| \lesssim \| f \|_{L^p(\mu)} \| g \|_{L^q(\mu)} \|  h \|_{L^{r'}(\mu)} \le 1.
\end{equation}

\subsection*{Step 2}
In this subsection we bound
$$
II_2 = \E_\omega \sum_{K \in \calD_0(\omega_3)}   \sum_{\substack{I \in \calD_0(\omega_1) \\ K \sim I}} \sum_{\substack{J \in \calD_0(\omega_2) \\ \ell(I) \leq \ell(J) \leq2^{\sigma}\ell(K) \\ d(J,K) \leq  d_{K,J}}} 
\bla T_{\mu, \Phi}(\Delta^1_I f_\omega , \Delta^2_J g_\omega), \Delta^3_K h_\omega \bra_\mu.
$$
To ease the notation, let us define $\tilde{\frakD}(\omega)$ to be the collection of triples 
$$
(I,J,K)\in \DI \times \DII \times \DIII
$$ 
such that
$K \sim I$, $\ell(I) \leq \ell(J) \leq2^{\sigma}\ell(K)$ and  $d(J,K) \leq  d_{K,J}$. Then, define $\frakD(\omega)$ to be those triples $(I',J',K')$ such that there exists $(I,J,K) \in \tilde\frakD(\omega)$ so that $I' \in \ch(I)$, $J' \in \ch(J)$ and $K' \in \ch(K)$. 
Notice that
\begin{equation}\label{eq:diag_frakD}
II_2 = 
\E_\omega\sum_{(I,J,K) \in \frakD(\omega)} 
\bla T_{\mu, \Phi}(D^1_{I} f_\omega 1_{I}b_1 , D^2_{J} g_\omega 1_{J}b_2), D^3_K h_\omega 1_K b_3 \bra_\mu.
\end{equation}
We will decompose
$$
\bla T_{\mu, \Phi}( 1_{I}b_1, 1_{J}b_2), 1_K b_3 \bra_\mu
$$
using surgery for the triple of cubes $(I, J, K)$.

\subsubsection*{Surgery for a triple $(I,J,K) \in \frak{D}(\omega)$}
We perform surgery on $(I,J,K) \in \frak{D}(\omega)$ with a parameter $\theta > 0$. 
This follows the standard surgery for a pair of cubes (see for example \cite{NTV}) and we modify it for the present purposes. However, we have a little additional modification because we want that certain cubes in the intersection $I \cap J \cap K$  have small boundaries. (This small boundary property is merely arranged for the purposes of the proof of Corollary \ref{cor:MainThm}, see Chapter \ref{sec:syn}. It is not explicitly needed here in this triple surgery.)

Let $j(\theta)\in\Z$ be such that $$2^{-21}\theta\leq 2^{j(\theta)}<2^{-20}\theta.$$
Let $\calD(\omega_4)$ be yet another random grid in $\R^n$, independent of all other grids considered. 
Define $$\mathcal{Q}:=\{Q\in\calD(\omega_4) \colon \ell(Q)=2^{j(\theta)}\ell(K)\},$$ 
and for $x\in \R^n$, let $Q(x)$ be the unique cube in $\calQ$ that contains $x$. 
Let $Q \in \calQ$ and consider the cube $(1-\theta)Q$. By letting the small boundary parameter $t=t(\theta)$ to be large enough, there exists a cube $S_Q$,  that is concentric with $Q$, satisfies $$(1-\theta)Q \subset S_Q \subset (1-\theta/2) Q$$ and has $t$-small boundary.

We define
\begin{equation*}
\begin{split}
I_{\partial}^{J,K}(\omega_4) =  I_{\partial}&:=\{x \in I  \colon d(Q(x),\partial J)<\theta \ell(J)/2\} \\
&\cup \{x \in I  \colon d(Q(x),\partial K)<\theta \ell(K)/2\}
\\
&\cup
\{x\in I\cap J \cap K \colon  x \not \in S_{Q(x)}\}.
\end{split}
\end{equation*}
Thus points in $I_{\partial}$ belong to $I$, and are either close to the boundary of $J$, to the boundary of $K$, or to the boundary of the grid $\calQ$. The set $I_{\partial}$ depends on the cubes $J$ and $K$. However, we have
\begin{equation}\label{eq:Ibad}
\begin{split}
I_{\partial}\subset I_{\textup{bad}}^{\omega_2, \omega_3, \omega_4}
:= &\bigcup_{
\substack{
J'\in \calD(\omega_2) \\ \ell(J') \sim \ell(I)}
}
\{x \in I \colon d(x,\partial J') <\theta \ell(J')\} \\
&\cup 
\bigcup_{
\substack{
K'\in \calD(\omega_3) \\ \ell(K') \sim \ell(I)}
}
\{x \in I \colon d(x,\partial K') <\theta \ell(K')\} \\
&\cup\bigcup_{\substack{Q \in\calD(\omega_4)\\ \ell(Q) \sim 2^{j(\theta)}\ell(I)}}\{x\in I  \colon d(x,\partial Q)<\theta\ell(Q)\},
\end{split}
\end{equation}
which depends only on $I$ and $\omega_2,\omega_3, \omega_4$.

We set
\begin{equation*}
I_{\textup{sep}}^{J,K}(\omega_4) =  I_{\textup{sep}}:=I\setminus (I_{\partial}\cup (J \cap K)),
\end{equation*}
the part of $I$ strictly separated from either $J$ or $K$. Finally, we have
\begin{equation*}
  I_{\Delta}^{J,K}(\omega_4) =  I_{\Delta} := I\setminus(I_{\partial}\cup I_{\textup{sep}})
  =\bigcup_i L_{I,i},
\end{equation*}
where each $L_{I,i}$ is of the form $L_{I,i}=S_Q \cap I \cap J \cap K$ for some $Q\in \calQ$, and $\#i \lesssim_{\theta} 1$.

We have the partition
\begin{equation*}
  I=I_{\textup{sep}}\cup I_{\partial}\cup I_{\Delta} = I_{\textup{sep}}\cup I_{\partial}\cup \bigcup_i L_{I,i} ,
\end{equation*}
and in a completely analogous manner also
\begin{equation*}
  J=J_{\textup{sep}}\cup J_{\partial}\cup J_{ \Delta} =  J_{\textup{sep}}\cup J_{\partial} \cup  \bigcup_j L_{J,j}
\end{equation*}
and
\begin{equation*}
  K=K_{\textup{sep}}\cup K_{\partial}\cup K_{ \Delta} =  K_{\textup{sep}}\cup K_{\partial} \cup  \bigcup_k L_{K,k}.
\end{equation*}
A key observation is that all $L_{I,i} \subset I\cap J\cap K$ appearing in the first union are cubes (of the form $S_Q$ for $Q\in \calQ$) unless they are close to $\partial I$, and they are never close to the boundary of  $J$ or $K$. A similar statement is valid for the cubes appearing in the other unions (related to $J$ or $K$), and therefore
it follows that if a cube $L$ appears in all of the above unions (or just in two of them), then $L = S_Q$ for some $Q\in \calQ$. Such cubes also satisfy
$5L \subset I \cap J \cap K$ by construction.

We now continue with Step 2.
We fix $\omega = (\omega_1, \omega_2, \omega_3) \in \Omega \times \Omega \times \Omega$ and $\omega_4 \in \Omega$. We will find an upper bound
for
$$
II_2(\omega) := \sum_{(I,J,K) \in \frakD(\omega)} 
\bla T_{\mu, \Phi}(D^1_{I} f_\omega 1_{I}b_1 , D^2_{J} g_\omega 1_{J}b_2), D^3_K h_\omega 1_K b_3 \bra_\mu
$$
that depends on these random parameters. We will then take expectations.

Let $(I,J,K) \in \frakD(\omega)$, and let $I=I_{\text{sep}} \cup I_\partial \cup I_\Delta$ be the decomposition from above.  Since every $y \in I_{\text{sep} }$ satisfies
$$\max(d(y,J), d(y,K)) \gtrsim_{\theta} \ell(I),$$ we have
$$
\big| \bla T_{\mu, \Phi}(D^1_I f_\omega 1_{I_{sep}}b_1 , D^2_J g_\omega 1_Jb_2), D^3_K h_\omega 1_K b_3\bra_\mu \big | 
\lesssim \big| D^1_I f_\omega D^2_J g_\omega D^3_K h_\omega \big| \frac{\mu(I) \mu(J) \mu(K)}{\ell(I)^{2m}}.
$$

Notice that 
\begin{equation*}
\big| D^1_I f_\omega D^2_J g_\omega  \big| \frac{\mu(I) \mu(J)}{\ell(I)^{2m}}
\lesssim M_{\mu,m} (D^1_I f_\omega 1_I)(x) M_{\mu,m} (D^2_J g_\omega 1_J)(x)
\end{equation*}
for every $x \in K$, since the cubes are of comparable size and close to each other. 
Hence
\begin{equation*}
\begin{split}
\sum_{(I,J,K) \in \frakD(\omega)}   &\big| D^1_I f_\omega D^2_J g_\omega D^3_K h_\omega \big|\frac{\mu(I) \mu(J) \mu(K)}{\ell(I)^{2m}} \\
& \lesssim \sum_{(I,J,K) \in \frakD(\omega)} 
\int M_{\mu,m} (D^1_I f_\omega 1_I) M_{\mu,m} (D^2_J g_\omega 1_J) |D^3_K h_\omega| 1_K d \mu,
\end{split}
\end{equation*}
which is dominated via H\"older's inequality by
\begin{equation*}
\begin{split}
\Big \| &\Big( \sum_{(I,J,K) \in \frakD(\omega)} M_{\mu,m} (D^1_I f_\omega 1_I)^2 \Big)^{1/2} \Big\|_{L^p(\mu)} 
\Big\|     M_{\mu,m}(M_{\mu, \calD_0(\omega_2)} g_\omega ) \Big \|_{L^q(\mu)} \\
& \cdot \Big \| \Big( \sum_{(I,J,K) \in \frakD(\omega)}  |D^3_K h_\omega |^2 1_K \Big)^{1/2} \Big\|_{L^{r'}(\mu)}.
\end{split}
\end{equation*}
Since for every $I \in  \DI$ the number of triples such that $(I,J,K) \in \frakD(\omega)$ is uniformly bounded, and similarly for every $K \in \DIII$, this last expression is in turn dominated by 
$\|f \|_{L^p(\mu)}\|g \|_{L^q(\mu)}\|h \|_{L^{r'}(\mu)}.$ We have now controlled
$$
\sum_{(I,J,K) \in \frakD(\omega)} 
\bla T_{\mu, \Phi}(D^1_{I} f_\omega 1_{I_{\textup{sep}}}b_1 , D^2_{J} g_\omega 1_{J}b_2), D^3_K h_\omega 1_K b_3 \bra_\mu
$$
uniformly in $\omega \in \Omega \times \Omega \times \Omega$ and $\omega_4 \in \Omega$ (recall that here $I_{\textup{sep}} = I_{\textup{sep}}^{J,K}(\omega_4)$).

We will next control the sum, where $I_{\textup{sep}}$ is replaced with $I_{\partial}$. The point is to estimate this using the a priori boundedness of $T_{\mu, \Phi}$, and  apply the fact that on average (with respect to $\omega_2$, $\omega_3$ and $\omega_4$) the function $$\Big(\sum_{I} |D^1_I f|^21_{I_{\textup{bad}}^{\omega_2, \omega_3, \omega_4}} \Big)^{1/2}$$ has a small norm. The a priori boundedness can be utilised via the next lemma.
\begin{lemma}\label{lem:MZ}
For a sequence $(f_i)_{i \in \Z}$ of functions and $s \in (1, \infty)$ define
$$
\| (f_i)_{i\in \Z} \|_{L^s(\mu \colon \ell^2)}:=
\Big \| \Big( \sum_{i\in \Z} |f_i|^2 \Big)^{1/2} \Big \|_{L^s(\mu)}.
$$
If $(f_i)_{i} \subset L^p(\mu), (g_i)_i \subset L^q(\mu)$ and $(h_i)_{i} \subset L^{r'}(\mu)$ are sequences of functions, then
\begin{equation*}
\begin{split}
\Big|
\sum_i \bla &T_{\mu, \Phi} (f_i,g_i), h_i \rangle_\mu \Big | \\
& \lesssim \| T_{\mu, \Phi} \| \| (f_i)_i \|_{L^p(\mu \colon \ell^2)}
\| (g_i)_i \|_{L^q(\mu \colon \ell^2)}
\| (h_i)_i \|_{L^{r'}(\mu \colon \ell^2)}.
\end{split}
\end{equation*}
\end{lemma}

This statement is in the spirit of a classical theorem by Marcinkiewicz and Zygmund. We recall a quick proof using random signs.

\begin{proof}[Proof of Lemma \ref{lem:MZ}]
Let $(f_i)_{i} \subset L^p(\mu), (g_i)_i \subset L^q(\mu)$ and $(h_i)_{i} \subset L^{r'}(\mu)$. We may suppose that only finitely many functions in these sequences are non-zero. Let $\{\varepsilon_i\}_{i}$ and $\{ \varepsilon_i'\}_i$ be two collections of independent random signs and denote by $\E$ and $\E'$ the related expectations, correspondingly.

Taking the random signs into use we have
\begin{equation*}
\sum_i \bla T_{\mu, \Phi} (f_i,g_i), h_i \rangle_\mu 
= \E \Big \langle \sum_i \varepsilon_i T_{\mu, \Phi}(f_i,g_i), \sum_k \varepsilon_k h_k \Big \rangle,
\end{equation*} 
where further
\begin{equation*}
\sum_i \varepsilon_i T_{\mu, \Phi}(f_i,g_i)
= \E ' T_{\mu, \Phi} \Big( \sum_i \varepsilon_i \varepsilon_i' f_i, \sum_j \varepsilon_j'g_j \Big).
\end{equation*}
Thus, if we write $F_\varepsilon:= \sum_i \varepsilon_i \varepsilon_i' f_i$, $G_\varepsilon:= \sum_j \varepsilon_j' g_j$ and $H_\varepsilon := \sum_k \varepsilon_k h_k$,   then
\begin{equation*}
\begin{split}
\Big|
\sum_i \bla T_{\mu, \Phi} (f_i,g_i), h_i \rangle_\mu \Big |
&= \Big | \E \E' \big \langle T_{\mu, \Phi} ( F_\varepsilon, G_\varepsilon ), H_\varepsilon \big \rangle_\mu \Big| \\
& \leq \E \E' \| T_{\mu, \Phi}\| 
\| F_\varepsilon  \|_{L^p(\mu)}
\| G_\varepsilon  \|_{L^q(\mu)}
\| H_\varepsilon  \|_{L^{r'}(\mu)},
\end{split}
\end{equation*}
which by H\"older's inequality is at most
\begin{equation*}
 \| T_{\mu, \Phi}\| 
\big( \E \E' \| F_\varepsilon  \|_{L^p(\mu)}^p \big)^{1/p}
\big( \E \E'\| G_\varepsilon  \|_{L^q(\mu)}^q \big)^{1/q}
\big( \E \E'  \| H_\varepsilon  \|_{L^{r'}(\mu)}^{r'} \big)^{1/r'}.
\end{equation*}

Applying Khintchine inequality there holds for example that
\begin{equation*}
\begin{split}
\E \E' \| F_\varepsilon  \|_{L^p(\mu)}^p
= \E \int \E' \big| \sum_i \varepsilon_i \varepsilon_i' f_i \big|^p d \mu 
\sim \int \big(  \sum_i |f_i |^2 \big)^{p/2} d \mu.
\end{split}
\end{equation*}
Doing the same computation for $G_\varepsilon$ and $H_\varepsilon$ proves the lemma.
\end{proof}

With help of Lemma \ref{lem:MZ} we have
\begin{equation*}
\begin{split}
\Big|  &\sum_{(I,J,K) \in \frakD(\omega)}
 \bla T_{\mu, \Phi}(D^1_I f_\omega 1_{I_\partial}b_1 , D^2_J g_\omega 1_J b_2), D^3_K h_\omega 1_K b_3 \bra_\mu \Big| \\
&\lesssim \| T_{\mu, \Phi} \| \Big \| \Big( \sum_{(I,J,K) \in \frakD(\omega)}  |D^1_I f_\omega 1_{I_\partial}|^2 \Big)^{1/2} \Big\|_{L^p(\mu)} 
\Big \| \Big( \sum_{(I,J,K) \in \frakD(\omega)}  |D^2_J g_\omega 1_{J}|^2 \Big)^{1/2} \Big\|_{L^q(\mu)} \\
& \quad \quad  \quad  \! \times \Big \| \Big( \sum_{(I,J,K) \in \frakD(\omega)}  |D^3_K h_\omega 1_{K}|^2 \Big)^{1/2} \Big\|_{L^{r'}(\mu)} \\
& \lesssim \| T_{\mu, \Phi} \|
 \Big \| \Big( \sum_{I \in \calD_0'(\omega_1)} |D^1_I f 1_{I_{\textup{bad}}^{\omega_2, \omega_3, \omega_4}}|^2 \Big)^{1/2} \Big\|_{L^p(\mu)}
\| g \|_{L^q(\mu)} \| h \|_{L^{r'}(\mu)}.
\end{split}
\end{equation*}
This is of the right form so that we can take averages at the end. Therefore, we now move on to the term where
$I_{\partial}$ is replaced with $I_{\Delta}$.

We have arrived at the term
$$
\sum_{(I,J,K) \in \frakD(\omega)} 
\bla T_{\mu, \Phi}(D^1_{I} f_\omega 1_{I_{\Delta}}b_1 , D^2_{J} g_\omega 1_{J}b_2), D^3_K h_\omega 1_K b_3 \bra_\mu.
$$
Continuing in the natural way, we can dominate $II_2(\omega)$
with the sum of 
$$C(\theta)\|f\|_{L^p(\mu)} \|g\|_{L^q(\mu)} \|h\|_{L^{r'}(\mu)},$$
$$
\| T_{\mu, \Phi} \|
 \Big \| \Big( \sum_{I \in \calD_0'(\omega_1)} |D^1_I f 1_{I_{\textup{bad}}^{\omega_2, \omega_3, \omega_4}}|^2 \Big)^{1/2} \Big\|_{L^p(\mu)}
\| g \|_{L^q(\mu)} \| h \|_{L^{r'}(\mu)},
$$
two corresponding terms where the ``bad'' square function appears in $g$ or $h$, and
\begin{equation}\label{eq:intersection}
\Big| \sum_{(I,J,K) \in \frakD(\omega)} 
\bla T_{\mu, \Phi}(D^1_{I} f_\omega 1_{I_{\Delta}}b_1 , D^2_{J} g_\omega 1_{J_{\Delta}}b_2), D^3_K h_\omega 1_{K_{\Delta}} b_3 \bra_\mu \Big|.
\end{equation}

Suppose $(I,J,K) \in \frakD(\omega)$. We further split $I_\Delta= \bigcup_iL_{I,i}$, $J_\Delta= \bigcup_jL_{J,j}$ and $K_\Delta= \bigcup_k L_{K,k}$. If
$(i,j,k)$ is such that $L_{I,i}$, $L_{J,j}$ and $L_{K,k}$ are not all equal, then separation between two of these is $\gtrsim_\theta \ell(I)$, and we have
\begin{equation*}
\big |\langle T_{\mu, \Phi}(1_{L_{I,i}},1_{L_{j,J}}),1_{L_{K,k}} \rangle_\mu \big|
\lesssim_\theta \frac{\mu(L_{I,i})\mu(L_{J,j})\mu(L_{K,k})}{\ell(I)^{2m}}
\lesssim \mu(I \cap J \cap K).
\end{equation*}
On the other hand, if $(i,j,k)$ is such that $L_{I,i}=L_{J,j}=L_{K,k} =: L$, then by the observations made during the construction of the surgery, we have that $L$ is a cube (with small boundary) and $5L \subset I \cap J \cap K$. In this case we can use the weak boundedness assumption to have
$$
\big |\langle T_{\mu, \Phi}(1_{L},1_{L}),1_{L} \rangle_\mu \big|
\lesssim \mu(5L)
\leq \mu(I \cap J \cap K).
$$
Since there are only $\lesssim_\theta 1$ cubes in the splittings of $I_\Delta$, $J_\Delta$ and $K_\Delta$, we have shown that
\begin{equation}\label{eq:all_Delta}
\big |\langle T_{\mu, \Phi}(1_{I_\Delta},1_{J_\Delta}),1_{K_\Delta} \rangle_\mu \big|
\lesssim_\theta \mu(I \cap J \cap K).
\end{equation}

Using \eqref{eq:all_Delta}, we have
\begin{equation*}
\begin{split}
\eqref{eq:intersection}
&\lesssim_\theta \sum_{(I,J,K) \in \frakD(\omega)} \big | D^1_{I} f_\omega  D^2_{J} g_\omega  D^3_K h_\omega \big| \mu(I \cap J \cap K) \\
& \lesssim \Big \| \Big( \sum_{(I,J,K) \in \frakD(\omega)}  |D^1_I f_\omega|^2 1_{I \cap J \cap K} \Big)^{1/2} \Big\|_{L^p(\mu)}
\Big \| \sup_{(I,J,K) \in \frakD(\omega)}  |D^2_J g_\omega| 1_{I \cap J \cap K} \Big\|_{L^q(\mu)} \\
 &\ \ \! \times \Big \| \Big( \sum_{(I,J,K) \in \frakD(\omega)}  |D^3_K h_\omega|^2 1_{I \cap J \cap K} \Big)^{1/2} \Big\|_{L^{r'}(\mu)} \\
&\lesssim \|f \|_{L^p(\mu)}\|g \|_{L^q(\mu)}\|h \|_{L^{r'}(\mu)},
\end{split}
\end{equation*}
where we used that $$\sup_{(I,J,K) \in \frakD(\omega)}  |D^2_J g_\omega| 1_{I \cap J \cap K} \lesssim M_{\mu, \calD_0(\omega_2)} g_\omega.$$

Combining everything, and recalling that the functions have norm at most $1$, we have shown that
\begin{align*}
|II_2| &= |\E_{\omega} II_2(\omega)| = |\E_{\omega_4} \E_{\omega} II_2(\omega)|  \\
& \lesssim C(\theta) + \| T_{\mu, \Phi} \| \E_{\omega_4} \E_{\omega} \Big[ 
\Big \| \Big( \sum_{I \in \calD_0'(\omega_1)} |D^1_I f|^2 1_{I_{\textup{bad}}^{\omega_2, \omega_3, \omega_4}} \Big)^{1/2} \Big\|_{L^p(\mu)} \\
&\quad \quad \quad \quad \quad \quad  \quad \quad \quad \quad + 
\Big \| \Big( \sum_{J \in \calD_0'(\omega_2)} |D^2_J g|^2 1_{J_{\textup{bad}}^{\omega_1, \omega_3, \omega_4}} \Big)^{1/2} \Big\|_{L^q(\mu)} \\
& \quad \quad \quad \quad \quad \quad \quad \quad \quad \quad +
\Big \| \Big( \sum_{K \in \calD_0'(\omega_3)} |D^3_K h|^2 1_{K_{\textup{bad}}^{\omega_1, \omega_2, \omega_4}} \Big)^{1/2} \Big\|_{L^{r'}(\mu)} \Big].
\end{align*}

The above averages of bad square functions can be dominated with $c(\theta)$, where $\lim_{\theta \to 0} c(\theta) = 0$, as stated in Lemma \ref{lem:badSF}.
In $L^2$ this is easy and follows from the work of Nazarov--Treil--Volberg.
In $L^p$ such bounds have usually been obtained using some fairly heavy machinery (an improved contraction principle), see \cite{Hy:vec}.
However, the interpolation technique used in \cite{LV} to control the good and bad parts of functions (as in Lemma \ref{lem:ave_bad_norm})
seems to lend an easier proof here also. For some reason this was not already used in \cite{LV}, however. We record this simpler proof idea below.
\begin{lemma}\label{lem:badSF}
The estimate
\begin{equation}\label{eq:interp_surg}
\E_{\omega_2, \omega_3, \omega_4} \Big \| \Big( \sum_{I \in \calD_0'(\omega_1)} |D^1_I f|^2 1_{I_{\textup{bad}}^{\omega_2, \omega_3, \omega_4}} \Big)^{1/2} \Big\|_{L^p(\mu)}
\lesssim c(\theta,p) \| f \|_{L^p(\mu)}, \quad f \in L^p(\mu),
\end{equation}
holds. Here $\lim_{\theta \to 0} c(\theta,p) = 0$.
\end{lemma}

\begin{proof}
For $f \in L^1_{\text{loc}}(\mu)$, define the operator $\calS_\calB$ by
$$
\calS_\calB f(x, \omega_2, \omega_3, \omega_4) := \Big( \sum_{I \in \calD_0'(\omega_1)} |D^1_I f|^2 1_{I_{\textup{bad}}^{\omega_2, \omega_3, \omega_4}}(x) \Big)^{1/2}.
$$
For a $\mu \times \bbP \times \bbP \times \bbP$-measurable function $\varphi$ write 
$$
\| \varphi \|_{L^p(\mu \times  \bbP \times \bbP \times \bbP)}
:= \Big(\E_{\omega_2}\E_{\omega_3}\E_{\omega_4} \int |\varphi(x,\omega_2,\omega_3, \omega_4) |^p d \mu \Big)^{1/p}.
$$

We will show via interpolation that
$$
 \|\calS_\calB f \|_{L^p(\mu \times  \bbP \times \bbP \times \bbP)} \lesssim c(\theta,p) \| f \|_{L^p(\mu)}.
$$
This concludes the proof of the lemma, since the left hand side of \eqref{eq:interp_surg} is at most  $ \|\calS_\calB f \|_{L^p(\mu \times  \bbP \times \bbP \times \bbP)}$  by H\"older's inequality.

First, notice that trivially
$$
\Big( \sum_{I \in \calD_0'(\omega_1)} |D^1_I f|^2 1_{I_{\textup{bad}}^{\omega_2, \omega_3, \omega_4}}(x) \Big)^{1/2}
\le \Big( \sum_{I \in \calD_0'(\omega_1)} |D^1_I f|^2 1_{I}(x) \Big)^{1/2},
$$
and so  $$ \|\calS_\calB f \|_{L^p(\mu \times  \bbP \times \bbP \times \bbP)} \lesssim_p \| f \|_{L^p(\mu)}, \qquad p \in (1,\infty).$$

Considering the $L^2$ estimate, we have
$$
\|\calS_\calB f \|_{L^2(\mu \times  \bbP \times \bbP \times \bbP)}^2
=\sum_{I \in \calD_0'(\omega_1)} |D^1_I f|^2 \E_{\omega_2}\E_{\omega_3}\E_{\omega_4}\mu( I_{\textup{bad}}^{\omega_2, \omega_3, \omega_4} ),
$$
where it is standard that the average above is $\le c(\theta) \mu(I)$ (see \cite{NTV}).
Using this it follows that
$$
\|\calS_\calB f \|_{L^2(\mu \times  \bbP \times \bbP \times \bbP)}^2 \lesssim c(\theta) \|f\|_{L^2(\mu)}^2.
$$
\end{proof}
In synthesis, in Step 2 we have shown that
\begin{equation}\label{eq:resStep2}
|II_2| \lesssim C(\theta) + c(\theta)\| T_{\mu, \Phi}\|.
\end{equation}
We will not fix the parameter $\theta$ yet, since more surgeries will appear. Therefore, we are ready to move to Step 3.

\subsection*{Step 3}
What is left after steps 1 and 2 is

\begin{equation*}
II_3 = \E_\omega \sum_{K \in \calD_0(\omega_3)}   \sum_{\substack{I \in \calD_0(\omega_1) \\ K \sim I}} \sum_{\substack{J \in \calD_0(\omega_2) \\ \ell(J) > 2^{\sigma}\ell(K) \\ d(J,K) \le  d_{K,J}}} 
\bla T_{\mu, \Phi}(\Delta^1_I f_\omega , \Delta^2_J g_\omega), \Delta^3_K h_\omega \bra_\mu.
\end{equation*}
Recall the set  $\calD_h(\omega_3)$ and the notation $J_{K, l} \in \calD_0(\omega_2)$ for those $K \in \calD_h(\omega_3)$ and $l\in \Z$ that satisfy $2^{\sigma}\ell(K) \leq 2^l\ell(K) \leq 2^{u_0}$.
The existence of these cubes $J_{K,l}$ has been justified in the paraproduct section. We can now write
$$
II_3 = \E_\omega \sum_{K \in \calD_h(\omega_3)} \sum_{\substack{I \in \calD_0(\omega_1) \\ K \sim I}}   \sum_{l\colon\, 2^{\sigma}\ell(K) < 2^l\ell(K) \le 2^{u_0}}  
\bla T_{\mu, \Phi}(\Delta^1_I f_\omega , \Delta^2_{J_{K,l}} g_\omega), \Delta^3_K h_\omega \bra_\mu.
$$
Next, we perform the standard splitting
$$
\Delta^2_{J_{K,l}} g_\omega
= 1_{J_{K,l-1}^c}\Delta^2_{J_{K,l}} g_\omega
- (D^2_{J_{K,l-1}} g_\omega)1_{J_{K,l-1}^c}  b_2
+ (D^2_{J_{K,l-1}} g_\omega)b_2.
$$
Recalling that
$$
\sum_{l\colon\, 2^{\sigma}\ell(K) < 2^l\ell(K) \le 2^{u_0}}
D^2_{J_{K,l-1}}g_{\omega}
= E^2_{J_{K,\sigma}}g_\omega
$$
we get, after replacing $\Delta^2_{J_{K,l}} g_\omega$ with $(D^2_{J_{K,l-1}} g_\omega)b_2$ in $II_3$, the main term
$$
II_3' = \E_\omega \sum_{K \in \calD_h(\omega_3)} \sum_{\substack{I \in \calD_0(\omega_1) \\ K \sim I}} 
\bla T_{\mu, \Phi}(\Delta^1_I f_\omega , E^2_{J_{K,\sigma}}g_\omega b_2) , \Delta^3_K h_\omega \bra_\mu.
$$
However, we first need to deal with the two error terms:
$$
II_3'' =  \E_\omega \sum_{K \in \calD_h(\omega_3)} \sum_{\substack{I \in \calD_0(\omega_1) \\ K \sim I}}   \sum_{l=\sigma+1}^{u_0-\log_2 \ell(K)}
\bla T_{\mu, \Phi}(\Delta^1_I f_\omega ,1_{J_{K,l-1}^c}\Delta^2_{J_{K,l}} g_\omega), \Delta^3_K h_\omega \bra_\mu
$$
and
$$
II_3''' =  \E_\omega \sum_{K \in \calD_h(\omega_3)} \sum_{\substack{I \in \calD_0(\omega_1) \\ K \sim I}}  \sum_{l=\sigma+1}^{u_0-\log_2 \ell(K)}
\bla T_{\mu, \Phi}(\Delta^1_I f_\omega , (D^2_{J_{K,l-1}} g_\omega)1_{J_{K,l-1}^c}  b_2) , \Delta^3_K h_\omega \bra_\mu.
$$
These are very brief to deal with as in the paraproduct section. Indeed, applying Lemma \ref{lem:zeroave_sep} we have
\begin{equation*}
\begin{split}
\Big | \big \langle &T_{\mu, \Phi}(\Delta^1_I f_\omega , 1_{J_{K,l-1}^c}\Delta^2_{J_{K,l}} g_\omega), \Delta^3_K h_\omega \big \rangle_\mu \Big | \\
&\lesssim 2^{-l(1-\gamma)\alpha} \int M_{\mu, m}(\Delta^1_I f_\omega) M_{\mu, m}  (\Delta^2_{J_{K,l}} g_\omega) | \Delta^3_K h_\omega| d \mu  \\
& \lesssim 2^{-l(1-\gamma)\alpha}\int M_{\mu, m}(\Delta^1_I f_\omega) M_{\mu,m} (M_{\mu, \calD_0(\omega_2)} g_\omega) | \Delta^3_K h_\omega| d \mu
\end{split}
\end{equation*} 
and similarly
\begin{equation*}
\begin{split}
\Big | \big \langle &T_{\mu, \Phi}(\Delta^1_I f_\omega , (D^2_{J_{K,l-1}} g_\omega)1_{J_{K,l-1}^c}  b_2), \Delta^3_K h_\omega \big \rangle_\mu \Big | \\
& \lesssim 2^{-l(1-\gamma)\alpha}\int M_{\mu,m}(\Delta^1_I f_\omega) |D^2_{J_{K,l-1}} g_\omega | | \Delta^3_K h_\omega| d \mu \\
&\lesssim 2^{-l(1-\gamma)\alpha}\int M_{\mu,m}(\Delta^1_I f_\omega) M_{\mu, \calD_0(\omega_2)} g_\omega | \Delta^3_K h_\omega| d \mu.
\end{split}
\end{equation*}
These are summable to the right bound over the appropriate $K,I$ and $l$, and so $|II_3''| + |II_3'''| \lesssim 1$.

We continue with the main term $II_3'$.
Denote by $\tilde{\scrD}(\omega)$ those pairs $(I,K)$ where $K \in \calD_h(\omega_3)$, $I \in \DI$ and $K \sim I$. Then we define
$$
\scrD(\omega):= \big\{ (I',K') \colon I' \in \ch (I) \text{ and } K' \in \ch(K) \text{ for some } (I,K) \in \tilde{\scrD}(\omega)\big\}.
$$
We can write
$$
II_3' = \E_\omega \sum_{(I,K) \in \scrD(\omega)}
\bla T_{\mu, \Phi}(D^1_I f_\omega 1_I b_1 ,  (E^2_{J_{K^{(1)}, \sigma}}g_\omega) b_2), D^3_K h_\omega  1_Kb_3 \bra_\mu.
$$

We shall perform another surgery argument, this time with the pairs $(I,K) \in \scrD(\omega)$. This is standard, but we again make the little modification by which the cubes in the intersection $I \cap K$ have small boundaries.
(This time we actually use the small boundaries also in this proof, but they are also needed for the proof of Corollary \ref{cor:MainThm}.)

\subsubsection*{Surgery for a pair of cubes $(I, K) \in \scrD(\omega)$}
As before, we set $j(\theta)\in\Z$ so that $2^{-21}\theta\leq 2^{j(\theta)}<2^{-20}\theta$,
$$\calQ:=\{Q\in\calD(\omega_4)\colon \ell(Q)=2^{j(\theta)}\ell(K)\}$$ and let $Q(x)$ be the unique cube in $\calQ$ that contains $x$.
Let $Q \in \calQ$ and consider the cube $(1-\theta)Q$. By letting the small boundary parameter $t=t(\theta)$ to be large enough, there exists a cube $S_Q$,  that is concentric with $Q$, satisfies $$(1-\theta)Q \subset S_Q \subset (1-\theta/2) Q$$ and has $t$-small boundary. 

Now, we define the sets $I_{\partial}$, $I_{\textup{sep}}$, $I_{\Delta}$ and $I_{\textup{bad}}$: 
\begin{equation*}
I_{\partial}^K(\omega_4) = I_{\partial}:=\{x\in I \colon d(Q(x),\partial K)<\theta\ell(K)/2\} 
 \cup
 \{x\in I\cap K \colon  x \not \in  S_{Q(x)}\};
\end{equation*}
\begin{equation*}
\begin{split}
I_{\partial}\subset I_{\textup{bad}}^{\omega_3, \omega_4}
:= &\bigcup_{\substack{K'\in \calD(\omega_3)\\ \ell(K') \sim \ell(I)}}\{x \in I \colon d(x,\partial K')<\theta\ell(K')\} \\
&\cup
\bigcup_{\substack{Q \in \calD(\omega_4)\\ \ell(Q) \sim 2^{j(\theta)}\ell(I)}} \{x \in I \colon d(x,\partial Q)<\theta\ell(Q)\};
\end{split}
\end{equation*}
\begin{equation*}
  I_{\textup{sep}}^K(\omega_4) =  I_{\textup{sep}}:=I\setminus (I_{\partial} \cup K)
\end{equation*}
and
\begin{equation*}
I_{\Delta}^K(\omega_4) = I_{\Delta} := I\setminus(I_{\partial}\cup I_{\textup{sep}})
=\bigcup_{i} L_{I,i}.
\end{equation*}
Here each $L_{I,i}$ is of the form $S_Q \cap I \cap K$ for some $Q\in  \calQ$,  $\# i \lesssim_{\theta} 1$, and
$L_i=S_Q$ unless it is close to the boundary of $I$. As before, the same splitting is performed starting from
the cube $K$, and one observes that if $L = L_{I,i} = L_{K,k}$, then $L$ is a cube and $5L \subset I \cap K$.

We now continue with $II_3'$. We write $$II_3' = E_{\omega} II_3'(\omega) = E_{\omega_4} E_{\omega} II_3'(\omega),$$ where
$$
II_3'(\omega) = \sum_{(I,K) \in \scrD(\omega)}
\bla T_{\mu, \Phi}(D^1_I f_\omega 1_I b_1 ,  (E^2_{J_{K^{(1)}, \sigma}}g_\omega) b_2), D^3_K h_\omega  1_Kb_3 \bra_\mu.
$$
We fix $\omega = (\omega_1, \omega_2, \omega_3) \in \Omega \times \Omega \times \Omega$ and $\omega_4 \in \Omega$,
and estimate $II_3'(\omega)$ with a bound depending on these random parameters.

First, we replace the cubes $I$ with the separated parts $I_{\text{sep}}$. 
Let $(I,K) \in \scrD(\omega)$. Then
\begin{equation*}
\begin{split}
\big |\bla T_{\mu, \Phi}&(D^1_I f_\omega 1_{I^K_{\text{sep}}} b_1,  (E^2_{J_{K^{(1)}, \sigma}}g_\omega) b_2), D^3_K h_\omega 1_K b_3 \bra_\mu \big | \\
& \lesssim_{\theta} \frac{  | D^1_I f_\omega |   \mu(I) |E^2_{J_{K^{(1)}, \sigma}}g_\omega | | D^3_K h_\omega | \mu(K)}{\ell(I)^m},
\end{split}
\end{equation*}
where separation and the fact that
$$
\int \frac{d\mu(z)}{(\ell(I) + |x-z|)^{2m}} \lesssim \ell(I)^{-m}
$$
were used.
Because $K^{(1)} \sim I^{(1)}$, there holds for every $x \in K$ that
$$
\frac{| D^1_I f_\omega|\mu(I)}{\ell(I)^m} \lesssim M_{\mu,m}(D^1_If_\omega 1_I)(x).
$$
Since also $|E^2_{J_{K^{(1)}, \sigma}}g_\omega | \lesssim M_{\mu, \calD_0(\omega_2)} g_\omega(x)$ for every $x \in K$ we have
$$
\frac{  | D^1_I f_\omega |   \mu(I) |E^2_{J_{K^{(1)}, \sigma}}g_\omega | | D^3_K h_\omega | \mu(K)}{\ell(I)^m}
\lesssim \int M_{\mu,m}(D^1_I f_\omega 1_I) M_{\mu, \calD_0(\omega_2)} g_\omega |D^3_Kh_\omega|1_K d \mu.
$$
This can be summed over $(I,K) \in \scrD(\omega)$ as we have seen many times.

Next, we look at the term that arises when in $II_3'$ the cubes $I$ are replaced with $I_\partial$. As above with the other surgery argument, the point is to use the a priori boundedness of $T_{\mu, \Phi}$. To this end, we need the standard principal stopping cubes for the function $g_\omega$ in the grid $\calD_0(\omega_2)$. These are constructed
precisely as in the proof of Lemma \ref{lem:parest}. We denote this collection by $\calG = \calG(\omega)$.

For every $K  \in \calD_0(\omega_3)$ that appears in the pairs $(I,K) \in \scrD(\omega)$, define
$$
a_K:= \frac{E^2_{J_{K^{(1)}, \sigma}}g_\omega}{\langle | g_\omega | \rangle^\mu_{\pi_\calG (J_{K^{(1)}, \sigma})}}.
$$
By the stopping condition (and accretivity) we know that $|a_K| \lesssim 1$.
If $G \in \calG(\omega)$, define $\scrD_G(\omega)$ to be those pairs $(I,K) \in \scrD(\omega)$ such that $$\pi_\calG J_{K^{(1)},\sigma} =G.$$
Using the collections $\scrD_G(\omega)$ and the numbers $a_K$ we  have 
\begin{equation}\label{eq:diag_stop}
\begin{split}
\sum_{(I,K) \in \scrD(\omega)}
&\bla T_{\mu, \Phi}(D^1_I f_\omega  1_{I_\partial} b_1 , E^2_{J_{K^{(1)} , \sigma}} g_{\omega} b_2) , D^3_K h_\omega  1_K b_3 \bra_\mu \\
 &=\sum_{G \in \calG(\omega)} \langle | g_\omega | \rangle^\mu_{G} 
\sum_{(I,K) \in \scrD_G(\omega)}   
 \bla T_{\mu, \Phi}(D^1_I f_\omega  1_{I_\partial} b_1,  b_2), a_KD^3_K h_\omega  1_K b_3 \bra_\mu.
\end{split}
\end{equation}

Let $G \in \calG(\omega)$ and $(I,K) \in \scrD_{G}(\omega)$. Using the goodness of $K^{(1)}$ and the size estimate of the kernel we have
$$
\big|\bla T_{\mu, \Phi}( 1_{I_\partial} b_1, 1_{G^c}b_2),  1_K b_3\bra_\mu \big|
\lesssim \frac{\mu(I) \mu(K)}{(\ell(K)^\gamma \ell(G)^{1-\gamma})^m} 
\leq \frac{\mu(I) \mu(K)}{\ell(K)^m},
$$
whence, using as before that $K^{(1)} \sim I^{(1)}$  and  $K \subset G$, we have
\begin{equation*}
\begin{split}
\langle | g_\omega | \rangle^\mu_{G}\big|\bla &T_{\mu, \Phi}(D^1_I f_\omega 1_{I_\partial} b_1, 1_{G^c}b_2), a_K D^3_K h_\omega 1_K  b_3 \bra_\mu \big| \\
& \lesssim \int M_{\mu, m}(D^1_I f_\omega 1_I) M_{\mu, \calD_0(\omega_2)}g_\omega |D^3_K h_\omega| 1_K \,d \mu.
\end{split}
\end{equation*}
Since  this is summable over $(I,K) \in \scrD(\omega)$, to control the right hand side of \eqref{eq:diag_stop} it suffices to estimate
\begin{equation}\label{eq:diag_para_boundary_stop}
\begin{split}
\sum_{G \in \calG(\omega)} \langle | g_\omega | \rangle^\mu_{G} 
\sum_{(I,K) \in \scrD_G(\omega)}  
& \bla T_{\mu, \Phi}(D^1_I f_\omega  1_{I_\partial} b_1, 1_G b_2), a_KD^3_K h_\omega  1_K b_3\bra_\mu.
\end{split}
\end{equation}

Let $G \in \calG(\omega)$.  We enumerate $\scrD_G(\omega)$ by writing
$$
\scrD_G(\omega)=\big\{(I,K)^G_i\}_i.
$$ 
Let $(\varepsilon_i)_i$ be an independent sequence of random signs, and write $\E_\varepsilon$ for the corresponding expectation.
For the moment define the shorthands
$$
f_{(I,K)}:=D^1_I f_\omega 1_{I_\partial} \quad  
\text{ and } 
\quad  h_{(I,K)}:= a_K D^3_K h_\omega 1_K, \quad (I,K) \in \scrD_G(\omega).
$$
Applying this notation and the random signs there holds
\begin{equation*}
\begin{split}
\sum_{(I,K) \in \scrD_G(\omega)}  
 &\bla T_{\mu, \Phi}(D^1_I f_\omega  1_{I_\partial} b_1, 1_G b_2), a_KD^3_K h_\omega  1_K b_3\bra_\mu \\
 &=\sum_i \bla T_{\mu, \Phi} ( f_{(I,K)^G_i},1_Gb_2), h_{(I,K)^G_i} \bra_\mu \\
& = \E_\varepsilon \bla T_{\mu, \Phi} (\psi_{\varepsilon,G},1_Gb_2), \eta_{\varepsilon,G} \bra_\mu,
\end{split}
\end{equation*}
where
$$
\psi_{\varepsilon,G}:= \sum_i \varepsilon_i f_{(I,K)^G_i}
$$
and 
$$
\eta_{\varepsilon,G}:= \sum_i \varepsilon_i h_{(I,K)^G_i}.
$$
Using this we can write \eqref{eq:diag_para_boundary_stop} as
\begin{equation*}
\sum_{G \in \calG(\omega)} \langle | g_\omega | \rangle^\mu_{G} 
\E_\varepsilon \bla T_{\mu, \Phi} (\psi_{\varepsilon,G},1_Gb_2), \eta_{\varepsilon,G} \bra_\mu.
\end{equation*}

Now, Lemma \ref{lem:MZ} gives
\begin{equation*}
\begin{split}
|\eqref{eq:diag_para_boundary_stop}|
& \leq \E_\varepsilon \Big | \sum_{G \in \calG(\omega)}   
\bla T_{\mu, \Phi} (\psi_{\varepsilon,G}, \langle | g_\omega | \rangle^\mu_{G}1_Gb_2), \eta_{\varepsilon,G} \bra_\mu \Big | \\
& \lesssim  \E_\varepsilon  \| T_{\mu, \Phi} \|
\Big \| \Big( \sum_{G \in \calG(\omega)} | \psi_{\varepsilon,G}|^2 \Big)^{1/2 } \Big\|_{L^p(\mu)}
\Big \| \Big( \sum_{G \in \calG(\omega)} \big( \langle | g_\omega |\rangle_G ^\mu|\big)^21_G \Big)^{1/2 } \Big\|_{L^q(\mu)} \\
& \quad \quad  \quad   \quad \times
\Big \| \Big( \sum_{G \in \calG(\omega)} | \eta_{\varepsilon,G}|^2 \Big)^{1/2 } \Big\|_{L^{r'}(\mu)}.
\end{split}
\end{equation*}

Carleson embedding theorem implies that the middle factor related to $g_\omega$  is dominated by $\| g \|_{L^q(\mu)}$.  
Concerning the other two factors, we have
\begin{equation*}
\begin{split}
\E_\varepsilon\Big \| &\Big( \sum_{G \in \calG(\omega)} | \psi_{\varepsilon,G}|^2 \Big)^{1/2 } \Big\|_{L^p(\mu)}
\Big \| \Big( \sum_{G \in \calG(\omega)} | \eta_{\varepsilon,G}|^2 \Big)^{1/2 } \Big\|_{L^{r'}(\mu)} \\
&\leq \Big(\E_\varepsilon  \Big \| \Big( \sum_{G \in \calG(\omega)} | \psi_{\varepsilon,G}|^2 \Big)^{1/2 } \Big\|_{L^p(\mu)}^p \Big)^{1/p}
\Big(\E_\varepsilon  \Big \| \Big( \sum_{G \in \calG(\omega)} | \eta_{\varepsilon,G}|^2 \Big)^{1/2 } \Big\|_{L^{r'}(\mu)}^{r'} \Big)^{1/{r'}},
\end{split}
\end{equation*}
and
\begin{equation}\label{eq:use_KK}
\begin{split}
\E_\varepsilon  \Big \| \Big( \sum_{G \in \calG(\omega)} | \psi_{\varepsilon,G}|^2 \Big)^{1/2 } \Big\|_{L^p(\mu)}^p
&= \E_\varepsilon \Big \| \Big( \sum_{G \in \calG(\omega)} \Big| \sum_i \varepsilon_i f_{(I,K)^G_i}\Big|^2 \Big)^{1/2 } \Big\|_{L^p(\mu)}^p \\
& \sim \Big \| \Big( \sum_{G \in \calG(\omega)}  \sum_i  |f_{(I,K)^G_i}|^2 \Big)^{1/2 } \Big\|_{L^p(\mu)}^p.
\end{split}
\end{equation}
Here we used an $\ell^2$-valued version of the Kahane--Khintchine inequality -- see for example Theorem 6.2.4 in \cite{HNVW2}.
Writing out the definition of the functions $f_{(I,K)}$ the right hand side of \eqref{eq:use_KK} can be estimated up by
\begin{equation*}
\Big \| \Big( \sum_{I \in \calD'_0(\omega_1)}  |D^1_I f_\omega |^2 1_{I_{\textup{bad}}^{\omega_3, \omega_4}} \Big)^{1/2 } \Big\|_{L^p(\mu)}^p.
\end{equation*}
The corresponding estimate holds for the functions $\eta_{\varepsilon,G}$, and the resulting $L^{r'}(\mu)$ norm of the square sum of the terms $|D^3_K h_\omega |1_K$ is $\lesssim \| h \|_{L^{r'}(\mu)}$.

Putting the above steps together, we have shown that
$$
|\eqref{eq:diag_para_boundary_stop}|
\lesssim  \| T_{\mu, \Phi} \| \Big \| \Big( \sum_{I \in \calD'_0(\omega_1)}  |D^1_I f_\omega |^2 1_{I_{\textup{bad}}^{\omega_3, \omega_4}} \Big)^{1/2 } \Big\|_{L^p(\mu)}
\| g \|_{L^q(\mu)} \| h \|_{L^{r'}(\mu)}.
$$
We are left with the term
\begin{equation*}
 \sum_{(I,K) \in \scrD(\omega)}
\bla T_{\mu, \Phi}(D^1_I f_\omega 1_{I_{\Delta}} b_1 ,  (E^2_{J_{K^{(1)}, \sigma}}g_\omega) b_2), D^3_K h_\omega  1_Kb_3 \bra_\mu.
\end{equation*}

Repeating the arguments with $K$, we reduce to
\begin{equation}\label{eq:diaggg}
 \sum_{(I,K) \in \scrD(\omega)}
\bla T_{\mu, \Phi}(D^1_I f_\omega 1_{I_{\Delta}} b_1 ,  (E^2_{J_{K^{(1)}, \sigma}}g_\omega) b_2), D^3_K h_\omega  1_{K_{\Delta}}b_3 \bra_\mu.
\end{equation}
Let $(I,K) \in \scrD(\omega)$. We use the decompositions $I_{\Delta}= \bigcup_i L_{I,i}$ and $K_\Delta=\bigcup_kL_{K,k}$.
If $L_{I,i} \not = L_{K,k}$, then the separation between these is $\gtrsim_{\theta} \ell(I)$, whence
\begin{align*}
\big |\bla T_{\mu, \Phi}(D^1_I f_\omega 1_{L_{I,i}} &b_1 ,  (E^2_{J_{K^{(1)}, \sigma}}g_\omega) b_2), D^3_K h_\omega  1_{L_{K,k}}b_3 \bra_\mu \big | \\
&\lesssim_{\theta} \big | D^1_I f_\omega E^2_{J_{K^{(1)}, \sigma}} g_\omega D^3_K h_\omega \big| \frac{\mu(L_{I,i}) \mu(L_{K,k})}{\ell(I)^m} \\
& \lesssim \big | D^1_I f_\omega E^2_{J_{K^{(1)}, \sigma}} g_\omega D^3_K h_\omega \big| \mu(I \cap K) \\
& \lesssim \int | D^1_I f_\omega |1_I M_{\mu, \calD_0(\omega_2)} g_\omega | D^3_Kh_\omega| 1_K \, d \mu.
\end{align*}

If $L_{I,i}=L_{K,k}=L$, then $L$ is a cube that has a $t(\theta)$-small boundary and $5L \subset I \cap K$. Using separation, we have
\begin{align*}
\big |\bla T_{\mu, \Phi}(D^1_I f_\omega 1_{L }&b_1 ,  (E^2_{J_{K^{(1)}, \sigma}}g_\omega) 1_{(2L)^c} b_2), D^3_K h_\omega  1_{L}b_3 \bra_\mu \big | \\
& \lesssim \big | D^1_I f_\omega E^2_{J_{K^{(1)}, \sigma}} g_\omega D^3_K h_\omega \big| \frac{\mu(L) \mu(L)}{\ell(L)^m},
\end{align*}
which produces the same bound as the estimate before this one.

Since $L$ has a $t(\theta)$-small boundary, we have (similarly as for example in \eqref{eq:SmallBoundComp}) that
\begin{align*}
\big |\bla T_{\mu, \Phi}(D^1_I f_\omega 1_{L }&b_1 ,  (E^2_{J_{K^{(1)}, \sigma}}g_\omega) 1_{2L \setminus L} b_2), D^3_K h_\omega  1_{L}b_3 \bra_\mu \big | \\
& \lesssim_{\theta} \big | D^1_I f_\omega E^2_{J_{K^{(1)}, \sigma}} g_\omega D^3_K h_\omega \big| \mu(2L). \\
\end{align*}
This again leads to the same estimate as above.

Finally, weak boundedness gives
\begin{align*}
\big |\bla T_{\mu, \Phi}(D^1_I f_\omega 1_{L }&b_1 ,  (E^2_{J_{K^{(1)}, \sigma}}g_\omega) 1_{L} b_2), D^3_K h_\omega  1_{L}b_3 \bra_\mu \big | \\
& \lesssim \big | D^1_I f_\omega E^2_{J_{K^{(1)}, \sigma}} g_\omega D^3_K h_\omega| \mu(5L),
\end{align*}
which again yields the same estimate.

Since there are $\lesssim_{\theta} 1$ cubes $L_{I,i}$ and $L_{K,k}$, we have shown that
\begin{align*}
|\bla T_{\mu, \Phi}(D^1_I f_\omega 1_{I_{\Delta}}& b_1 ,  (E^2_{J_{K^{(1)}, \sigma}}g_\omega) b_2), D^3_K h_\omega  1_{K_{\Delta}}b_3 \bra_\mu | \\
&\lesssim_{\theta} \int | D^1_I f_\omega |1_I M_{\mu, \calD_0(\omega_2)} g_\omega | D^3_Kh_\omega| 1_{K} \, d \mu.
\end{align*}
This in turn shows that
$$
|\eqref{eq:diaggg}| \lesssim_{\theta} \|f\|_{L^p(\mu)} \|g\|_{L^q(\mu)} \|h\|_{L^{r'}(\mu)}  \le 1.
$$

We are now ready with Step 3, since we have proved that
\begin{align*}
|II_3'| &= |E_{\omega_4} E_{\omega} II_3'(\omega)|  \\
& \lesssim C(\theta) + \| T_{\mu, \Phi} \| \E_{\omega_4} \E_{\omega} \Big[ 
\Big \| \Big( \sum_{I \in \calD_0'(\omega_1)} |D^1_I f|^2 1_{I_{\textup{bad}}^{\omega_2, \omega_3, \omega_4}} \Big)^{1/2} \Big\|_{L^p(\mu)} \\
& \quad \quad \quad \quad \quad \quad \quad \quad \quad \quad +
\Big \| \Big( \sum_{K \in \calD_0'(\omega_3)} |D^3_K h|^2 1_{K_{\textup{bad}}^{\omega_1, \omega_2, \omega_4}} \Big)^{1/2} \Big\|_{L^{r'}(\mu)} \Big].
\end{align*}
As previously, this leads to the bound
\begin{equation*}
|II_3'| \lesssim C(\theta) + c(\theta) \| T_{\mu, \Phi} \|
\end{equation*}
using Lemma \ref{lem:badSF}. Therefore, we have shown that
\begin{equation}\label{eq:resStep3}
|II_3| \lesssim C(\theta) + c(\theta) \| T_{\mu, \Phi} \|.
\end{equation}

\subsubsection*{Synthesis of the diagonal}
Taking into account \eqref{eq:resStep1}, \eqref{eq:resStep2} and \eqref{eq:resStep3} we have shown that
$$
|II| \lesssim C(\theta) + c(\theta) \| T_{\mu, \Phi}\|,
$$
which implies that
\begin{equation}\label{eq:resDIAG}
|II| \le C + \| T_{\mu, \Phi}\|/100
\end{equation}
fixing $\theta$ small enough.

Taking into account all the above parts we have finally proved \eqref{eq:reducedtogood}. This ends the proof.
\end{proof}

\chapter{End point estimates}\label{sec:endpoint}

Before proving the weak type bound we need to recall the non-homogeneous Calder\'on--Zygmund decomposition of measures, see \cite{ToBook}.
\begin{lemma}\label{lem:decomp}
Let $\mu$ be a Radon measure in $\R^n$. For every $\nu \in M(\R^n)$ with compact support and
every $\lambda > 2^{n+1}\|\nu\|/\|\mu\|$, we have:
\begin{enumerate}
\item There exists a family of cubes $(Q_i)_i$ so that $\sum_i 1_{Q_i} \le C_n$ and
a function $f \in L^1(\mu)$ such that
\begin{align}
\label{cd1} |\nu|(Q_i) &> \frac{\lambda}{2^{n+1}} \mu(2Q_i), \\
\label{cd2}  |\nu|(\eta Q_i) \le \frac{\lambda}{2^{n+1}}& \mu(2\eta Q_i) \textup{ for } \eta > 2, \\
\label{cd3} \nu = f \,d\mu \textup{ in } \R^n \setminus \bigcup_i &Q_i, \textup{ with } |f| \le \lambda \, \mu\textup{-a.e.}
\end{align}
\item Suppose that for each $i$ we are given a $(6, \beta_0)$-$\mu$-doubling cube $R_i$ such that it is
concentric with $Q_i$ and $Q_i \subset R_i$. For each $i$ set
$$
w_i = \frac{1_{Q_i}}{\sum_k 1_{Q_k}}.
$$
Then there exists a family of functions $(\varphi_i)_i$ (of the form $\varphi_i = \alpha_i h_i$ for some constant $\alpha_i \in \C$ and non-negative function $h_i \ge 0$)
such that
\begin{align}
\label{cd4} \textup{spt}\,\varphi_i &\subset R_i, \\
\label{cd5} \int \varphi_i \,d\mu &= \int w_i \,d\nu, \\
\label{cd6} \sum_i |\varphi_i| \le B\lambda\,\,\,& (B \textup{ depends only on } \beta_0, n), \\
\label{cd7} \|\varphi_i\|_{L^{\infty}(\mu)}& \mu(R_i) \le 2|\nu|(Q_i).
\end{align}
\end{enumerate}
\end{lemma}
We are ready to prove the weak type bound. However, the next proposition is not enough, since we need the
weak type bound for $T_{\sharp}$ (the good lambda method requires this).
To get this we shall need to combine the next proposition with a certain formulation of Cotlar's inequality and some additional arguments.
We give the full details, since we are not aware of a reference covering this type of generality. 

\begin{proposition}\label{prop:weaktypehalfforT}
Let $\mu$ be a measure of order $m$ on $\R^n$ and $T$ be a bilinear $m$-dimensional SIO. Let
$1 < r, p,q < \infty$ be so that $1/p + 1/q = 1/r$, and suppose we have uniformly on $\varepsilon > 0$ that
$$
\|T_{\mu, \varepsilon}\|_{L^p(\mu) \times L^q(\mu) \to L^r(\mu)} \lesssim 1.
$$
Then we have uniformly on $\varepsilon > 0$ that
$$
\|T_{\varepsilon}\|_{M(\R^n) \times M(\R^n) \to L^{1/2,\infty}(\mu)} \lesssim 1.
$$
\end{proposition}
\begin{proof}
We are given $\varepsilon >0$, $\nu, \eta \in M(\R^n)$ and $\lambda > 0$, and want to prove that
$$
\mu(\{x \in \R^n \colon\, |T_{\varepsilon}(\nu, \eta)(x)| > \lambda\}) \lesssim \bigg(\frac{\|\nu\|\|\eta\|}{\lambda}\bigg)^{1/2}.
$$
Without loss of generality we can assume that $\|\nu\| = \|\eta\| = 1$. We can then also assume that $\lambda^{1/2} > 2^{n+1}/\|\mu\|$, since
otherwise the claim is trivial.

Let us first assume that $\nu$ and $\eta$ have compact support. Then we can perform the following decompositions using Lemma \ref{lem:decomp}.
Applying the lemma to the measure $\nu$ on the level $\lambda^{1/2}$ we get cubes
$(Q_{1,i})_i$ and a function $f_1 \in L^1(\mu)$ like in (1) of Lemma \ref{lem:decomp}. For each $i$ let $R_{1,i}$ be the smallest $(6, 6^{m+1})$-$\mu$-doubling
cube of the form $6^k Q_{1,i}$, $k \ge 0$. Such a cube exists by standard arguments, see \cite{ToBook}. Then let $w_{1,i}$ and $\varphi_{1,i}$ be like in (2) of Lemma \ref{lem:decomp}.
We write
$$
\nu = f_1\,d\mu + \sum_i \varphi_{1,i} \,d\mu + \sum_i (w_{1,i}\, d\nu- \varphi_{1,i} \,d\mu) =: g_1\,d\mu + \sum_i \beta_{1,i},
$$
where the function $g_1$ is defined by
$$
g_1 = f_1 + \sum_i \varphi_{1,i} 
$$
and the complex measure $\beta_{1,i}$ is defined by
$$
\beta_{1,i} = w_{1,i}\, d\nu- \varphi_{1,i} \,d\mu.
$$
We also set
$$
\beta_1 =  \sum_i \beta_{1,i}.
$$

It is easy to see the properties
\begin{equation}\label{eq:mes}
\mu\Big(\bigcup_i 2Q_{1,i}\Big) \lesssim \frac{1}{\lambda^{1/2}},
\end{equation}
$\|g_1\|_{L^{\infty}(\mu)} \lesssim \lambda^{1/2}$, $\|g_1\|_{L^1(\mu)} \lesssim 1$ and
\begin{equation}\label{eq:form1}
\|g_1\|_{L^u(\mu)}^u \le \|g\|_{L^{\infty}(\mu)}^{u-1} \|g\|_{L^1(\mu)} \lesssim (\lambda^{1/2})^{u-1}, \qquad 1 < u < \infty.
\end{equation}
Regarding the complex measures $\beta_{1,i}$ we have the following:
\begin{enumerate}[(a)]
\item spt$\,\beta_{1,i} \subset R_{1,i}$;
\item $\beta_{1,i}(R_{1,i}) = 0$;
\item $
\|\beta_{1,i}\| \le 2|\nu|(Q_{1,i}).
$
\end{enumerate}
Finally, the fact that $R_{1,i}$ is the \emph{smallest} $(6, 6^{m+1})$-$\mu$-doubling cube of the form $6^k Q_{1,i}$, $k \ge 0$, is utilised via the standard fact that it implies the estimate
\begin{equation}\label{eq:form2}
\int_{R_{1,i} \setminus Q_{1,i}} \frac{d\mu(x)}{|x-c_{Q_{1,i}}|^m} \lesssim 1.
\end{equation}
We then perform the exact same decomposition using the measure $\eta$. The notation for this decomposition is $Q_{2,j}$, $R_{2,j}$, $f_2$, $g_2$, $\varphi_{2,j}$ etc.

We will separately estimate the naturally appearing good--good, good--bad, bad--good, and bad--bad parts, denoted by $I_{gg}$, $I_{gb}$, $I_{bg}$ and $I_{bb}$ respectively.
For example, we have
$$
I_{gg} := \mu(\{x \in \R^n \colon\, T_{\mu, \varepsilon}(g_1, g_2)(x) > \lambda/4\}).
$$
In fact, the estimate for $I_{gg}$ is trivial using the assumed boundedness of $T_{\mu, \varepsilon}$ and \eqref{eq:form1}:
$$
I_{gg} \lesssim \lambda^{-r} [\|g_1\|_{L^p(\mu)} \|g_2\|_{L^q(\mu)} ]^r \lesssim \lambda^{-r} (\lambda^{1/2})^{2r-1} = \lambda^{-1/2}.
$$

\subsection*{Estimation of $I_{bg}$}We shall now estimate $I_{bg}$ -- the term $I_{gb}$ is handled symmetrically. Because of \eqref{eq:mes} it is enough to estimate
\begin{align*}
&\mu\Big(\Big\{x \in \R^n \setminus \bigcup_i 2Q_{1,i} \colon\, |T_{\varepsilon}(\beta_1, g_2\,d\mu)(x)| > \lambda/4\Big\}\Big) \\
&\lesssim \lambda^{-1} \sum_i \int_{\R^n \setminus 2Q_{1,i}} |T_{\varepsilon}(\beta_{1,i}, g_2\,d\mu)(x)|\,d\mu(x).
\end{align*}
To control this it is enough to fix $i$ and prove
\begin{equation}\label{eq:fixi}
\int_{\R^n \setminus 2Q_{1,i}} |T_{\varepsilon}(\beta_{1,i}, g_2\,d\mu)(x)|\,d\mu(x) \lesssim \lambda^{1/2} |\nu|(Q_{1,i}).
\end{equation}
We write $$\int_{\R^n \setminus 2Q_{1,i}} = \int_{\R^n \setminus 2R_{1,i}} + \int_{ 2R_{1,i} \setminus  2Q_{1,i}},$$ and estimate these separately.

We estimate the integral over $\R^n \setminus 2R_{1,i}$ first. Let $x \in \R^n \setminus 2R_{1,i}$.
Notice that if $d(x,R_{1,i}) > \varepsilon$, then
$$
\operatorname{spt}\, \beta_{1,i} \subset R_{1,i} \subset B(x,\varepsilon)^c.
$$
So if $d(x,R_{1,i}) > \varepsilon$, the fact that
$\beta_{1,i}(R_{1,i}) = 0$ together with the $y$-continuity of $K$ gives that
$$
|T_{\varepsilon}(\beta_{1,i}, g_2\,d\mu)(x)| \lesssim \lambda^{1/2} \frac{\ell(R_{1,i})^{\alpha}}{|x-c_{R_{1,i}}|^{m+\alpha}}  |\nu|(Q_{1,i}),
$$
where we also used that $\|g_2\|_{L^{\infty}(\mu)} \lesssim \lambda^{1/2}$ and $\|\beta_{1,i}\| \lesssim  |\nu|(Q_{1,i})$. This shows that
$$
\mathop{\int_{\R^n \setminus 2R_{1,i}}}_{d(x,R_{1,i}) > \varepsilon} |T_{\varepsilon}(\beta_{1,i}, g_2\,d\mu)(x)|\,d\mu(x) \lesssim \lambda^{1/2} |\nu|(Q_{1,i}).
$$
The size estimate of $K$ gives that
$$
|T_{\varepsilon}(\beta_{1,i}, g_2\,d\mu)(x)| \lesssim \lambda^{1/2}  |\nu|(Q_{1,i}) \frac{1}{\varepsilon^m}.
$$
Notice that if $x \in \R^n \setminus 2R_{1,i}$ satisfies $d(x, R_{1,i}) \le \varepsilon$, then $$|x-c_{R_{1,i}}| \lesssim d(x, R_{1,i}) \le \varepsilon.$$
This gives that
$$
\mathop{\int_{\R^n \setminus 2R_{1,i}}}_{d(x,R_{1,i}) \le \varepsilon} |T_{\varepsilon}(\beta_{1,i}, g_2\,d\mu)(x)|\,d\mu(x) \lesssim
\lambda^{1/2}  |\nu|(Q_{1,i}) \frac{\mu(B(c_{R_{1,i}}, C\varepsilon))}{\varepsilon^m} \lesssim \lambda^{1/2}  |\nu|(Q_{1,i}).
$$

We have shown that
$$
\int_{\R^n \setminus 2R_{1,i}} |T_{\varepsilon}(\beta_{1,i}, g_2\,d\mu)(x)|\,d\mu(x) \lesssim \lambda^{1/2} |\nu|(Q_{1,i}),
$$
and will next show that the same bound holds for
$$
\int_{ 2R_{1,i} \setminus  2Q_{1,i}} |T_{\varepsilon}(\beta_{1,i}, g_2\,d\mu)(x)|\,d\mu(x).
$$
For this it is enough to separately bound
$$
\int_{ 2R_{1,i} \setminus  2Q_{1,i}} |T_{\varepsilon}(w_{1,i}\,d\nu, g_2\,d\mu)(x)|\,d\mu(x) \,\textup{ and }\, \int_{ 2R_{1,i} } |T_{\mu, \varepsilon}(\varphi_{1,i}, g_2)(x)|\,d\mu(x).
$$
We begin with the first integral. The size estimate gives
$$
|T_{\varepsilon}(w_{1,i}\,d\nu, g_2\,d\mu)(x)| \lesssim \lambda^{1/2} |\nu|(Q_{1,i}) \frac{1}{|x-c_{Q_{1,i}}|^m}, \qquad x \not \in 2Q_{1,i},
$$
so we have by \eqref{eq:form2} that
$$
\int_{ 2R_{1,i} \setminus  2Q_{1,i}} |T_{\varepsilon}(w_{1,i}\,d\nu, g_2\,d\mu)(x)|\,d\mu(x) \lesssim  \lambda^{1/2} |\nu|(Q_{1,i}).
$$

Next, we bound
\begin{align*}
 \int_{ 2R_{1,i} } |T_{\mu, \varepsilon}(\varphi_{1,i}, g_2)(x)|\,d\mu(x) \le  \int_{ 2R_{1,i} }& |T_{\mu, \varepsilon}(\varphi_{1,i}, 1_{4R_{1,i}} g_2)(x)|\,d\mu(x) \\
&  + \int_{ 2R_{1,i} } |T_{\mu, \varepsilon}(\varphi_{1,i}, 1_{(4R_{1,i})^c} g_2)(x)|\,d\mu(x).
\end{align*}
Using the assumed boundedness of $T_{\mu, \varepsilon}$ and recalling that $R_{1,i}$ is doubling we see that
\begin{align*}
\int_{ 2R_{1,i} } |T_{\mu, \varepsilon}(\varphi_{1,i}, 1_{4R_{1,i}}g_2)|\,d\mu &\lesssim \mu(R_{1,i})^{1-1/r} \cdot \|\varphi_{1,i}\|_{L^{\infty}(\mu)} \mu(R_{1,i})^{1/p} \cdot \lambda^{1/2} \mu(R_{1,i})^{1/q} \\
&= \lambda^{1/2}\|\varphi_{1,i}\|_{L^{\infty}(\mu)} \mu(R_{1,i}) \lesssim \lambda^{1/2}|\nu|(Q_{1,i}).
\end{align*}
We move on to bounding
$$
\int_{ 2R_{1,i} } |T_{\mu, \varepsilon}(\varphi_{1,i}, 1_{(4R_{1,i})^c} g_2)(x)|\,d\mu(x).
$$
Notice that for $x \in 2R_{1,i}$ we have
$$
\int_{(4R_{1,i})^c} \int_{R_{1,i}} \frac{d\mu(y) d\mu(z)}{(|x-y|+|x-z|)^{2m}} \lesssim \mu(R_{1,i}) \int_{R_{1,i}^c} \frac{d\mu(z)}{|z-c_{R_{1,i}}|^{2m}} \lesssim \frac{\mu(R_{1,i})}{\ell(R_{1,i})^m} \lesssim 1,
$$
so that
$$
\int_{ 2R_{1,i} } |T_{\mu, \varepsilon}(\varphi_{1,i}, 1_{(4R_{1,i})^c} g_2)(x)|\,d\mu(x) \lesssim \lambda^{1/2} \|\varphi_{1,i}\|_{L^{\infty}(\mu)} \mu(R_{1,i}) \lesssim \lambda^{1/2}|\nu|(Q_{1,i}).
$$
Putting everything together we have shown \eqref{eq:fixi}, which ends our treatment of the term $I_{bg}$.

\subsection*{Estimation of $I_{bb}$}
Now we turn to estimate the final part
$$
I_{bb} = \mu\big(\big\{x \in \R^n \colon |T_{\varepsilon}(\beta_1,\beta_2)(x)| >\lambda/4\big\}\big).
$$
Let $\scrA:= \bigcup_i 2\Qi \cup \bigcup_j 2\Qj$. Since $\mu(\scrA) \lesssim  \lambda^{-1/2}$, it is enough to consider
$$
\mu\big(\big\{x \in \R^n\setminus \scrA \colon |T_{\varepsilon}(\beta_1,\beta_2)(x)| >\lambda/4\big\}\big).
$$

First, we divide $T_{\varepsilon}(\beta_1,\beta_2)$ into two symmetric parts according to the relative side lengths of the cubes $\Ri$ and $\Rj$. Namely, we have
$$
T_{\varepsilon}(\beta_1,\beta_2)
= \sum_i T_{\varepsilon}\bigg(\bei, \sum_{\begin{substack}{j\colon \\ \ell(\Ri) \leq \ell(\Rj)}\end{substack}}\bej\bigg)
+\sum_j T_{\varepsilon}\bigg(\sum_{\begin{substack}{i \colon \\ \ell(\Ri) > \ell(\Rj)}\end{substack}}\bei, \bej\bigg).
$$
These two terms are handled symmetrically, so we focus on the first. 
Define the sets of indices $$\calJ_i:= \big\{j \colon \ell(\Ri) \leq \ell(\Rj)\big\}.$$ 
We have
\begin{equation}\label{eq:main cases}
\begin{split}
\mu\Big(\Big\{x &\in \R^n \setminus \scrA \colon \sum_i \Big| T_{\varepsilon}\Big(\bei, \sum_{j \in \calJ_i}\bej\Big)(x) \Big| >\lambda/8 \Big\}\Big) \\
&\leq \mu\Big(\Big\{x \in \R^n \setminus \scrA \colon \sum_i 1_{(2\Ri)^c}(x) \Big|T_{\varepsilon}\Big(\bei, \sum_{j \in \calJ_i}\bej\Big)(x)\Big| >\lambda/16\Big\}\Big) \\
& + \mu\Big(\Big\{x \in \R^n \setminus \scrA \colon \sum_i 1_{2\Ri}(x) \Big|T_{\varepsilon}\Big(\bei, \sum_{j \in \calJ_i}\bej\Big)(x)\Big| >\lambda/16\Big\}\Big)\\
&=: I + II.
\end{split}
\end{equation}
These two cases will be handled separately.

We begin the estimation of $I$. We have
\begin{equation*}
\begin{split}
\lambda^{1/2}I 
&\lesssim  \int_{\R^n} \bigg( \sum_{i,j}
1_{(2\Ri)^c}1_{(2\Rj)^c}|T_{\varepsilon}(\bei,\bej)| \bigg)^{1/2} d\mu \\
&+  \int_{\R^n \setminus \scrA} \bigg( \sum_{i,j}
1_{(2\Ri)^c}1_{2\Rj}|T_{\varepsilon}(\bei,\bej)| \bigg)^{1/2} d \mu  \\
& =:I_a+I_b.
\end{split}
\end{equation*}
Notice that we dropped the restriction $j \in \calJ_i$.

To control $I_a$ consider some $x \in (2\Ri)^c \cap (2\Rj)^c$. Suppose first that $d(x, R_{1,i}) > \varepsilon$. Then, since $\Ri \subset \bar B(x,\varepsilon)^c$, we may estimate
\begin{align*}
|T_{\varepsilon}(\bei,\bej)(x)| &= \bigg| \int_{\Rj}\int_{\Ri} [K(x,y,z)-K(x,c_{\Ri}, z)]\,d\bei(y)\,d\bej(z)\bigg| \\
&\lesssim \ell(R_{1,i})^{\alpha} \frac{|\nu|(\Qi) |\eta|(\Qj)}{(|x-c_{\Ri}| + |x-c_{\Rj}|)^{2m+\alpha}} \\
&\lesssim \ell(R_{1,i})^{\alpha} \frac{|\nu|(\Qi)}{|x-c_{\Ri}|^{m+\alpha/2}}  \frac{|\eta|(\Qj)}{|x-c_{\Rj}|^{m+\alpha/2}}.
\end{align*}
We could also use the cancellation in $\beta_{2,j}$, and so it actually holds that
$$
|T_{\varepsilon}(\bei,\bej)(x)| \lesssim  \ell(\Ri)^{\alpha/2} \frac{|\nu|(\Qi)}{|x-c_{\Ri}|^{m+\alpha/2}} \cdot  \ell(\Rj)^{\alpha/2} \frac{|\eta|(\Qj)}{|x-c_{\Rj}|^{m+\alpha/2}}.
$$
This gives that
\begin{align*}
&\int_{\R^n} \bigg(\sum_{i,j}1_{(2\Ri)^c} 1_{\{d(\cdot,\Ri) > \varepsilon\}}  1_{(2\Rj)^c}  |T_{\varepsilon}(\bei,\bej)|\bigg)^{1/2}d \mu \\
\lesssim &\bigg(\int_{\R^n} \sum_i 1_{\Ri^c}(x)
\frac{\ell(\Ri)^{\alpha/2}|\nu|(\Qi)}{|x-c_{\Ri}|^{m+\alpha/2}}d \mu(x) \bigg)^{1/2} \\
\times &\bigg(\int_{\R^n} \sum_j 1_{\Rj^c}(x)
\frac{\ell(\Rj)^{\alpha/2}|\eta|(\Qj)}{|x-c_{\Rj}|^{m+\alpha/2}}d \mu(x) \bigg)^{1/2} \\
\lesssim  &\Big( \sum_i |\nu|(\Qi) \Big)^{1/2} 
\Big(\sum_j |\eta|(\Qj)\Big)^{1/2} \lesssim 1.
\end{align*}
The size estimate gives for $x \in (2\Rj)^c$ that
$$
|T_{\varepsilon}(\bei,\bej)(x)| \lesssim |\nu|(\Qi) \cdot \frac{|\eta|(\Qj)}{(\varepsilon +  |x-c_{\Rj}|)^{2m}}.
$$
Notice that if $x \in \R^n \setminus 2R_{1,i}$ satisfies $d(x, R_{1,i}) \le \varepsilon$, then $|x-c_{R_{1,i}}| \lesssim d(x, R_{1,i}) \le \varepsilon$.
This gives that
\begin{align*}
&\int_{\R^n} \bigg(\sum_{i,j}1_{(2\Ri)^c} 1_{\{d(\cdot,\Ri) \le \varepsilon\}}  1_{(2\Rj)^c}  |T_{\varepsilon}(\bei,\bej)|\bigg)^{1/2}d \mu \\
&\lesssim \int_{\R^n} \bigg(\sum_{i} 1_{B(c_{\Ri}, C\varepsilon)}(x) |\nu|(\Qi) \bigg)^{1/2} \bigg(\sum_{j} \frac{|\eta|(\Qj)}{(\varepsilon +  |x-c_{\Rj}|)^{2m}}  \bigg)^{1/2}d\mu(x) \\
&\lesssim \bigg(\sum_{i} \varepsilon^m  |\nu|(\Qi) \bigg)^{1/2} \bigg(\sum_{j} |\eta|(\Qj) \varepsilon^{-m}  \bigg)^{1/2} \lesssim 1.
\end{align*}
Combining everything, we have shown that $I_a \lesssim 1$.

We turn our attention to the term $I_b$. We estimate
\begin{align*}
I_b \le  \int_{\R^n}& \bigg( \sum_{i,j} 1_{(2\Ri)^c}1_{2\Rj}|T_{\varepsilon}(\bei,\phij\,d\mu)| \bigg)^{1/2} d \mu \\
&+  \int_{\R^n} \bigg( \sum_{i,j} 1_{(2\Ri)^c}1_{2\Rj \setminus 2\Qj}|T_{\varepsilon}(\bei,\wj\,d\eta)| \bigg)^{1/2} d \mu = I_b' + I_b''.
\end{align*}
In the second term we were able to change $1_{2\Ri}$ into $1_{2\Ri \setminus 2\Qi}$, since  the integral in $I_b$ is over $\R^n \setminus \scrA$. After this change we omitted the restriction of the integral to the complement of $\scrA$.

To deal with $I_b'$ we have to once again perform the usual trickery involving the truncation parameter $\varepsilon$.
If $x \in (2\Ri)^c$ satisfies $d(x, \Ri) > \varepsilon$, then the H\"older estimate in the $y$-variable yields
$$
|T_{\varepsilon}(\bei,\phij\,d\mu)(x)| \lesssim  \frac{\ell(\Ri)^{\alpha} |\nu| (\Qi)}{ |x-c_{\Ri}|^{m+\alpha}} \| \phij \|_{L^\infty(\mu)}.
$$
This leads to the bound
\begin{align*}
\int_{\R^n} &\bigg( \sum_{i,j} 1_{(2\Ri)^c} 1_{\{d(\cdot, R_{1,i}) > \varepsilon\}} 1_{2\Rj}|T_{\varepsilon}(\bei,\phij\,d\mu)| \bigg)^{1/2} d \mu \\
&\lesssim \bigg(\int_{\R^n} \sum_i 1_{\Ri^c}(x)
\frac{\ell(\Ri)^{\alpha}|\nu|(\Qi)}{|x-c_{\Ri}|^{m+\alpha}}\,d \mu(x) \bigg)^{1/2}  \\
&\times  \bigg(\int_{\R^n} \sum_j \|\phij \|_{L^\infty(\mu)} 1_{2\Rj}(x)\, d\mu(x) \bigg)^{1/2} \lesssim 1.
\end{align*}
On the other hand, the size estimate gives
$$
|T_{\varepsilon}(\bei,\phij\,d\mu)(x)| \lesssim \frac{|\nu|(\Qi)}{\varepsilon^m} \|\phij \|_{L^\infty(\mu)},
$$
which leads to the bound
\begin{align*}
\int_{\R^n} &\bigg( \sum_{i,j} 1_{(2\Ri)^c} 1_{\{d(\cdot, R_{1,i}) \le \varepsilon\}} 1_{2\Rj}|T_{\varepsilon}(\bei,\phij\,d\mu)| \bigg)^{1/2} d \mu \\
&\lesssim  \bigg(\int_{\R^n} \sum_i 1_{B(c_{\Ri}, C\varepsilon)}(x)\frac{|\nu|(\Qi)}{\varepsilon^m} \,d \mu(x) \bigg)^{1/2} \\
&\times \bigg(\int_{\R^n} \sum_j \|\phij \|_{L^\infty(\mu)} 1_{2\Rj}(x)\, d\mu(x) \bigg)^{1/2} \lesssim 1.
\end{align*}
This shows that $I_b' \lesssim 1$.

Let us now control $I_b''$. For $x \in (2\Ri)^c \cap (2\Qj)^c$ satisfying $d(x, \Ri) > \varepsilon$, we have using the H\"older estimate in the $y$-variable that
$$
|T_{\varepsilon}(\bei,\wj\,d\eta)(x)| \lesssim \frac{\ell(\Ri)^{\alpha} |\nu| (\Qi)}{ |x-c_{\Ri}|^{m+\alpha}} \frac{ |\eta| (\Qj)}{ |x-c_{\Qj}|^{m}},
$$
which gives the bound
\begin{align*}
 \int_{\R^n} &\bigg( \sum_{i,j} 1_{(2\Ri)^c}1_{\{d(\cdot, R_{1,i}) > \varepsilon\}}1_{2\Rj \setminus 2\Qj}|T_{\varepsilon}(\bei,\wj\,d\eta)| \bigg)^{1/2} d \mu \\
 &\lesssim \bigg(\int_{\R^n} \sum_i 1_{\Ri^c}(x) \frac{\ell(\Ri)^{\alpha} |\nu| (\Qi)}{ |x-c_{\Ri}|^{m+\alpha}} \,d \mu(x) \bigg)^{1/2} \\
 &\times \bigg(\int_{\R^n} \sum_j  1_{2\Rj \setminus \Qj} \frac{ |\eta| (\Qj)}{ |x-c_{\Qj}|^{m}} \,d \mu(x) \bigg)^{1/2} \lesssim 1.
\end{align*}
Here we used \eqref{eq:form2} to estimate the integrals over $2\Rj \setminus \Qj$. For $x \in (2\Qj)^c$ we have using the size estimate that
$$
|T_{\varepsilon}(\bei,\wj\,d\eta)(x)| \lesssim \frac{|\nu|(\Qi)}{\varepsilon^m} \frac{ |\eta| (\Qj)}{ |x-c_{\Qj}|^{m}},
$$
which leads to the bound
\begin{align*}
 \int_{\R^n} &\bigg( \sum_{i,j} 1_{(2\Ri)^c}1_{\{d(\cdot, R_{1,i}) \le \varepsilon\}}1_{2\Rj \setminus 2\Qj}|T_{\varepsilon}(\bei,\wj\,d\eta)| \bigg)^{1/2} d \mu \\
 &\lesssim \bigg(\int_{\R^n} \sum_i 1_{B(c_{\Ri}, C\varepsilon)}(x)\frac{|\nu|(\Qi)}{\varepsilon^m} \,d \mu(x) \bigg)^{1/2} \\
&\times \bigg(\int_{\R^n} \sum_j  1_{2\Rj \setminus \Qj} \frac{ |\eta| (\Qj)}{ |x-c_{\Qj}|^{m}} \,d \mu(x) \bigg)^{1/2} \lesssim 1.
 \end{align*}
Therefore, $I_b'' \lesssim 1$, and so $I_b \lesssim 1$. We have shown that
 $$
 I \lesssim \lambda^{-1/2}(I_a + I_b) \lesssim  \lambda^{-1/2}.
 $$
 
It remains to show that
$$
II =  \mu\Big(\Big\{x \in \R^n \setminus \scrA \colon \sum_i 1_{2\Ri}(x) \Big|T_{\varepsilon}\Big(\bei, \sum_{j \in \calJ_i}\bej\Big)(x)\Big| >\lambda/16\Big\}\Big) \lesssim \lambda^{-1/2}.
$$
To reduce things somewhat, notice that
\begin{equation*}
\begin{split}
 \lambda^{1/2}\mu\Big(\Big\{x \in & \R^n \setminus \scrA \colon \sum_i 1_{2\Ri}(x)
 \Big|T_{\varepsilon}\Big(\bei, \sum_{\substack{j \in \calJ_i \\ 2\Ri \cap 2\Rj = \emptyset}}\bej\Big)(x)\Big| >\lambda/32\Big\}\Big) \\
& \lesssim \int_{\R^n \setminus \scrA} \bigg( \sum_{i,j} 1_{2\Ri}1_{(2\Rj)^c}|T_{\varepsilon}(\bei,\bej)| \bigg)^{1/2} d \mu \lesssim 1,
\end{split}
\end{equation*}
where we used that the appearing term is similar with $I_b$ from above.
Define 
$$
\wt{\calJ_i}:=\big\{ j \in \calJ_i \colon2\Ri \cap 2\Rj \not= \emptyset\big\}.
$$  
After splitting $$\bei= \phii d \mu+ \wi d \nu$$ and $$\bej= \phij d \mu+ \wj d \eta,$$ what remains to be done  is to estimate the following four terms:
\begin{align*}
&II_a:=\mu\Big(\Big\{x \in  \R^n \setminus \scrA \colon \sum_i 1_{2\Ri}(x) \Big|T_{\mu, \varepsilon}\Big(\phii , \sum_{j \in \wt{\calJ_i}}\phij\Big)(x)\Big| >\lambda/128\Big\}\Big), \\
&II_b:=\mu\Big(\Big\{x \in  \R^n \setminus \scrA \colon \sum_i 1_{2\Ri}(x) \Big|T_{\varepsilon}\Big(\phii d \mu, \sum_{j \in \wt{\calJ_i}}\wj d \eta\Big)(x)\Big| >\lambda/128\Big\}\Big), \\
&II_c:=\mu\Big(\Big\{x \in  \R^n \setminus \scrA \colon \sum_i 1_{2\Ri}(x) \Big|T_{\varepsilon}\Big(\wi d \nu, \sum_{j \in \wt{\calJ_i}}\phij d \mu\Big)(x)\Big| >\lambda/128\Big\}\Big), \\
& II_d:=\mu\Big(\Big\{x \in  \R^n \setminus \scrA \colon \sum_i 1_{2\Ri}(x) \Big|T_{\varepsilon}\Big(\wi d \nu, \sum_{j \in \wt{\calJ_i}}\wj d \eta\Big)(x)\Big| >\lambda/128\Big\}\Big). \\
\end{align*}

Estimate $$II_a \le II_a' + II_a'',$$ where
$$
II_a' = \mu\Big(\Big\{x \in  \R^n \setminus \scrA \colon \sum_i 1_{2\Ri}(x) \Big|T_{\mu, \varepsilon}\Big(\phii , 1_{4\Ri}\sum_{j \in \wt{\calJ_i}}\phij\Big)(x)\Big| >\lambda/256\Big\}\Big),
$$
and $II_a''$ is defined in the obvious way with $1_{4\Ri}$ replaced by $1_{(4\Ri)^c}$ inside $T_{\mu, \varepsilon}$. We have using the boundedness of $T_{\mu, \varepsilon}$
that
\begin{align*}
II_a' &\lesssim  \lambda^{-1} \int \sum_i 1_{2\Ri} \Big|T_{\mu, \varepsilon}\Big(\phii , 1_{4\Ri}\sum_{j \in \wt{\calJ_i}}\phij\Big)\Big| \,d \mu \\
&\lesssim  \lambda^{-1} \sum_i  \mu(\Ri)^{1-1/r} \| \phii \|_{L^\infty(\mu)} \mu(\Ri)^{1/p} \Big\| \sum_j |\phij| \Big\|_{L^{\infty}(\mu)} \mu(\Ri)^{1/q} \\
&\lesssim \lambda^{-1/2} \sum_i \mu(\Ri)\| \phii \|_{L^\infty(\mu)} \lesssim \lambda^{-1/2}.
\end{align*}
Notice that we used that $\sum_j |\phij| \lesssim \lambda^{1/2}$. Using the size estimate we see that for $x \in 2\Ri$ it holds that
\begin{align*}
\Big|T_{\mu, \varepsilon}\Big(\phii , 1_{(4\Ri)^c}\sum_{j \in \wt{\calJ_i}}\phij\Big)\Big| &\lesssim \| \phii \|_{L^\infty(\mu)} \Big\| \sum_j |\phij| \Big\|_{L^{\infty}(\mu)}
\int_{\Ri^c} \int_{\Ri} \frac{d\mu(y)\,d\mu(z)}{|z-c_{\Ri}|^{2m}} \\
&\lesssim \lambda^{1/2} \| \phii \|_{L^\infty(\mu)} \frac{\mu(\Ri)}{\ell(\Ri)^m} \lesssim \lambda^{1/2} \| \phii \|_{L^\infty(\mu)}.
\end{align*}
Therefore, we have
$$
II_a'' \lesssim \lambda^{-1} \int \sum_i 1_{2\Ri} \lambda^{1/2} \| \phii \|_{L^\infty(\mu)}\,d\mu \lesssim \lambda^{-1/2},
$$
and this completes the proof of the fact that $II_a \lesssim \lambda^{-1/2}$.

For $II_b$ notice that the size estimate gives for $x \in \R^n \setminus \scrA \subset (2\Qj)^c$ that
$$
|T_{\varepsilon}(\phii\,d\mu, \wj\,d\eta)(x)| \lesssim \frac{ \| \phii \|_{L^\infty(\mu)} | \eta |(\Qj)}{|x-c_{\Qj}|^m}.
$$
Notice also that if $j \in \wt{\calJ_i}$, then $2\Ri \subset 6\Rj$. Therefore, we get
\begin{align*}
\lambda^{1/2}II_b \lesssim  \int_{\R^n} \bigg( \sum_{i,j} 1_{2\Ri}(x)1_{6\Rj \setminus \Qj}(x)\frac{ \| \phii \|_{L^\infty(\mu)} | \eta |(\Qj)}{|x-c_{\Qj}|^m} \bigg)^{1/2} d \mu(x).
\end{align*}
From here the estimate is concluded as before, using first H\"older's inequality in $L^2(\mu)$, and then estimating the resulting two integrals with the help of equations \eqref{cd7} and \eqref{eq:form2}. This shows that $II_b \lesssim \lambda^{-1/2}$. The estimate $II_c \lesssim \lambda^{-1/2}$ is concluded with essentially same arguments.
Regarding the term $II_d$, we have
\begin{equation*}
|T_{\varepsilon}( \wi d \nu, \wj d \eta)(x)|
\lesssim \frac{ |\nu|(\Qi) |\eta|(\Qj)}{|x-c_{\Qi}|^m |x-c_{\Qj}|^m}, \quad x \in (2\Qi)^c \cap (2 \Qj)^c.
\end{equation*}
This allows to estimate $II_d$ with similar steps as we used with $II_b$.

This finally almost concludes the proof. It remains to note that it is straightforward to drop the assumption that $\nu$ and 
$\eta$ have compact support. The argument goes quite similarly as in the linear case (see e.g. \cite{ToBook}).
\end{proof}

Next, we prove a version of Cotlar's inequality.
\begin{proposition}\label{prop:cotbas}
Let $\mu$ be a measure of order $m$ on $\R^n$ and $T$ be a bilinear $m$-dimensional SIO.
Let $\delta > 0$ and suppose that
$$
\|T_{\delta}\|_{M(\R^n) \times M(\R^n) \to L^{1/2,\infty}(\mu)} \lesssim 1.
$$
Then for all $\nu_1, \nu_2 \in M(\R^n)$ and $x \in \R^n$ we have uniformly on $\varepsilon > \delta$ that
$$
|T_{\varepsilon}(\nu_1, \nu_2)(x)| \lesssim N_{\mu, 1/4}(T_{\delta}(\nu_1, \nu_2))(x) + M_{\mu}\nu_1(x)M_{\mu}\nu_2(x).
$$
\end{proposition}
\begin{proof}
Fix $\nu_1, \nu_2 \in M(\R^n)$, $x \in \R^n$ and $\varepsilon_0 > 0$. We will estimate $|T_{\varepsilon_0}(\nu_1, \nu_2)(x)|$.
For convenience, we will throughout the proof denote restrictions of measures, like $\nu \rest A$, with $1_A \nu$.
Choose the smallest $k$ so that $B(x, 5^k\varepsilon_0)$ is $(5, 5^{m+1})$ doubling with respect to $\mu$. Set $\varepsilon = 5^k \varepsilon_0$. 

We will begin by controlling $|T_{\varepsilon_0}(\nu_1, \nu_2)(x) - T_{\varepsilon}(\nu_1, \nu_2)(x)|$ -- which is a standard argument for experts.
Notice that
\begin{align*}
|T_{\varepsilon_0}&(\nu_1, \nu_2)(x) - T_{\varepsilon}(\nu_1, \nu_2)(x)| \\
&\lesssim \iint_{\varepsilon_0 < \max(|x-y|,|x-z|) \le \varepsilon} \frac{d|\nu_1|(y)\,d|\nu_2|(z)}{(|x-y| + |x-z|)^{2m}} \\
&\lesssim M_{\mu} \nu_2(x)  \int_{\varepsilon_0 < |x-y| \le \varepsilon} \frac{d|\nu_1|(y)}{|x-y|^{m}}
+ M_{\mu} \nu_1(x)  \int_{\varepsilon_0 < |x-z| \le \varepsilon} \frac{d|\nu_2|(z)}{|x-z|^{m}}.
\end{align*}
These terms are completely symmetric, so it suffices to deal with the first.
Notice that
$$
\int_{|x-y| = \varepsilon} \frac{d|\nu_1|(y)}{|x-y|^{m}} \le \frac{{|\nu_1|(B(x,2\varepsilon))}}{\varepsilon^m} \lesssim M_{\mu} \nu_1(x).
$$
We bound
\begin{align*}
\int_{\varepsilon_0 < |x-y| < \varepsilon} \frac{d|\nu_1|(y)}{|x-y|^{m}} &\le \sum_{j=0}^{k-1} \int_{5^j\varepsilon_0 \le |x-y| < 5^{j+1}\varepsilon} \frac{d|\nu_1|(y)}{|x-y|^{m}} \\
&\le  \sum_{j=0}^{k-1} (5^j\varepsilon_0)^{-m} |\nu_1|(B(x,5^{j+1}\varepsilon_0)) \\
&\le M_{\mu} \nu_1(x) \sum_{j=0}^{k-1} (5^j\varepsilon_0)^{-m} \mu(B(x,5^{j+1}\varepsilon_0)).
\end{align*}
Since
$$
\mu(B(x,5^{j+1}\varepsilon_0)) \le (5^{-m-1})^{k-j-1}\mu(B(x,\varepsilon)) \lesssim (5^{-m-1})^{k-j}\varepsilon^m = (5^{-m-1})^{k-j}(2^k\varepsilon_0)^m,
$$ 
it follows that
$$
\sum_{j=0}^{k-1} (5^j\varepsilon_0)^{-m} \mu(B(x,5^{j+1}\varepsilon_0)) \lesssim \sum_{j=0}^{\infty} \Big(\frac{1}{5}\Big)^j \lesssim 1.
$$
We have shown that
$$
|T_{\varepsilon_0}(\nu_1, \nu_2)(x)| \lesssim |T_{\varepsilon}(\nu_1, \nu_2)(x)| + M_{\mu} \nu_1(x)M_{\mu} \nu_2(x),
$$
and so are reduced to bounding $|T_{\varepsilon}(\nu_1, \nu_2)(x)|$. 

For a fixed $w \in B(x,\varepsilon)$ write
\begin{align*}
T_{\varepsilon}(\nu_1, \nu_2)(x) =  T_{\varepsilon}&(\nu_1, \nu_2)(x) - T_{\delta}(1_{B(x,2\varepsilon)^c}\nu_1, \nu_2)(w) \\
&+ T_{\delta}(\nu_1, \nu_2)(w) - T_{\delta}(1_{B(x,2\varepsilon)}\nu_1, \nu_2)(w).
\end{align*}
Now for every $w \in B(x,\varepsilon)$ we have $B(w,\delta) \subset B(x,2\varepsilon)$ so that we can dominate
$$
|T_{\varepsilon}(\nu_1, \nu_2)(x) - T_{\delta}(1_{B(x,2\varepsilon)^c}\nu_1, \nu_2)(w)|
$$
with the sum of
\begin{align*}
\int_{\R^n}& \int_{B(x,2\varepsilon)^c} |K(x,y,z)-K(w,y,z)|\,d|\nu_1|(y)| \,d|\nu_2|(z) \\
&\lesssim \int_{B(x,2\varepsilon)^c} \int_{\R^n} \frac{\varepsilon^{\alpha}}{(|x-y| + |x-z|)^{2m+\alpha}} \,d|\nu_2|(z)| \,d|\nu_1|(y) \\
&\lesssim M_{m}\nu_2(x) \cdot \varepsilon^{\alpha} \int_{B(x,\varepsilon)^c} \frac{d|\nu_1|(y)}{|x-y|^{m+\alpha}}
\lesssim  M_{m}\nu_1(x) M_{m}\nu_2(x) \lesssim M_{\mu}\nu_1(x) M_{\mu}\nu_2(x).
\end{align*}
and
\begin{align*}
\mathop{\iint_{\max(|x-y|, |x-z|) > \varepsilon}}_{y \in B(x, 2\varepsilon)}  |K(x,y,z)| \,d|\nu_1|(y) \,d|\nu_2|(z) 
&\lesssim \int_{B(x,2\varepsilon)} \int \frac{d|\nu_2|(z)}{(\varepsilon + |x-z|)^{2m}}\,d|\nu_1|(y) \\
&\lesssim \frac{|\nu_1|(B(x,2\varepsilon))}{\varepsilon^m}M_{\mu} \nu_2(x) \\
& \lesssim M_{\mu}\nu_1(x) M_{\mu}\nu_2(x).
\end{align*}

The above shows that for all $w \in B(x,\varepsilon)$ we have
$$
|T_{\varepsilon}(\nu_1, \nu_2)(x)| \lesssim M_{\mu}\nu_1(x) M_{\mu}\nu_2(x) + |T_{\delta}(\nu_1, \nu_2)(w)| + |T_{\delta}(1_{B(x,2\varepsilon)}\nu_1, \nu_2)(w)|.
$$
It follows from this by raising to the power $1/4$, averaging over $w \in B(x,\varepsilon)$ and raising to power $4$ that
$$
|T_{\varepsilon}(\nu_1, \nu_2)(x)| \lesssim M_{\mu}\nu_1(x) M_{\mu}\nu_2(x) + I + II,
$$
where
$$
I := \bigg( \frac{1}{\mu(B(x,\varepsilon))} \int_{B(x,\varepsilon)} |T_{\delta}(\nu_1, \nu_2)(w)|^{1/4}\,d\mu(w) \bigg)^4
$$
and
$$
II :=  \bigg( \frac{1}{\mu(B(x,\varepsilon))} \int_{B(x,\varepsilon)} |T_{\delta}(1_{B(x,2\varepsilon)} \nu_1, \nu_2)(w)|^{1/4}\,d\mu(w) \bigg)^4.
$$
Since $\mu(B(x,5\varepsilon)) \le 5^{m+1}\mu(B(x,\varepsilon))$ we have
$$
I \lesssim N_{\mu, 1/4}(T_{\delta}(\nu_1, \nu_2))(x).
$$

It only remains to estimate the term $II$, which we begin by estimating
\begin{align*}
II &\lesssim  \bigg( \frac{1}{\mu(B(x,\varepsilon))} \int_{B(x,\varepsilon)} |T_{\delta}(1_{B(x,2\varepsilon)} \nu_1, 1_{B(x,2\varepsilon)^c}\nu_2)(w)|^{1/4}\,d\mu(w) \bigg)^4 \\
&+ \bigg( \frac{1}{\mu(B(x,\varepsilon))} \int_{B(x,\varepsilon)} |T_{\delta}(1_{B(x,2\varepsilon)} \nu_1, 1_{B(x,2\varepsilon)}\nu_2)(w)|^{1/4}\,d\mu(w) \bigg)^4 = II' + II''.
\end{align*}
Notice that for all $w \in B(x,\varepsilon)$ we have
\begin{align*}
|T_{\delta}(1_{B(x,2\varepsilon)} \nu_1, 1_{B(x,2\varepsilon)^c}\nu_2)(w)| &\lesssim \int_{B(x,2\varepsilon)} \int_{B(x,2\varepsilon)^c} \frac{d|\nu_2|(z)}{|w-z|^{2m}}\,d|\nu_1|(y) \\
&\lesssim |\nu_1|(B(x,2\varepsilon)) \int_{B(x,\varepsilon)^c} \frac{d|\nu_2|(z)}{|z-x|^{2m}} \\
&\lesssim M_m\nu_2(x) \cdot \frac{|\nu_1|(B(x,2\varepsilon)}{\varepsilon^m} \lesssim M_{\mu}\nu_1(x) M_{\mu}\nu_2(x)
\end{align*}
so that
$$
II' \lesssim M_{\mu}\nu_1(x) M_{\mu}\nu_2(x).
$$

We are left with $II''$, which we will handle using the assumption
$$\|T_{\delta}\|_{M(\R^n) \times M(\R^n) \to L^{1/2,\infty}(\mu)} \lesssim 1,$$
the doubling property of the ball $B(x,\varepsilon)$ and some Kolmogorov type arguments.
We have
\begin{align*}
&\int_{B(x,\varepsilon)} |T_{\delta}(1_{B(x,2\varepsilon)} \nu_1, 1_{B(x,2\varepsilon)}\nu_2)(w)|^{1/4}\,d\mu(w)  \\
&= \frac{1}{4} \int_0^{\infty} \lambda^{-3/4} \mu(\{w \in B(x,\varepsilon)\colon\,  |T_{\delta}(1_{B(x,2\varepsilon)} \nu_1, 1_{B(x,2\varepsilon)}\nu_2)(w)| > \lambda\})\,d\lambda \\
&\lesssim \int_0^{A} \lambda^{-3/4} \,d\lambda \cdot \mu(B(x,\varepsilon)) + \int_A^{\infty} \lambda^{-5/4} \,d\lambda \cdot |\nu_1|(B(x,2\varepsilon))^{1/2} |\nu_2|(B(x,2\varepsilon))^{1/2} \\
&\lesssim A^{1/4} \mu(B(x,\varepsilon))  + A^{-1/4}|\nu_1|(B(x,2\varepsilon))^{1/2} |\nu_2|(B(x,2\varepsilon))^{1/2} \\
&\lesssim |\nu_1|(B(x,2\varepsilon))^{1/4} |\nu_2|(B(x,2\varepsilon))^{1/4} \mu(B(x,\varepsilon))^{1/2},
\end{align*}
where we used the choice
$$
A = \frac{ |\nu_1|(B(x,2\varepsilon)) |\nu_2|(B(x,2\varepsilon))}{\mu(B(x,\varepsilon))^2}.
$$
This gives
\begin{align*}
II'' &\lesssim \frac{|\nu_1|(B(x,2\varepsilon))}{\mu(B(x,\varepsilon))} \frac{|\nu_2|(B(x,2\varepsilon))}{\mu(B(x,\varepsilon))} \\
&\lesssim \frac{|\nu_1|(B(x,2\varepsilon))}{\mu(B(x,2\varepsilon))} \frac{|\nu_2|(B(x,2\varepsilon))}{\mu(B(x,2\varepsilon))} 
\lesssim M_{\mu}\nu_1(x) M_{\mu}\nu_2(x),
\end{align*}
and completes the proof.
\end{proof}

Finally, we can get the weak type bound for $T_{\sharp}$.
\begin{proposition}\label{prop:weaktypehalf}
Let $\mu$ be a measure of order $m$ on $\R^n$ and $T$ be a bilinear $m$-dimensional SIO. Let
$1 < r, p,q < \infty$ be so that $1/p + 1/q = 1/r$, and suppose we have uniformly on $\varepsilon > 0$ that
$$
\|T_{\mu, \varepsilon}\|_{L^p(\mu) \times L^q(\mu) \to L^r(\mu)} \lesssim 1.
$$
Then we have
$$
\|T_{\sharp}\|_{M(\R^n) \times M(\R^n) \to L^{1/2,\infty}(\mu)} \lesssim 1.
$$
\end{proposition}
\begin{proof}
Fix $\nu_1, \nu_2 \in M(\R^n)$.
It suffices to prove that
$$
\sup_{\delta >0} \sup_{\lambda > 0}  \lambda \mu(\{x \in \R^n\colon\, T_{\sharp,\delta}(\nu_1, \nu_2)(x) > \lambda\})^2 \lesssim \|\nu_1\|\|\nu_2\|.
$$
Fix $\delta  >0$. We know by Proposition \ref{prop:weaktypehalfforT} that
$$
\|T_{\delta}\|_{M(\R^n) \times M(\R^n) \to L^{1/2,\infty}(\mu)} \lesssim 1.
$$
In particular, we have by Proposition \ref{prop:cotbas} that
$$
T_{\sharp, \delta}(\nu_1, \nu_2)(x) \lesssim N_{\mu, 1/4}(T_{\delta}(\nu_1, \nu_2))(x) + M_{\mu}\nu_1(x)M_{\mu}\nu_2(x), \qquad x \in \R^n.
$$

Recall that
$$
\sup_{\lambda > 0} \lambda \mu(\{x \in \R^n\colon\, M_{\mu}\nu_1(x)M_{\mu}\nu_2(x) > \lambda\})^2 \lesssim \|\nu_1\|\|\nu_2\|.
$$
This is probably easiest to see by using the facts that
$$
\|f g \|_{L^{1/2,\infty}(\mu)} \lesssim \|f \|_{L^{1,\infty}(\mu)}\|g \|_{L^{1,\infty}(\mu)}
$$
and $$\|M_{\mu}\|_{M(\R^n) \to L^{1,\infty}(\mu)} \lesssim 1.$$

Thus, it suffices to prove that
$$
\sup_{\lambda > 0} \lambda \mu(\{x \in \R^n\colon\, N_{\mu, 1/4}(T_{\delta}(\nu_1, \nu_2))(x) > \lambda\})^2 \lesssim \|\nu_1\|\|\nu_2\|.
$$
We will next use the easy fact that
$$
\sup_{\lambda > 0} \lambda \mu(\{x \in \R^n\colon\, N_{\mu}\nu(x) > \lambda\}) \le |\nu|(\{x \in \R^n\colon\, N_{\mu}\nu(x) > \lambda\}), \qquad \nu \in M(\R^n).
$$
This sharper form of the weak $(1,1)$ inequality is the only reason why the non-homogeneous, non-centered maximal function $N_{\mu}$ is important in the Cotlar's inequality.
Fix $\lambda > 0$ and set
\begin{align*}
H &= \{x \in \R^n\colon\, N_{\mu, 1/4}(T_{\delta}(\nu_1, \nu_2))(x) > \lambda\} \\
&= \{x \in \R^n\colon\, N_{\mu}(|T_{\delta}(\nu_1, \nu_2)|^{1/4})(x) > \lambda^{1/4}\}.
\end{align*}
We now have
$$
\mu(H) \le \frac{1}{\lambda^{1/4}} \int_H |T_{\delta}(\nu_1, \nu_2)|^{1/4}\,d\mu,
$$
where
$$
\int_H |T_{\delta}(\nu_1, \nu_2)|^{1/4}\,d\mu \lesssim \|\nu_1\|^{1/4} \|\nu_2\|^{1/4} \mu(H)^{1/2}
$$
using again that $$\|T_{\delta}\|_{M(\R^n) \times M(\R^n) \to L^{1/2,\infty}(\mu)} \lesssim 1$$ and the Kolmogorov type
argument from the proof of Proposition \ref{prop:cotbas}. Therefore, we have
$$
\mu(H) \lesssim \bigg(\frac{\|\nu_1\|\|\nu_2\|}{\lambda}\bigg)^{1/2},
$$
and we are done.
\end{proof}

\chapter{Bilinear good lambda method}

In this chapter we aim to prove Theorem \ref{thm:goodlambda} -- a certain very useful good lambda type result.
The proof of this in the bilinear setting turns out not to be very different from the linear setting presented in Theorem 2.22 in \cite{ToBook}.
For the convenience of the reader we give most of the details here.

Before proving the good lambda, we need to recall the following version of Whitney covering especially useful for non-doubling measures.
This is originally from \cite{ToBook}, but the version with small boundary cubes as here appears in \cite{MMT1}.
\begin{lemma} \label{lem:whitney}
If $\Omega\subset\R^n$ is open, $\Omega\neq\R^n$, then $\Omega$
can be decomposed as
$$\Omega = \bigcup_{i\in I} Q_i, $$
where $Q_i$, $i\in I$, are closed dyadic cubes with disjoint interiors such
that for some constants $R>20$ and $D_0\geq1$, depending only on $n$, the following holds:
\begin{itemize}
\item[(i)] $10Q_i \subset \Omega$ for each $i\in I$.
\item[(ii)] $R Q_i \cap \Omega^{c} \neq \varnothing$ for each $i\in I$.
\item[(iii)] For each cube $Q_i$, there are at most $D_0$ cubes $Q_j$
such that $10Q_i \cap 10Q_j \neq \varnothing$. Further, for such cubes $Q_i$, $Q_j$, we have $\ell(Q_i)\sim
\ell(Q_j)$.
\end{itemize}
Moreover, if $t$ is a large enough dimensional constant, $\mu$ is a positive Radon measure on $\R^n$ with
$\mu(\Omega)<\infty$, there is a family of cubes $\{\wt Q_j\}_{j\in S}$, with $S\subset I$, so that
$$Q_j\subset \wt Q_j\subset 1.1 Q_j,$$ satisfying the following:
\begin{itemize}
\item[(a)] Each cube $\wt Q_j$, $j\in S$, is $(9,2D_0)$-doubling and has $t$-small boundary.
\item[(b)] The cubes $\wt Q_j$, $j\in S$, are pairwise disjoint.
\item[(c)]
\begin{equation} \label{bqht22}
\mu\biggl( \,\bigcup_{j\in S} \wt Q_j \biggr) \geq \frac1{8D_0}\,
\mu(\Omega).
\end{equation}
\end{itemize}
\end{lemma}

We are ready for the main result of this chapter.
\begin{theorem}\label{thm:goodlambda}
Let $\mu$ be a measure of order $m$ in $\R^n$.
Let $\beta>0$ and $t > 0$ be big enough numbers, depending only on the dimension $n$, and assume $\theta \in (0,1)$. Suppose for each $(5,\beta)$-doubling cube $Q$
with $t$-small boundary there exists a subset $G_Q \subset Q$  such that
$\mu(G_Q)\geq \theta \mu(Q)$ and $T_{\sharp}\colon M(\R^n) \times M(\R^n) \to L^{1/2,\infty}(\mu \rest G_Q)$ is bounded with a uniform constant independent of $Q$.
Then we have that
$T_{\mu, \sharp}\colon L^p(\mu) \times L^q(\mu) \to L^r(\mu)$ boundedly for all $1 < p, q < \infty$ and $1/2 < r < \infty$ satisfying  $1/r = 1/p + 1/q$
with a constant depending on $r,p,q$ and the preceding constants.
\end{theorem}
\begin{proof}
Let us fix two functions $f, g \in L^1(\mu)$ with compact support. 
For $\lambda > 0$ let
$$
\Omega_{\lambda} = \{T_{\mu, \sharp}(f,g) > \lambda\}.
$$
Because the functions $f$ and $g$ have compact support, $\Omega_\lambda$ is a bounded set, and accordingly of finite $\mu$-measure. Let us also check the fact that it is open, and for this suppose $x \in \Omega_\lambda$. Thus, there exists $\rho_0>0$ so that $T_{\mu, \rho_0}(f,g)(x)> \lambda$. Because the measure $\mu \times \mu$ is locally finite, we can find a slightly larger $\rho> \rho_0$ so that $$T_{\mu, \rho}(f,g)(x)> \lambda$$ and
$$
\mu \times \mu \big(\{(y,z) \colon \max(|x-y|,|x-z|)=\rho\}\big) =0.
$$
Then it follows from the dominated convergence theorem that
$$
T_{\mu, \rho}(f,g)(x') \to T_{\mu, \rho}(f,g)(x), \quad x' \to x,
$$
which shows that $T_{\mu, \sharp}(f,g)(x')> \lambda$ if $|x'-x|$ is small enough. Hence $\Omega_\lambda$ is an open set.

Now we can use
Lemma \ref{lem:whitney} to write
$$
\Omega_{\lambda} = \bigcup_{i \in I} Q_i,
$$
and also to extract the collection $\{\wt Q_j\}_{j\in S}$, where $S\subset I$, so that all the properties of the lemma hold.
For $j \in S$ let us write $P_j = \wt Q_j$. The cubes $P_j$ have $t$-small boundary and are $(9, 2D_0)$-doubling, in particular $(5, 2D_0)$-doubling.
So assuming that the parameter $\beta$ from the assumptions is larger than $2D_0$, we have by assumption that there exists $G_{P_j} \subset P_j$ so that
$$\mu(G_{P_j}) \ge \theta\mu(P_j)$$ and $$T_{\sharp}\colon\, M(\R^n) \times M(\R^n) \to L^{1/2,\infty}(\mu \rest G_{P_j})$$ boundedly with a constant $A$ that is uniform in $j \in S$. For $j \in S$
denote $G_j = G_{P_j}$.

The idea is to prove using the previous cubes that given $\varepsilon, \lambda > 0$ there exists $\delta = \delta(\varepsilon, \theta, A) = \delta(\varepsilon)  > 0$ (
$\theta$ and $A$ are fixed constants from the assumptions) so that
\begin{equation}\label{eq:gl}
\mu(\{x \colon\, T_{\mu, \sharp}(f,g)(x) > (1+\varepsilon)\lambda, \, M_{\mu}^{\mathcal{Q}}f(x)M_{\mu}^{\mathcal{Q}}g(x) \le \delta \lambda\})
\le \Big(1- \frac{\theta}{16D_0} \Big)\mu(\Omega_{\lambda}).
\end{equation}
This is enough to conclude the whole proof by standard considerations, but we shall quickly recall the necessary steps later.

By exploiting the fact that now
\begin{align*}
\mu\Big( \Omega_{\lambda} \setminus \bigcup_{j \in S} P_j\Big) + \sum_{j \in S} \mu(P_j \setminus G_j)
\le \Big( 1 - \frac{\theta}{8D_0}\Big) \mu(\Omega_{\lambda}),
\end{align*}
we are reduced to proving
\begin{equation}\label{eq:form5}
\sum_{j \in S} \mu(\{x \in G_j\colon\, T_{\mu, \sharp}(f,g)(x) > (1+\varepsilon)\lambda, \, M_{\mu}^{\mathcal{Q}}f(x)M_{\mu}^{\mathcal{Q}}g(x) \le \delta \lambda\}) \le \frac{\theta}{16D_0}\mu(\Omega_{\lambda})
\end{equation}
if $\delta = \delta(\varepsilon) > 0$ is small enough.
We will do this by showing for every fixed $j$, if $x \in P_j$ is such that
$$
T_{\mu, \sharp}(f,g)(x) > (1+\varepsilon)\lambda \qquad \textup{and} \qquad M_{\mu}^{\mathcal{Q}}f(x)M_{\mu}^{\mathcal{Q}}g(x) \le \delta \lambda
$$
and $\delta$ is small enough,
then
\begin{equation}\label{eq:GL_claim1}
T_{\mu, \sharp}(f1_{2P_j}, g1_{2P_j})(x) > \frac{\varepsilon}{2}\lambda.
\end{equation}
This implies \eqref{eq:form5} (for $\delta(\varepsilon)$ small enough) by the fact that  $$T_{\sharp}\colon\, M(\R^n) \times M(\R^n) \to L^{1/2,\infty}(\mu \rest G_{P_j}).$$
This calculation is done almost exactly as in the linear case, but let us quickly check it. So we assume the above pointwise bound.
For a fixed $j \in S$ notice that now
\begin{align*}
\mu(\{x \in G_j\colon\, &T_{\mu, \sharp}(f,g)(x) > (1+\varepsilon)\lambda, \, M_{\mu}^{\mathcal{Q}}f(x)M_{\mu}^{\mathcal{Q}}g(x) \le \delta \lambda\})  \\
& \le \mu (\{ x \in G_j\colon\, T_{\mu, \sharp}(f1_{2P_j}, g1_{2P_j})(x) > \varepsilon \lambda / 2\}) \\
&\le \Big( \frac{2A}{\varepsilon\lambda} \Big)^{1/2} \Big( \int_{2P_j} |f|\,d\mu\Big)^{1/2}  \Big( \int_{2P_j} |g|\,d\mu\Big)^{1/2}.
\end{align*}
We can assume that there exists $x_0 \in P_j$ such that $$M_{\mu}^{\mathcal{Q}}f(x_0)M_{\mu}^{\mathcal{Q}}g(x_0)  \le \delta \lambda,$$ and estimate
\begin{align*}
 \Big( \int_{2P_j} |f|\,&d\mu\Big)^{1/2}  \Big( \int_{2P_j} |g|\,d\mu\Big)^{1/2}  \\ 
 &\le  (\mu(10Q_j) M_{\mu}^{\mathcal{Q}}f(x_0))^{1/2} (\mu(10Q_j) M_{\mu}^{\mathcal{Q}}g(x_0))^{1/2} 
 \le  (\delta \lambda)^{1/2} \mu(10Q_j).
 \end{align*}
So we have
\begin{align*}
\sum_{j \in S} &\mu(\{x \in G_j\colon\, T_{\mu, \sharp}(f,g)(x) > (1+\varepsilon)\lambda, \, M_{\mu}^{\mathcal{Q}}f(x)M_{\mu}^{\mathcal{Q}}g(x) \le \delta \lambda\}) \\
&\le \Big( \frac{2A\delta }{\varepsilon} \Big)^{1/2} \sum_{j \in I} \mu(10Q_j) \le  D_0\Big( \frac{2A\delta }{\varepsilon} \Big)^{1/2}  \mu(\Omega_{\lambda})
\le \frac{\theta}{16D_0}\mu(\Omega_{\lambda})
\end{align*}
for $\delta$ small enough.

It only remains to prove the pointwise lower bound.
So suppose $x \in P_j$ is such that 
$$
T_{\mu, \sharp}(f,g)(x) > (1+\varepsilon)\lambda \qquad \textup{and} \qquad M_{\mu}^{\mathcal{Q}}f(x)M_{\mu}^{\mathcal{Q}}g(x) \le \delta \lambda.
$$
Using the $\tilde T$ notation there holds
\begin{equation*}
\begin{split}
T_{\mu,\sharp} (f,g)(x)
\leq
T_{\mu,\sharp} (1_{2P_j}f,1_{2P_j}g)(x) 
+ \tilde{T}_{\mu, \sharp} \big(1_{(2P_j \times 2P_j)^c} f \otimes g\big)(x).
\end{split}
\end{equation*}
Therefore, \eqref{eq:GL_claim1} follows  once we show that the latter term on the right is at most $(1+\varepsilon/2)\lambda$. To do this, we fix an arbitrary $\rho_0 > 0$ and show that $$\tilde{T}_{\mu, \rho_0} \big(1_{(2P_j \times 2P_j)^c} f \otimes g\big)(x) \leq (1+ \varepsilon/2)\lambda.$$
Notice first that
\begin{equation*}
\begin{split}
 \iint \displaylimits_{2 \diam (RP_j) \geq\max(|x-y|,|x-z|)}  &\frac{ 1_{(2P_j \times 2P_j)^c}(y,z) |f(y)g(z)|}{ (|x-y|+|x-z|)^{2m}}  \, d\mu(y) \, d\mu(z) \\
& \lesssim   \iint \displaylimits_{\substack{|x-y| \leq 2 \diam(RP_j) \\ |x-z| \leq 2 \diam(RP_j)}} \frac{|f(y)g(z) |}{\ell(P_j)^{2m}} \, d\mu(y) \, d\mu(z) \\
& \lesssim M_{\mu}^{\mathcal{Q}}f(x)M_{\mu}^{\mathcal{Q}}g(x)
\leq \delta \lambda.
\end{split}
\end{equation*}
Thus, it is enough the define $\rho:= \max(\rho_0, 2 \diam(RP_j))$ and consider 
$$
\tilde{T}_{\mu, \rho} \big(1_{(2P_j \times 2P_j)^c} f \otimes g \big)(x)
=\tilde{T}_{\mu, \rho} ( f \otimes g )(x)
=T_{\mu, \rho} ( f , g )(x),
$$
where the  first equality holds because 
$$ 
\big\{ (y,z) \colon \max(|x-y|,|x-z|)>\rho \big\}
\subset (2P_j \times 2P_j)^c.
$$

By the properties of the Whitney cubes there exists a point $x' \in RP_j \cap \Omega_\lambda^c$, and we can estimate
\begin{equation*}
\begin{split}
\big|T_{\mu, \rho} ( f , g )(x) \big|
&\leq \big|T_{\mu, \rho} ( f , g )(x) - T_{\mu, \rho} ( f , g )(x')\big| 
+\big |T_{\mu, \rho} ( f , g )(x') \big | \\
&\leq \big|T_{\mu, \rho} ( f , g )(x) - T_{\mu, \rho} ( f , g )(x')\big|
+ \lambda.
\end{split}
\end{equation*}
Thus,  it suffices to estimate the difference, which can be dominated with
\begin{equation}\label{eq:GL_difference}
\begin{split}
&\Big| \iint \displaylimits_{\max(|x-y|,|x-z|)>\rho} \big(K(x,y,z)-K(x',y,z) \big)  f(y) g(z) \, d\mu(y) \, d\mu(z)\Big| \\
&+ \Big| \iint \displaylimits_{\max(|x-y|,|x-z|)>\rho} K(x',y,z)  f(y) g(z) \, d\mu(y) \, d\mu(z)-T_{\mu, \rho} ( f , g )(x') \Big|.
\end{split}
\end{equation}

Applying kernel estimates, the first term in \eqref{eq:GL_difference} can be dominated with
\begin{equation*}
\iint \displaylimits_{\max(|x-y|,|x-z|)>\rho} \frac{|x-x'|^\alpha |f(y) g(z)|}{ (|x-y|+|x-z|)^{2m+\alpha}} \, d\mu(y) \, d\mu(z)
\lesssim M_{\mu}^{\mathcal{Q}}f(x)M_{\mu}^{\mathcal{Q}}g(x)
\leq \delta \lambda.
\end{equation*}
Notice that the symmetric difference
$$
\big\{ (y,z) \colon \max(|x-y|,|x-z|)>\rho \big\} \triangle \big\{ (y,z) \colon \max(|x'-y|,|x'-z|)>\rho \big\}
$$
is contained in the set 
$$
\big\{ (y,z) \colon \max(|x-y|,|x-z|) \sim  \max(|x'-y|,|x'-z|) \sim \rho\}.
$$
Thus, the second term in \eqref{eq:GL_difference} can be dominated with
$$
\iint \displaylimits_{\max(|x-y|,|x-z|) \sim \rho} \frac{|f(y) g(z)|}{ (|x-y|+|x-z|)^{2m}} \, d\mu(y) \, d\mu(z)
\lesssim  M_{\mu}^{\mathcal{Q}}f(x)M_{\mu}^{\mathcal{Q}}g(x)
\leq \delta \lambda.
$$

Combining the above arguments we have shown that there exists a constant $C$ such that
$$
\tilde{T}_{\mu, \rho_0} \big(1_{(2P_j \times 2P_j)^c} f \otimes g\big)(x)
\leq C \delta \lambda+ \lambda.
$$
Hence, if $\delta(\varepsilon)$ is chosen to be small enough, it is seen that $$C \delta \lambda+ \lambda \leq (1+ \varepsilon/2)\lambda,$$ and accordingly \eqref{eq:GL_claim1} is satisfied.

We have shown \eqref{eq:gl}. The claim follows from this by standard arguments, but requires a moderate amount of approximation.
For convenience, we outline these details now. Let $1 < p, q < \infty$ and $1/2 < r < \infty$ satisfy  $1/r = 1/p + 1/q$.
Since we do not know that $\| T_{\mu, \sharp} (f,g) \| _{L^{r}(\mu)} < \infty$, we first do the following.
Define $$h_k = \inf(k, T_{\mu, \sharp} (f,g)), \qquad k \ge 1.$$ Suppose $R > 0$ is such that $B(0,R)$ contains the supports of $f$ and $g$.
Then for $x \in B(0,2R)^c$ there holds that
$$
h_k(x) \lesssim \frac{\|f\|_{L^1(\mu)} \|g\|_{L^1(\mu)}}{|x|^{2m}}.
$$
Notice that $2mr > m$ so that $$\int_{B(0,2R)^c} |x|^{-2mr} \,d\mu(x) < \infty.$$ It follows that $\| h_k \| _{L^{r}(\mu)} < \infty$.
Moreover, the good lambda inequality \eqref{eq:gl} is also true with $T_{\mu, \sharp} (f,g)$ replaced by $h_k$ everywhere 
(this follows from \eqref{eq:gl} directly using the definition of $h_k$). Using this good lambda inequality, the fact that
$$\|M_{\mu}^{\mathcal{Q}}fM_{\mu}^{\mathcal{Q}}g\|_{L^r(\mu)} \lesssim \|f\|_{L^p(\mu)} \|g\|_{L^q(\mu)}$$ and $\| h_k \| _{L^{r}(\mu)} < \infty$, we easily see using the
distributional formula for the $L^r(\mu)$ norm that
$$
\|h_k\|_{L^r(\mu)} \lesssim \|f\|_{L^p(\mu)} \|g\|_{L^q(\mu)}.
$$
Letting $k \to \infty$ we get
\begin{equation}\label{eq:coml1}
\| T_{\mu, \sharp} (f,g) \| _{L^{r}(\mu)} \lesssim \| f \|_{L^p(\mu)} \| g\|_{L^q(\mu)}.
\end{equation}
Recall that $f$ and $g$ were $L^1(\mu)$ functions with compact supports.
Let us now extend this to all functions $f \in L^p(\mu)$ and $g \in L^q(\mu)$.
To this end, choose some arbitrary real numbers $M, \rho>0$. Using \eqref{eq:basicbound} we see that there exists a constant $C(M, \rho)$ so that
\begin{equation}\label{eq:GL_apriori}
\|1_{B(0,M)} T_{\mu, \sharp, \rho} (f,g) \| _{L^{r}(\mu)} \leq C(M, \rho) \| f \|_{L^p(\mu)} \| g\|_{L^q(\mu)}, 
\quad f \in L^p(\mu),g \in L^q(\mu).
\end{equation}
Approximating with $L^1(\mu)$ functions with compact supports and using \eqref{eq:coml1} we see that the above holds
with a constant independent of $M$ and $\rho$.  Finally, letting $M \to \infty$ and $\rho \to 0$ we get \eqref{eq:coml1} for all $f \in L^p(\mu)$ and $g \in L^q(\mu)$.
\end{proof}

\chapter{Proof of the main theorem}\label{sec:syn}

We are ready to prove our main theorem.
\begin{proof}[Proof of Theorem \ref{thm:main}]
Fix a $(5,b)$-doubling cube $Q_0$ with $t$-small boundary, and set $\sigma = \mu \rest Q_0$.
Now, the measure $\sigma$ is of order $m$, $$\sigma(\R^n \setminus Q_0) = 0$$ and for $t_0 := tb$ we have
$$
\sigma(\{x \in Q_0\colon d(x, \partial Q_0) \le \lambda \ell(Q_0)\}) \le t_0 \lambda \sigma(Q_0)
$$
for all $\lambda > 0$. Notice that $b_i \in L^{\infty}(\sigma)$, $i = 1,2,3$, and
$$
|\langle b_i \rangle_Q^{\sigma}| =  |\langle b_i \rangle_Q^{\mu}| \gtrsim 1 \qquad \textup{for all cubes } Q \subset Q_0.
$$
Moreover, for all $1$-Lipschitz functions $\Phi \colon \R^n \to [0,\infty)$, truncation parameters $\delta > 0$ and for all cubes $Q$ satisfying $5Q \subset Q_0$ we have
$$
|\langle T_{\sigma, \Phi, \delta}(1_Qb_1, 1_Qb_2), 1_Qb_3\rangle_{\sigma}|
= |\langle T_{\mu, \Phi, \delta}(1_Qb_1, 1_Qb_2), 1_Qb_3\rangle_{\mu}| 
 \lesssim \mu(5Q) = \sigma(5Q).
$$
Corollary \ref{cor:cot} gives that
$$
\int_{Q_0 } [S_{\sigma, \sharp}(b, b')]^{s/2} \,d\sigma
= \int_{Q_0 } [S_{\mu, \sharp}(1_{Q_0}b, 1_{Q_0}b')]^{s/2} \,d\mu \lesssim \mu(5Q_0) \sim \mu(Q_0) = \sigma(Q_0)
$$
for all the choices $$(S, b, b') \in \{(T, b_1, b_2), (T^{1*}, b_3, b_2), (T^{2*}, b_1, b_3)\}.$$

We are now in the position to use Theorem \ref{thm:bigpieceTb} to find a set $G \subset Q_0$ so that $$\mu(G) = \sigma(G) \sim \sigma(Q_0) = \mu(Q_0)$$
and uniformly over $\varepsilon > 0$ it holds that
$$
\|T_{\mu \rest G, \varepsilon}\|_{L^4(\mu \rest G) \times L^4(\mu \rest G) \to L^2(\mu \rest G)} = \|T_{\sigma \rest G, \varepsilon}\|_{L^4(\sigma \rest G) \times L^4(\sigma \rest G) \to L^2(\sigma \rest G)} \lesssim 1.
$$
It follows from Proposition \ref{prop:weaktypehalf} that
$$
\|T_{\sharp}\|_{M(\R^n) \times M(\R^n) \to L^{1/2,\infty}(\mu \rest G)} \lesssim 1.
$$
Since $Q_0$ was an arbitrary $(5,b)$-doubling cube with $t$-small boundary, the good lambda method gives that
for all $1 <p, q < \infty$ and $1/2 < r < \infty$ satisfying
$1/p + 1/q = 1/r$ we have that
$$
\|T_{\mu, \sharp}\|_{L^p(\mu) \times L^q(\mu) \to L^r(\mu)} \lesssim 1.
$$
Therefore, we are done.
\end{proof}

\begin{proof}[Proof of Corollary \ref{cor:MainThm}]
Notice that in the proof of the big piece $Tb$ theorem we used the weak boundedness property only in cubes that have $t_0$-small boundary  for some large $t_0$ (which only depends on some constants appearing in the statement). This implies that the main theorem, Theorem \ref{thm:main}, holds in the stronger form with weak boundedness only in cubes that have $t_0$-small boundary.

Therefore, to prove Corollary \ref{cor:MainThm} we need to verify the weak boundedness property in cubes that have $t_0$-small boundary.
Fix one such cube $Q$.
Then Corollary \ref{cor:cot} combined with Remark \ref{rem:StrongCorCot} implies that
\begin{equation*}
\int_Q |T_{\mu, \sharp}(1_Qb_1, 1_Qb_2)|d \mu \lesssim  \mu(5Q).
\end{equation*}

Let $\Phi$ be any $1$-Lipschitz function. Equation \eqref{eq:suppmain} gives that
\begin{equation*}
\begin{split}
T_{\mu, \Phi,\sharp}(1_Qb_1,1_Qb_2)(x) 
&\le T_{\mu, \sharp}(1_Qb_1,1_Qb_2)(x) + CM_{\mu} (1_Qb_1)(x) M_{\mu} (1_Qb_2)(x) \\
& \le T_{\mu, \sharp}(1_Qb_1,1_Qb_2)(x) +C.
\end{split}
\end{equation*}
Let $\delta>0$. Now we have 
\begin{equation*}
\begin{split}
|\langle T_{\mu, \Phi, \delta}(1_Qb_1, 1_Qb_2), 1_Qb_3\rangle_{\mu}|
&\lesssim \int_Q T_{\mu, \Phi, \sharp}(1_Qb_1, 1_Qb_2) d \mu \\
& \lesssim \int_Q T_{\mu, \sharp}(1_Qb_1, 1_Qb_2) d \mu+\mu(Q) \\
& \lesssim \mu(5Q).
\end{split} 
\end{equation*}
This concludes the proof.
\end{proof}

\chapter{Weakening the kernel estimates: modified Dini-condition}\label{sec:dini}

In this chapter we describe how, as a byproduct of our new summation methods, we can do all of our dyadic summations
in the core $Tb$ argument (proof of Theorem \ref{thm:bigpieceTb}) under some weaker regularity conditions on the underlying bilinear kernel $K$.
This means that one replaces the H\"older conditions like
$$
|K(x,y,z)-K(x',y,z)| \le C_K\frac{|x-x'|^{\alpha}}{(|x-y| + |x-z|)^{2m+\alpha}}
$$
whenever $|x-x'| \le \max(|x-y|, |x-z|)/2$, with some more general conditions (the size condition is unaltered).
In practice, one replaces the modulus of continuity $t \mapsto t^{\alpha}$, $\alpha \in (0,1]$, with some more general
modulus of continuity.

The remarks in this chapter are motivated by the recent paper by A. Grau de la Herr\'an and T. Hyt\"onen \cite{GH}.
There they show that (linear) $T1$ theorems can be proved, even in the non-homogeneous situation, assuming quite
weak modified Dini-conditions (the condition assumed is the sharpest known even in the homogeneous situation). Previously,
non-homogeneous $T1$ theory had been only developed with the modulus of continuity $t \mapsto t^{\alpha}$. In fact, this does
require some thought as previously everything was based on some summation arguments like the ones alluded to in Remark
\ref{rem:matrixlemma}, and some suitably small choice of the parameter $\gamma$ in the definition of goodness. The point we want to make here
is that we can extremely straightforwardly consider more general kernels in $Tb$ summation arguments as a byproduct of the fact we no longer operate as in Remark \ref{rem:matrixlemma}.
Our observation is completely independent of the considerations in \cite{GH}, and we cannot (at least without working harder) quite reach the sharpness considered there, but quite close.

Let us get to the details. Suppose $\psi$ is a modulus of continuity, i.e let $\psi\colon [0, \infty) \to [0, \infty)$ be an increasing function such that $\psi(0)=0$ and $\psi(s+t) \le \psi(s)+\psi(t)$.
The regularity demanded of $\psi$ will be of the modified Dini-type i.e.
\begin{equation}\label{eq:bDini}
\int_0^1 \psi(t) \Big( 1 + \log \frac{1}{t} \Big)^{\beta} \frac{dt}{t} < \infty
\end{equation}
for some $\beta$ (the case $\beta = 0$ is the Dini-condition). We will need $\beta=1$, while $\beta = 1/2$ is the sharpest condition known to be enough for (linear) $T1$ summation arguments (for the non-homogeneous result see \cite{GH}, and for the history
of the homogeneous results see the introduction of the said article). So for what follows we assume that
$$
\int_0^1 \psi(t) \Big( 1 + \log \frac{1}{t} \Big) \frac{dt}{t} < \infty,
$$
and that
$$
K\colon (\R^n \times \R^n \times \R^n) \setminus \Delta \to \C, \qquad \Delta := \{(x,y,z) \in \R^n \times \R^n \times \R^n\colon\, x = y = z\},
$$
is an $m$-dimensional bilinear Calder\'on--Zygmund kernel with a modulus of continuity $\psi$. More precisely, this means that for some constant $C_K < \infty$ we have
$$
|K(x,y,z)| \le \frac{C_K}{(|x-y|+|x-z|)^{2m}},
$$
$$
|K(x,y,z)-K(x',y,z)| \le \frac{C_K}{(|x-y|+|x-z|)^{2m}}\psi\Big(\frac{|x-x'|}{|x-y|+|x-z|}\Big)
$$
whenever $|x-x'| \le \max(|x-y|,|x-z|)/2$,
$$
|K(x,y,z)-K(x,y',z)| \le  \frac{C_K}{(|x-y|+|x-z|)^{2m}} \psi\Big(\frac{|y-y'|}{|x-y|+|x-z|}\Big)
$$
whenever $|y-y'| \le \max(|x-y|, |x-z|)/2$, and
$$
|K(x,y,z)-K(x,y,z')| \le \frac{C_K}{(|x-y|+|x-z|)^{2m}} \psi\Big(\frac{|z-z'|}{|x-y|+|x-z|}\Big)
$$
whenever $|z-z'| \le \max(|x-y|, |x-z|)/2$.

We now comment on the modifications required in the proof of Theorem \ref{thm:bigpieceTb}. The key auxiliary estimate using the continuity of the kernel
in that proof is Lemma \ref{lem:zeroave_sep}. With our current kernels it takes the following form.

\begin{lemma}\label{lem:zeroave_sep_D}
Suppose $A \subset \R^n$ is a bounded set and $h_0$ is a function supported on $A$ such that $\int h_0 \,d\mu =0$.
Suppose also that $t \ge 2$ and $B \subset \R^{2n}$ is a set satisfying
$$
B \subset \{(y,z) \in \R^{2n}\colon\, \inf_{x \in A} \max(|x-y|, |x-z|) \ge t d(A)\}.
$$
Then we have for $f_0,g_0 \in L^1_{\text{loc}}(\mu)$ that
$$
\big |\bla \tilde{T}_{\mu,\Phi}(1_{B} f_0 \otimes g_0), h_0 \bra_\mu \big|
\lesssim \int_0^{t^{-1}} \psi(s) \frac{ds}{s} \int M_{\mu, m}(f_0,g_0) | h_0 | \,d\mu.
$$
\end{lemma}

\begin{proof}
Using the zero average of $h_0$ and the $x$-continuity of the kernel we have
\begin{equation*}
\begin{split}
\big |\bla \tilde{T}_{\mu,\Phi}&(1_{B} f_0 \otimes g_0), h_0 \bra_\mu \big| \\
&\lesssim \iiint \frac{\psi \Big(\frac{d(A)}{|x-y|+|x-z|}\Big)}{(|x-y|+|x-z|)^{2m}}1_B(y,z)|f_0(y)g_0(z)h_0(x)| \, d\mu(y) \, d\mu(z) \, d\mu(x) \\
& \lesssim \iiint \frac{\psi \Big(\frac{d(A)}{td(A)+|x-y|+|x-z|}\Big)}{(td(A)+|x-y|+|x-z|)^{2m}}|f_0(y)g_0(z)h_0(x)| \, d\mu(y) \, d\mu(z) \, d\mu(x). 
\end{split}
\end{equation*}
For $x \in \R^n$ one sees that
\begin{equation*}
\begin{split}
\iint \frac{\psi \Big(\frac{d(A)}{td(A)+|x-y|+|x-z|}\Big)}{(td(A)+|x-y|+|x-z|)^{2m}}&|f_0(y)g_0(z)| \, d\mu(y) \, d\mu(z) \\
&\lesssim \sum_{k=0}^\infty \psi(2^{-k}t^{-1}) \calM_{\mu,m} (f_0,g_0)(x).
\end{split}
\end{equation*}
To conclude the proof, it remains to notice that
$$
\sum_{k=0}^\infty \psi(2^{-k}t^{-1})
\sim \int_0^{t^{-1}} \psi(s) \frac{ds}{s}.
$$
\end{proof}

To demonstrate how this is used in the summation arguments, where we previously used Lemma \ref{thm:bigpieceTb},
we estimate (a part of) the separated sum from the $Tb$ argument. All the other applications are completely analogous. We estimate the term
\begin{equation*}
\begin{split}
\sum_{\substack{K \in \calD_0(\omega_3) \\ \ell(K) < 2^{u_0}}} &\sum_{l=0}^{u_0-\log_2\ell(K)}
\Big| \Big \langle T_{\mu, \Phi}  ( \varphi_{K,l} , E^2_{2^{l-1}\ell(K)}g_\omega ), \Delta^3_K h_\omega \Big \rangle_\mu \Big| \\
&=\sum_{l=0}^\infty \sum_{\substack{K \in \calD_0(\omega_3) \\ \ell(K) < 2^{u_0} \\ \ell(K) \le 2^{u_0-l}}}
\Big| \Big \langle T_{\mu, \Phi}  ( \varphi_{K,l} , E^2_{2^{l-1}\ell(K)}g_\omega ), \Delta^3_K h_\omega \Big \rangle_\mu \Big|
\end{split}
\end{equation*}
from the subsection "The separated part" in the proof of Theorem \ref{thm:bigpieceTb}.

For a fixed $l$ Lemma \ref{lem:zeroave_sep_D} gives that
\begin{equation*}
\begin{split}
\sum_{\substack{K \in \calD_0(\omega_3) \\ \ell(K) < 2^{u_0} \\ \ell(K) \le 2^{u_0-l}}}
&\Big| \Big \langle T_{\mu, \Phi}  ( \varphi_{K,l} , E^2_{2^{l-1}\ell(K)}g_\omega ), \Delta^3_K h_\omega \Big \rangle_\mu \Big| \\
&\lesssim \sum_{\substack{K \in \calD_0(\omega_3) \\ \ell(K) < 2^{u_0} \\ \ell(K) \le 2^{u_0-l}}}
\int_0^{2^{-l(1-\gamma)}} \psi(s) \frac{ds}{s}
\int M_{\mu,m}( D_{2^l \ell(K)}^1 f_\omega, M_{\mu, \calD_0(\omega_2)} g_\omega ) | \Delta_K^3h_\omega | \,d\mu \\
& \lesssim \int_0^{2^{-l(1-\gamma)}} \psi(s) \frac{ds}{s}
\| f \|_{L^p(\mu)}\| g \|_{L^q(\mu)}\| h \|_{L^{r'}(\mu)}.
\end{split}
\end{equation*}
To conclude the estimate it remains to consider
$$
\sum_{l=0}^\infty \int_0^{2^{-l(1-\gamma)}} \psi(s) \frac{ds}{s} 
= \int_0^1 \sum_{l=0}^\infty 1_{(0,2^{-l(1-\gamma)}]}(s) \psi(s) \frac{ds}{s}. 
$$
Notice that if  $2^{-l(1-\gamma)} \ge s$ for $s \in (0,1)$ then 
$$
l \le \frac{1}{(1-\gamma) \log 2} \log 1/s \lesssim \log (1/s).
$$
This shows that
$$
\sum_{l=0}^\infty \int_0^{2^{-l(1-\gamma)}} \psi(s) \frac{ds}{s} 
\lesssim \int_0^1 \psi(s) \Big( 1 + \log \frac{1}{s} \Big) \frac{ds}{s},
$$
which concludes our demonstration.

\subsection*{Standard estimates}
It is well known that the Dini condition  (case $\beta=0$ in \eqref{eq:bDini}) is sufficient for much of the classical theory of singular integrals. 
Apart from some summations in the  $Tb$ argument (discussed above), this is also the case for the current paper. That is, everything else in this paper can be made to
work with just the Dini--condition, in particular all the results of this paper work with the modified Dini--condition with $\beta = 1/2$.
To see this, in most of the estimates where we have used the H\"older continuity of the kernel 
one can just use the Dini continuity instead and apply estimates of the form
\begin{equation}\label{eq:DecayDini}
\int
 \frac{\psi\Big(\frac{s}{s+|x-y|}\Big)}{(s+|x-y|)^{m}} \, d \mu (y)  
 \lesssim \int_0^1 \psi(t) \frac{d t}{t}, \quad x \in \R^n, s>0.
 \end{equation}

However, there is one place where we need to argue a little differently with the Dini--condition, namely, 
in a part of the proof of Proposition \ref{prop:weaktypehalfforT} (which stated that $T_{\mu, \varepsilon} \colon L^p \times L^q \to L^r$ 
implies $T_\varepsilon \colon M(\R^n) \times M(\R^n) \to L^{1/2, \infty}(\mu)$). 
We explicitly show this modification here, because we are unaware of a reference covering the details in our non-homogeneous situation.
For a homogeneous version of Proposition \ref{prop:weaktypehalfforT} with the Dini condition, see \cite{LP}.

Suppose we are in the set-up of Proposition \ref{prop:weaktypehalfforT}. Let $I$ be the term defined in Equation \eqref{eq:main cases}, that is,
$$
I=\mu\Big(\Big\{x \in \R^n \setminus \scrA \colon \sum_i 1_{(2\Ri)^c}(x) \Big|T_{\varepsilon}\Big(\bei, \sum_{j \in \calJ_i}\bej\Big)(x)\Big| >\lambda/16\Big\}\Big).
$$
We show that $I \lesssim \lambda^{-1/2}$ assuming that the kernel of $T_\varepsilon$ satisfies only the Dini regularity condition. 
The other estimates in the proof of Proposition \ref{prop:weaktypehalfforT} where we used H\"older continuity of the kernel can be directly handled with \eqref{eq:DecayDini}.

First we record the property of the Calder\'on-Zygmund decomposition that
\begin{equation}\label{eq:SmallMaximal}
\sup_{r>0} \frac{|\eta|(B(x,r))}{r^m} \lesssim \lambda^{1/2}, \quad \text{for all } x \in \R^n \setminus \bigcup_j 2Q_{2,j}.
\end{equation}
To see this, fix some $x \in \R^n \setminus \bigcup_j 2Q_{2,j}$ and $r>0$.
Recall that by the property \eqref{cd3} of the Calder\'on-Zygmund decomposition $1_{\R^n \setminus \bigcup_j \Qj} \eta = f_2 d \mu$ for a function $f_2$ with $\| f_2 \|_{L^\infty(\mu)} \le \lambda^{1/2}.$
If $B(x,r) \cap \bigcup_j \Qj = \emptyset$, then 
$$
|\eta|(B(x,r)) = \int_{B(x,r)} |f_2| \, d \mu \lesssim \| f_2\|_{L^\infty(\mu)} \mu(B(x,r)) \lesssim \lambda^{1/2} r^m.
$$

Suppose then $B(x,r) \cap \Qj  \not = \emptyset$ for some $j$. Because $x \not \in 2\Qj$ there holds $r \ge \ell(\Qj)/2$,
whence 
$$
B(x,r) \subset Q(c_{\Qj}, Cr)
$$
for some absolute constant $C$ that we can assume to satisfy $C \ge 3$. Recall that $Q(x,r)$ denotes the cube with center $x$ and sidelength $2r$.
Since $Cr \ge 3 \ell(\Qj)/2$, the property \eqref{cd2} of 
the Calder\'on-Zygmund decomposition 
gives that 
$$
|\eta|(Q(c_{\Qj}, Cr)) \lesssim \lambda^{1/2} \mu(Q(c_{\Qj}, 2Cr)) \lesssim \lambda^{1/2}r^m.
$$
Combining these shows that $|\eta| (B(x,r)) \lesssim \lambda^{1/2}r^m$, and so we have proved \eqref{eq:SmallMaximal}.

Now we turn to estimate $I$ and begin with
\begin{equation*}
I \le \lambda^{-1} \sum_i 
\int_{\R^n \setminus \scrA} 1_{(2\Ri)^c} \sum_j | T_{\varepsilon}(\bei, \bej)| \, d \mu.
\end{equation*}
Fix some $i$ for the moment and suppose $x \in (2\Ri \cup \scrA )^c$ is such that $d(x,\Ri) > \varepsilon$. Recall that
$\| \bej \| \lesssim |\eta|(\Qj)$ for every $j$.
We see by applying the zero average of $\bei$  that
\begin{equation*}
\begin{split}
|T_{\varepsilon}(\bei, \bej)(x)|
&\lesssim \| \bei \| \| \bej \| \frac{\psi\Big(\frac{\ell(\Ri)}{|x-c_{\Ri}|+|x-c_{\Rj}|}\Big)}{(|x-c_{\Ri}|+|x-c_{\Rj}|)^{2m}} \\
& \lesssim \| \bei \| \int \frac{\psi\Big(\frac{\ell(\Ri)}{|x-c_{\Ri}|}\Big) }{(|x-c_{\Ri}|+|x-z|)^{2m}}  1_{\Qj}\, d |\eta|(z),
\end{split}
\end{equation*}
where in the last step we used the fact that $x \in \scrA^c \subset (2\Qj)^c$.
Summing this over $j$ and using the fact that the cubes $\{\Qj\}_j$ have bounded overlap shows that
\begin{equation*}
\begin{split}
\sum_j |T_{\varepsilon}(\bei, \bej)(x)|
&\lesssim \| \bei \| \int \frac{\psi\Big(\frac{\ell(\Ri)}{|x-c_{\Ri}|}\Big) }{(|x-c_{\Ri}|+|x-z|)^{2m}} \, d |\eta|(z) \\
& \lesssim \| \bei \|  \frac{\psi\Big(\frac{\ell(\Ri)}{|x-c_{\Ri}|}\Big) }{|x-c_{\Ri}|^{m}} \lambda^{1/2},
\end{split}
\end{equation*}
where we applied \eqref{eq:SmallMaximal}. Thus
\begin{equation*}
\begin{split}
\sum_i 
\int_{\R^n \setminus \scrA} &1_{(2\Ri)^c\cap \{d(\cdot,\Ri) > \varepsilon\}} \sum_j | T_{\varepsilon}(\bei, \bej)| \, d \mu \\
& \lesssim  \lambda^{1/2} \sum_i \| \bei \| \int_{(2 \Ri)^c} \frac{\psi\Big(\frac{\ell(\Ri)}{|x-c_{\Ri}|}\Big) }{|x-c_{\Ri}|^{m}}
d \mu(x)
\lesssim \lambda^{1/2},
\end{split}
\end{equation*}
since $\sum_j \| \bei \| \lesssim 1$.

Suppose now $x \in (2\Ri \cup \scrA)^c$ is such that $d(x, \Ri) \le \varepsilon$. 
This implies that $\ell(\Ri) \le 2\varepsilon$, and so $x \in B(c_{\Ri}, C\varepsilon)$.
The size estimate of the kernel gives
$$
|T_{\varepsilon}(\bei, \bej)(x)| 
\lesssim \frac{\| \bei \| \| \bej \|}{(\varepsilon +|x-c_{\Ri}|+|x-c_{\Rj}|)^{2m}},
$$
which can similarly as above be summed over $j$ to conclude that
$$
\sum_j |T_{\varepsilon}(\bei, \bej)(x)|
\lesssim   \frac{\| \bei \| \lambda^{1/2}}{\varepsilon^{m}}.
$$
These observations show that
\begin{equation*}
\begin{split}
\sum_i 
\int_{\R^n \setminus \scrA}  & 1_{(2\Ri)^c\cap \{d(\cdot,\Ri) \le \varepsilon\}} \sum_j | T_{\varepsilon}(\bei, \bej)| \, d \mu \\
& \lesssim \sum_i \int 1_{B(c_{\Ri}, C \varepsilon)} \frac{\| \bei \| \lambda^{1/2}}{\varepsilon^{m}} \, d \mu
\lesssim \lambda^{1/2},
\end{split}
\end{equation*}
which concludes the proof of $I \lesssim \lambda^{-1/2}$.

\chapter{Briefly about square functions}\label{sec:SF}

There is a natural class of bilinear square functions for which the above theory can also be developed. Let us formulate the setting and a boundedness theorem for these, and very briefly 
discuss the proof -- and why it is significantly simpler than the Calder\'on--Zygmund case from above.

For each $t > 0$ and for two complex measures $\nu_1$ and $\nu_2$ on $\R^n$ we define
$$
\theta_t(\nu_1,\nu_2)(x) = \int_{\R^n} \int_{\R^n} s_t(x,y,z) \,d\nu_1(y) \,d\nu_2(z), \qquad x \in \R^n.
$$
The above integral converges absolutely if $\nu_1, \nu_2 \in M(\R^n)$ and the kernel satisfies \eqref{eq:sizeSF} from below.
For some $m, \alpha >0$ the kernels $s_t \colon \R^n \times \R^n \times \R^n \to \C$ are assumed to satisfy the size condition
\begin{equation}\label{eq:sizeSF}
|s_t(x,y,z)| \lesssim \frac{t^{2\alpha}}{(t+|x-y|)^{m+\alpha}(t+|x-z|)^{m+\alpha}},
\end{equation}
the $x$-H\"older condition
\begin{equation}\label{eq:xholSF}
|s_t(x,y,z) - s_t(x',y,z)| \lesssim \frac{t^{\alpha}|x-x'|^{\alpha}}{(t+|x-y|)^{m+\alpha}(t+|x-z|)^{m+\alpha}}
\end{equation}
whenever $|x-x'| < t/2$, the $y$-H\"older condition
\begin{equation}\label{eq:yholSF}
|s_t(x,y,z) - s_t(x,y',z)| \lesssim \frac{t^{\alpha}|y-y'|^{\alpha}}{(t+|x-y|)^{m+\alpha}(t+|x-z|)^{m+\alpha}}
\end{equation}
whenever $|y-y'| < t/2$, and the $z$-H\"older condition
\begin{equation}\label{eq:zholSF}
|s_t(x,y,z) - s_t(x,y,z')| \lesssim \frac{t^{\alpha}|z-z'|^{\alpha}}{(t+|x-y|)^{m+\alpha}(t+|x-z|)^{m+\alpha}}
\end{equation}
whenever $|z-z'| < t/2$. The bilinear vertical square function is defined by setting
$$
BV(\nu_1,\nu_2)(x) = \Big( \int_0^{\infty} |\theta_t(\nu_1,\nu_2)(x)|^2 \,\frac{dt}{t} \Big)^{1/2}, \qquad x \in \R^n.
$$

For a Borel measure $\mu$ of order $m$ we set for
$f,g \in \bigcup_{p \in [1,\infty]} L^p(\mu)$ and $x \in \R^n$ that
$$
\theta_t^{\mu}(f,g)(x) = \theta_t(f\,d\mu, g\,d\mu)(x) = \int_{\R^n} \int_{\R^n} s_t(x,y,z)f(y)g(z) \,d\mu(y) \,d\mu(z)
$$
and
$$
BV_{\mu}(f,g)(x) = BV(f\,d\mu, g\,d\mu)(x) = \Big( \int_0^{\infty} |\theta_t^{\mu}(f,g)(x)|^2 \,\frac{dt}{t} \Big)^{1/2}.
$$
The above definitions make sense also when $\mu$ is finite. We use the notation
$BV_{\mu}^A(f,g)(x)$ or $BV^A(\nu_1,\nu_2)(x)$ to mean that the integration $\int_0^{\infty}$ is replaced with $\int_0^A$ for some constant $A$.

For simplicity, we only state a $T1$ theorem, instead of a $Tb$ theorem, in the square function setting.
\begin{theorem}\label{thm:mainSF}
Let $\mu$ be a measure of order $m$ on $\R^n$, and $C_0 < \infty$, $\delta_0 < 1$ and $l > 0$ be given constants.
Let $\beta > 0$ and $C_1$ be large enough (depending only on $n$). Suppose that for every
$(2,\beta)$-doubling cube $Q \subset \R^n$ with $C_1$-small boundary there exists $H_Q \subset \R^n$
such that $\mu(H_Q) \le \delta_0\mu(Q)$ and
\begin{equation}\label{eq:weakT1}
\sup_{\lambda > 0} \lambda^l \mu(\{x \in Q \setminus H_Q\colon\, BV_{\mu}^{\ell(Q)}(1_Q, 1_Q)(x) > \lambda\}) \le C_0\mu(Q).
\end{equation}
Then $BV_{\mu}\colon L^p(\mu) \times L^q(\mu) \to L^r(\mu)$ boundedly for all $1 < p, q < \infty$ and $1/2 < r < \infty$ satisfying  $1/r = 1/p + 1/q$.
\end{theorem}
To prove this theorem one can just prove a big piece $Tb$ (or $T1$) and then apply it in conjunction with the bilinear good lambda method.
That is, one can completely skip all the difficulties involving adapted Cotlar type inequalities and transferring testing conditions to maximal truncations -- this is simply because
no maximal truncations appear, which is related to the fact that 
$$
\Big( \int_A^{\infty} |\theta_t^{\mu}(f,g)(x)|^2 \,\frac{dt}{t} \Big)^{1/2}
\le \Big( \int_0^{\infty} |\theta_t^{\mu}(f,g)(x)|^2 \,\frac{dt}{t} \Big)^{1/2}
$$ 
for $A \ge 0$. This is one of the main things why the square function
setting is easier, and this aspect would be greatly amplified in the context of local $Tb$ theorems (which we do not consider in this paper).
In addition, it is also the case that it is much simpler to prove the corresponding big piece $T1$. There are multiple reasons for this, but they include at least the following:
\begin{itemize}
\item Suppression is much easier and does not involve the sophisticated Lipschitz suppression;
\item The paraproduct is simpler and the overall dyadic structure is simpler since one can use just one dyadic grid;
\item Probabilistic arguments are easier since good functions are not needed (only good Whitney regions);
\item The diagonal is trivial;
\item There are less symmetries since no adjoints appear.
\end{itemize}
If one desires to read this much more approachable argument, one can find the full details in the first (v1) arXiv version of the current paper.

\appendix

\backmatter

\bibliographystyle{amsalpha}

\printindex

\end{document}